\documentclass[13pt]{article}
\usepackage{amsfonts}

\usepackage{color}
\usepackage{pdfsync}

\newcommand{\beq}{\begin{equation}}
\newcommand{\enq}{\end{equation}}

\newtheorem{Theorem}{Theorem}[section]
\newtheorem{Lemma}[Theorem]{Lemma}
\newtheorem{Definition}[Theorem]{Definition}
\newtheorem{Remark}[Theorem]{Remark}

\newcommand{\beqa}{\begin{eqnarray}}
\newcommand{\edm}{\end{displaymath}}
\newcommand{\eeq}{\end{equation}}
\newcommand{\eeqa}{\end{eqnarray}}

\newcommand{\br}{\begin{Remark}}
\newcommand{\er}{\end{Remark}}

\newcommand{\bqa}{\begin{eqnarray}}
\newcommand{\eqa}{\end{eqnarray}}

\newcommand{\non}{\nonumber}

\RequirePackage{amsfonts}
\usepackage{stmaryrd}
\usepackage{amssymb}
\usepackage{mathrsfs}
\usepackage{amsmath}
\usepackage{subeqnarray}
\usepackage{verbatim}
\usepackage{cases}
\usepackage{amsfonts}
\usepackage{amssymb}
\usepackage{amssymb}
\usepackage{amssymb}
\usepackage{hyperref}
\begin{document}

\newpage
\pagenumbering{arabic} \setcounter{page}{1}

\begin{center}
	{\Large Pullback dynamics of 2D incompressible non-autonomous Navier-Stokes
		equation on Lipschitz-like domain}\\\vspace{0.25in}\ \ \  \  Xin-Guang Yang   \footnote[2]{College of Mathematics and Information Science, Henan Normal University, Xinxiang, 453007,
		P. R. China. Email: yangxinguang@hotmail.com.}\ \ \ \  Yuming Qin    \footnote[1]{Department of Applied Mathematics, Donghua University, Shanghai, 201620, P. R. China. Email: yuming\_qin@hotmail.com.}\ \ \ \
	To Fu Ma  \footnote[3]{Instituto de Ci\^encias Matem\'aticas e da Computa\c c\~ao, Universidade de S\~ao Paulo-Campus de S\~ao Carlos,
		Caixa Postal 668, 13560-970 S\~ao Carlos SP, Brazil. Email: matofu@icmc.usp.br.}\ \ \ \
	Yongjin Lu\footnote[4] {Department of Mathematics and Economics, Virginia State University, Petersburg, VA 23806, USA. Email: ylu@vsu.edu}
	\ \ \ \ \vspace{0.06in}
\end{center}
\vspace{0.06in}
\begin{abstract}
	This paper concerns the tempered pullback dynamics of 2D incompressible non-autonomous Navier-Stokes equation with non-homogeneous boundary condition on Lipschitz-like domain.
	With the presence of a time-dependent external force $f(t)$ which only needs to be pullback translation bounded, we establish the existence of a minimal pullback attractor with respect to a universe of tempered sets for the corresponding non-autonomous dynamical system. We then give estimate on the finite fractal dimension of the attractor based on trace formula. Under the additional assumption that the external force is the sum of a stationary force and a non-autonomous perturbation, we also prove the upper semi-continuity of the attractors as the non-autonomous perturbation vanishes. Lastly, we also investigate the regularity of these attractors when smoother initial data is given. Our results are new even in the case of smooth domains.
\end{abstract}

\section{Introduction}\label{sec1}
\setcounter{equation}{0}
\ \ \ \
The incompressible Navier-Stokes equation (NSE for short) is a well-known physical model in
 hydrodynamics and it plays a key role in understanding turbulence in science and engineer. One of the main motivations to study dynamical systems for
the Navier-Stokes equation is the attempt to gain further understanding of viscous fluid turbulence. Our objective in this paper is to study pullback dynamics of an incompressible 2D fluid flow in a Lipschitz-like domain. The fluid is subject to a prescribed tangential boundary velocity and a time-dependent external force. The motion of fluid is thus described by the following non-autonomous 2D Navier-Stokes equation with non-homogeneous boundary condition:

\begin{equation} \label{a1}
	\left\{
	\begin{array}{ll}
		\displaystyle \frac{\partial u}{\partial t} - \nu\Delta u+(u\cdot\nabla)u+\nabla p = f (t,x), \quad x \in \Omega, \; t \ge \tau, \smallskip \\
		{\rm div} \, u = 0, \quad x \in \Omega, \; t \ge \tau, \smallskip \\
		u = \varphi,   \quad \varphi \cdot n=0, \quad x \in \partial \Omega, \; t \ge \tau,  \smallskip \\
		u(\tau,x) = u_{\tau}(x), \quad x \in \Omega,
	\end{array}
	\right.
\end{equation}
where $\Omega$ is a bounded domain in $\mathbb{R}^2$ with Lipschitz boundary $\partial \Omega$,
and $\tau \in \mathbb{R}$ is an initial time.
Variables $u$ represents the fluid velocity field, $p$ denotes the pressure, and
$\nu$ is the kinematic viscosity. In addition, $n$ represents the exterior unit normal vector to $\partial \Omega$, $\varphi = \varphi(x)$ is a
prescribed tangential boundary velocity,
and $f(t,x)$ is a time-dependent forcing term.

Many results on the mathematical analysis for 2D incompressible Navier-Stokes equations are available in the literature. We first recall some results for {\it homogeneous Dirichlet or periodic boundary in smooth domain}.

 The  local and global existence of weak solutions for $N=2,3$-dimensional case were investigated by
 Leray \cite{leray1933,leray193401} and Hopf \cite{hopf1951}, which is called Leray-Hopf weak solution.  From 1960s, Ladyzhenskaya \cite{L2}, Temam \cite{te} and some other mathematicians have established the global existence of weak and strong solution with uniqueness for incompressible Navier-Stokes equation in space dimension two. More results about well-posedness in $N$-dimension, one can refer to \cite{cf1}, \cite{fmrt}, \cite{lion} and references therein.

  The study of infinite dimensional dynamical systems of autonomous partial differential equations was paid much attention in 1980s and has since been fast developed into a field that generates vast amount of results in the literature. One of the important topics in this area is the research of asymptotic behavior for {\it autonomous} Navier-Stokes equations in {\it smooth domain}. If the domain $\Omega$ is {\it bounded with smooth boundary}, the existence of compact global attractor was established by Babin and Vishik \cite{bv92}, Constantin, Foias and Temam \cite{cft}, Foias, Manley,Ilylin \cite{ily}, Rosa and Temam \cite{fmrt}, Ladyzhenskaya \cite{lady}, Robinson \cite{ro}, \cite{ro11}, Sell and You \cite{sel}, Temam \cite{te88} and so on. If the external force satisfies appropriate assumption, the fractal and Hausdorff dimension of global attractor has finite bound depending on Grashof number $G$ (see \cite{cf1}, \cite{cfmt},
 \cite{cft}, \cite{dg}, \cite{fmrt}, \cite{liu}, \cite{ro}, \cite{ro11} and references therein).
 {\it However, the inertial manifold for 2D Navier-Stokes equation is still open}.

 Three conditions, invariance, attracting and
 compact, are necessary for a set to become an attractor. However, if the external force is time-independent, invariance of an attractor is lost. In this case, if the external force is translation compact, then it leads to the existence of uniform attractor. Applying the existence theory to {\it non-autonomous} 2D NSE, the uniform attractors and its finite fractal dimension were derived by Chepyzhov and Vishik \cite{cv3}, Lu \cite{lu} for {\it homogeneous Dirichlet or periodic boundary condition}, and by Miranville and Wang \cite{mw97} for inhomogeneous boundary condition.
 Unfortunately, the forward invariance and uniform estimates in some cases are not easy to achieve, the idea is to use pullback invariance to instead, which leads to the pullback attractor theory. The pullback dynamics for same problem of 2D {\it non-autonomous} NSE
 are obtained in \cite{clr2006}, \cite{clr}, \cite{ci}, \cite{gmr}, \cite{llr} and references therein.

{\it  Another situation to consider is the non-regular domain.} If the domain $\Omega$ is {\it an unbounded smooth manifold} where the Poincar\'e's inequality holds, the existence of global and pullback attractor were obtained from energy equation method by Rosa \cite{ros}, Langa, {\L}ukaszewicz and Real \cite{llr}. If the domain is {\it channel-like or stripe-like}, the existence of solutions and its long dynamics for 2D NSE were obtained by Babin and Vishik \cite{bv}, Temam and Wang \cite{tw}, Zelik \cite{ze}.

As for the 2D {\it autonomous} Navier-Stokes equations with {\it non-homogenous boundary} $u|_{\partial \Omega} = \varphi$, {\it on smooth domain},
the existence of global attractors with finite fractal dimension was established by Miranville and Wang \cite{mw96,mw97}. Their arguments
are based on the construction of a {\it background flow} $\psi$ satisfying
\begin{equation} \label{bgf}
	{\rm div} \, \psi =0 \mbox{ in } \Omega \quad {\rm and}  \quad \psi = \varphi \mbox{ on } \partial \Omega,
\end{equation}
which suffices for estimating fractal dimension. Early constructions of such background flows were presented in \cite{mir,tw}.

This {\it non-homogeneous} problem, with $f=f(x)$, was extended to {\it Lipschitz domains (non-smooth)} by Brown, Perry and Shen \cite{bps},
who constructed an appropriate background flow (satisfying (\ref{bgf})),
by taking into account the possible non-smooth character of the boundary $\partial \Omega$.
Then, writing $v=u-\psi$,
problem (\ref{a1}) becomes
\begin{equation} \label{a1-equiv}
	\left\{
	\begin{array}{ll}
		\displaystyle \frac{\partial v}{\partial t} - \nu\Delta v + (v \cdot \nabla)v  + (v \cdot \nabla) \psi + (\psi \cdot \nabla)v + \nabla p  \smallskip \\
		\qquad\qquad\qquad\qquad\qquad\qquad\qquad\qquad = f + \nu \Delta \psi -   (\psi \cdot \nabla)\psi, \quad x \in \Omega, \; t \ge \tau, \smallskip \\
		{\rm div} \, v = 0, \quad x \in \Omega, \; t \ge \tau, \smallskip \\
		v = 0,  \quad x \in \partial \Omega ,     \smallskip \\
		v(\tau,x) = v_{\tau}(x), \quad x \in \Omega,
	\end{array}
	\right.
\end{equation}
which is now a homogeneous problem.
We note that if $v$ is a solution of problem (\ref{a1-equiv}) with initial data
$v_{\tau} = u_{\tau} - \psi$, then $u=\psi + v$ is a
solution of problem (\ref{a1}) with initial data $u_{\tau}$.
The Hadamard well-posedness of problem based on stationary system (see \cite{fkv,shen}) and
the existence of a global attractor of finite fractal dimension
of problem (\ref{a1-equiv}) were proved in \cite{bps}, with respect to weak solutions and $f=f(x)$.

For system
 \eqref{a1-equiv} with {\it non-autonomous} external force satisfying {\it translation compact condition}, Wu and Zhong \cite{wz06} proved the existence of uniform attractor in space $H$.
In \cite{ze15}, the translation compactness condition ensures the forward dynamics to satisfy some dissipative evolutionary equations.

In this work, we consider pullback asymptotic dynamical behavior of problem \eqref{a1} with time-dependent force $f= f(t,x)$ which only satisfies {\it pullback translation bounded} condition. The major contributions of this paper are that:
	
	(1). The boundary condition is non-homogeneous. Therefore, in order to convert the problem into one that has homogeneous boundary condition, as mentioned above, we have to introduce a background function $\psi$, the solution to a non-homogenous Stokes problem. Estimates of $\psi$ on a Lipschitz-like domain established in \cite{bps} have to be critically invoked.
	
	(2). The underlying domain is Lipschitz-like and thus non-smooth. The lack of smoothness of the domain requires more delicate estimates that depend on the geometric properties of $\Omega$.
	
	 (3). The external force is assumed to be only pullback translation bounded, a major departure from the usual pullback translation compact condition (also as uniformly pullback tempered) assumed in the literature. The weaker compactness of the external force makes the proof more technically challenging.
	
	  (4). We give a rather complete study of pullback dynamics of \eqref{a1}: we establish not only the existence of a minimal family of pullback attractors, but also estimates on the finite fractal dimension of these attractors, upper-semicontinuity of these attractors (by assuming that the external force follows a special form), and regularity of these attractors. As the foundation of all other results, the first result, the existence of a minimal family of pullback attractors, has an important feature that the tempered universe $\mathcal{D}_{\mu}$ ($0<\mu\leq \mu_0=\frac{\nu\lambda_1}{2}$) where we seek the minimal pullback attractor is \emph{not necessarily bounded} and it contains the universe of fixed bounded sets in $H$. The minimality property then implies that uniqueness of the attractor. Proof of these different aspects of the pullback dynamics demand different techniques:

 (a) For example, to establish existence of a minimal family of pullback attractors, the asymptotic compactness, the technically most intensive component of the proof, is obtained from a decomposition method of the phase space. First introduced by Ma, Wang and Zhong in \cite{wz}, to show that a process is pullback $\mathcal{D}$-asymptotically compact (see Definition \ref{de8.28} in Section \ref{ap} below), one only needs to check the so-called pullback-$\mathcal{D}$-condition (MWZ): the existence of a fixed pullback time, up to which, the projection of solution to the finite dimensional subspace of the phase space is bounded, while the orthogonal complemented projection to the infinite dimensional subspace is sufficiently small in appropriate topology.
	
	   (b) To establish estimates on the fractal dimension of the dynamics, we apply the trace formula to describe the box counting which is based on a quasi-differential compact approximated operator.  By using this technique, our estimate for the fractal dimension of the pullback attractors has a $1/2$ lower power on the generalized Reynold number than that for the forward global attractor obtained in \cite{bps}.
	
	    (5). We assume the external force to have a form of $f(t,x) = \sigma_0(x)+\varepsilon \sigma(t,x)$ rather than $\varepsilon g(t,x)$ when studying the upper semi-continuity property of the pullback attractors. The former form of the external force allows us to study the continuity of the pullback attractors when the external force approaches to a possibly non-zero stationary state while the latter form only admits a zero final state of the external force.

(6). As far as we know, our results, the existence of minimal pullback attractors for 2D NSE, are new even in the case where $\Omega$ has smooth boundary $\partial \Omega$.

For the sake of completeness, we will also prove the well-posedness result for \eqref{a1}. The rest of the article will be organized as follows: in Section 2, we first review some preliminaries on studying NSE and the background flow \eqref{bgf} and then summarize the main results of this article: existence of global weak solution to \eqref{a1}, existence of minimal and unique family of pullback $\mathcal{D}_{\mu}$-attractors in $H$, estimates on finite fractal dimension of pullback attractors, continuity of pullback attractors and regularity of pullback attractors. In Section 3, we will prove these results in order.

\section{Main results and comments}\label{mr}
\setcounter{equation}{0} In this section, we shall state our main results, i.e., the existence of weak solution, tempered pullback dynamics of process for system \eqref{a1}.

\subsection{Preliminaries}\label{sec2}
\setcounter{equation}{0}
\ \ \ \  Denote $E:=\{u|u\in(C^{\infty}_0(\Omega))^2, div u=0\}$, $H$ is the
closure of $E$ in $(L^2(\Omega))^2$, $(\cdot,\cdot)$ and $|\cdot|$ denote the inner product and
norm in $H$ respectively, i.e., \begin{eqnarray}
	&&(u,v)=\sum^2_{j=1}\int_{\Omega}u_j(x)v_j(x)dx, \ \forall\ u, v\in
	(L^2(\Omega))^2\ \ \mbox{and}\ \  |u|^2=(u,u).\nonumber\end{eqnarray}
$V$ is the
closure of $E$ in $(H^1(\Omega))^2$ topology, $((\cdot,\cdot))$ and $\|\cdot\|$ denote
the inner product and norm in $V$ respectively, i.e., \begin{eqnarray}&&((u,v))=\sum^2_{i,
		j=1}\int_{\Omega}\frac{\partial u_j}{\partial
		x_i}\frac{\partial v_j}{\partial x_i}dx,\ \forall\ u, v\in
	(H^1_0(\Omega))^2\ \ \mbox{and}\ \ \|u\|^2=((u,u)).\nonumber\end{eqnarray}
$H'$ and $V'$ are
dual spaces of $H$ and $V$ respectively, where the injections $V\hookrightarrow H\equiv H'\hookrightarrow V'$ are
dense and continuous. The norm $\|\cdot\|_{*}$ and $\langle\cdot\rangle$ denote the
norm in $V'$ and the dual product between $V$ and $V'$ respectively.

Let $P$ be the Helmholz-Leray orthogonal projection operator  from $(L^2(\Omega))^3$ onto $H$. We define $A:=-P\Delta $ to be the Stokes operator with domain $D(A)=(H^{2}(\Omega))^3\bigcap
(H^{1}_{0}(\Omega))^3$, then the operator $A: V\rightarrow V'$ has the property $\langle Au,v\rangle=((u,v))$ for all $u, \ v\in V$ which is an isomorphism from $V$ into $V'$.
$\{\lambda_j\}^{\infty}_{j=1}\ (0<\lambda_1\leq\lambda_2\leq\cdots)$
are eigenvalues of operator $A$ for the eigenvalue problem
$A u=\lambda u,\ u|_{\partial\Omega}=0$ in $L^2(\Omega)$. $\{\omega_j\}^{\infty}_{j=1}$ is an orthonormal
basis of operator $A$ corresponding to $\{\lambda_j\}^{\infty}_{j=1}$,
 i.e., $A\omega_j=\lambda_j\omega_j$.

We define the bilinear and trilinear operators as (see \cite{te88})
\begin{eqnarray}
	&&B(u,v):=P((u\cdot\nabla)v),\ \ \forall \ u,v \in E,\label{5-1}\\
	&&b(u,v,w)=(B(u,v),w)=\sum^2_{i,
		j=1}\int_{\Omega}u_i\frac{\partial v_j}{\partial x_i}
	w_jdx\label{5}\end{eqnarray}
respectively, where $B(u,v)$ is a linear continuous operator and $b(u,v,w)$ satisfies
\begin{equation}\left\{
	\begin{array}{llll}
		b(u,v,v)&=& 0,\ &\forall\ u\in V,\ \ v\in (H^1_0(\Omega))^2,\\
		b(u,v,w)&=& -b(u,w,v),\ &\forall\ u,v,w\in V,\\
		|b(u,v,w)|&\leq&
		 C|u|^{\frac{1}{2}}\|u\|^{\frac{1}{2}}\|v\||w|^{\frac{1}{2}}\|w\|^{\frac{1}{2}},\ &\forall\ u, v, w\in V.\end{array}\right.\label{q7}\end{equation}

The fractal operator $A^{s}$ $(s\in \mathbb{C})$ is defined as (see \cite{te})
\begin{align}
&A^{s}f=\sum_{j}\lambda_{j}^{s}(f,\omega_j)\omega_{j},\ s\in \mathbb{C},\ j\in \mathbb{R}, \label{2}\\
&V^s=D(A^{s})=\Big\{g\in H: A^{s}g\in H,\ \sum_{i=1}^{+\infty}\lambda_i^{2\alpha}|(u,\omega_i)|^2<+\infty\Big\},\label{3}\\
&\|A^{\sigma}u\|=\Big(\sum^{+\infty}_{i=1}\lambda_i^{2\sigma}|(u,\omega_i)|^2\Big)^{1/2}.\label{3-1}
\end{align}
$D(A^{s})$ denotes the domain of $A^{s}$ with the inner product and the norm $\|\cdot\|_{s}$ as
\begin{equation}
(u,v)_{V^s}=(A^{\frac{s}{2}}u, A^{\frac{s}{2}}v),\quad\ \|u\|_{V^s}^2=(u,u)_{V_s}
\end{equation}
and $V^{\frac{s+1}{2}}=D(A^{\frac{s+1}{2}})$ with the norm $\|\cdot\|_{V^{\frac{s+1}{2}}}$. Especially, $V=V^1$, $V^{2}=W:=(H^{2}(\Omega))^3\bigcap (H^{1}_{0}(\Omega))^3$.
Moreover, $A^{s}$ satisfies (see \cite{bps})
\begin{eqnarray}
	&&\|u\|_{L^4}\leq C_{1}|A^{\frac{1}{4}}u|,\ \forall u\in D(A^{\frac{1}{4}}).\label{4}\end{eqnarray}
The Gagliardo-Nirenberg interpolation inequality: \begin{eqnarray} &&|A^{1/2}u|^2\leq
	C_{2}|A^{1/4}u||A^{3/4}u|,\ \ \forall \ u\in D(A^{3/4}).\label{6}\end{eqnarray}
Hardy's inequality:\begin{eqnarray}
	\int_{\Omega}\frac{|u(x)|^{2}}{[dist(x,\partial\Omega)]^{2}}dx\leq C_{3}\int_{\Omega}|\nabla u(x)|^{2}dx,\ \forall u \in V.\label{7}\end{eqnarray}

\subsection{Background functions for the Stokes problem on Lipschitz domains}
\label{sec3}
\ \ \ \ The main idea we use in this paper is to transform \eqref{a1} into a problem that has homogeneous boundary condition. To this end, we introduce a new function $v=u-\psi$, where $u$ is the solution to the original problem \eqref{a1} with non-homogeneous boundary condition $\varphi$ and the background function $\psi$ is the solution to the following system that shares the same boundary condition $\varphi$ as \eqref{a1}:

\begin{eqnarray}
	\left\{\begin{array}{ll} \mbox{div}\ \psi=0,\ \ \ \ in\ \ \Omega,\\
		\psi=\varphi,\ \ \ \ \varphi\cdot n=0\ \ on\ \ \partial\Omega,\end{array}\right.\label{s3}
\end{eqnarray}
The idea is motivated by Miranville and Wang \cite{mw96} and \cite{mw97}. Brown, Perry and Shen \cite{bps} extended the well-posedness of global solution of the non-homogeneous boundary problem to non-smooth domains (such as Lipschitz-like domain) of 2D Navier-Stokes equation by critically invoking estimates of the Stokes problem
\begin{equation}\left\{
	\begin{array}{lll}
		-\triangle \hat{u}+\nabla q= 0,\ \ in\ \ \ \Omega,\\
		 \mbox{div}\ \hat{u}=0,\ \ in\ \ \ \Omega,\\
		\hat{u}=\varphi,\ \ a.e.\ on\ \partial\Omega $\ $in$\ $ the$\ $ sense$\ $ of$\ $ nontangential$\ $ convergence.\end{array}\right.\label{8}\end{equation}

\begin{Definition}\label{de3.1}
	(The Lipschitz-like domain, see \cite{bps}, \cite{fkv}) We call a bounded set $\Omega\subset \mathbb{R}^d$ to be a Lipschitz-like domain if its boundary $\partial\Omega$ can be covered by finitely many balls $B_i=B(Q_i,r_0)$ centered at the point $Q_i\in\partial\Omega$ such that for each ball $B_i$, there exists a rectangular coordinate system and a Lipschitz function $\Psi: \mathbb{R}^{d-1}\rightarrow \mathbb{R}$ with $$B(Q_i,3r_0)\cap \Omega=\{(x_1,x_2,\cdots,x_d)|x_d>\Psi_i(x_1,x_2,\cdots, x_{d-1})\}\cap\Omega.$$
An example of a Lipschitz-like domain is the square $\Omega=(0,1)\times(0,1)$.	
\end{Definition}

\begin{Lemma} \label{le3.4}
	(See e.g., \cite{bps}) (1) The background function $\psi$ which is a solution of Stokes problem satisfies the following estimates in a Lipschitz-like domain
	\begin{eqnarray} &&\sup_{x\in\Omega}|\psi(x)|+\sup_{x\in\Omega}|\nabla\psi(x)|dist(x,\partial\Omega)\leq C_{4}\|\varphi\|_{L^{\infty}(\partial\Omega)},\label{c5}\\
		 &&\||\nabla\psi|dist(\cdot,\partial\Omega)^{1-\frac{1}{p}}\|_{L^{p}(\Omega)}\leq C_{5}\|\varphi\|_{L^{p}(\partial\Omega)},\ 2\leq p\leq\infty,\label{d6}\\
		&&\|\psi\|_{L^{\infty}(\Omega)}\leq C_{4}\|\varphi\|_{L^{\infty}(\partial\Omega)}.\label{b1}\end{eqnarray}
	
	(2) If $\psi$ satisfies the problem \begin{equation}\left\{
		\begin{array}{lll}
			 \mbox{div}\ \psi =0,\ x\in\Omega,\\
			\psi=\hat{u},\ if\ x\in\{x\in \Omega;\ dist(x,\partial\Omega)<C'_{1}\varepsilon\},\\
			\psi=\varphi,\ on\ \partial\Omega\ in\ sense\ of\ nontangential\ convergence.\end{array}\right.\label{d6-1}\end{equation}
	Then we have
	\begin{eqnarray} &&\Delta\psi=\nabla(q\eta_{\varepsilon})+F,\label{d7}\end{eqnarray}
	where
	\begin{eqnarray} &&Supp ~\psi\subset\{x\in \bar{\Omega};\ dist(x,\partial\Omega)<C'_{2}\varepsilon\},\label{d8}\\
		&&\|F\|_{L^{2}(\Omega)}\leq \frac{C}{\varepsilon^{\frac{3}{2}}}\|\varphi\|_{L^{2}(\partial\Omega)},\ \nabla q=\triangle \hat{u},\label{d9}\\
		&&F=0,\ \ \ x\in \{x|dist(x,\partial\Omega)<C'_{1}\varepsilon\ or\  dist(x,\partial\Omega)>C'_{2}\varepsilon\},\label{d10}\end{eqnarray}
	here $\varepsilon$ is an arbitrary fixed positive constant.
\end{Lemma}

\subsection{Main result I: Existence of global weak solution}\label{sub4.1}
\ \ \ \ Let $\psi$ be a background flow function which satisfies the Stokes problem \eqref{8} and $v(t,x)=u(t,x)-\psi(x)$, since the Stokes operator is compact, then \eqref{a1-equiv} is equivalent to the following homogeneous boundary value problem
\begin{equation}\left\{
\begin{array}{lll}
v_t-\nu\Delta v+(v(t) \cdot \nabla)v+(v(t) \cdot \nabla)\psi+(\psi \cdot \nabla)v+\nabla (p-\nu q\eta_{\varepsilon})=
\bar{f}-(\psi\cdot\nabla)\psi,\ &\\
\mbox{div}\ v=0,\ &\\
v|_{\partial\Omega}=0,\ &\\
v(\tau,x)=v_{\tau}(x)=u_{\tau}(x)+\psi(x),\ & \end{array}\right.\label{c1}\end{equation}
where
$\bar{f}=f(x,t)+\nu F(x)$.

Let $v_{\tau}\in H\ \mbox{or} \  V$, applying the Leray projector $P$ to problem \eqref{c1}, using the divergence free condition, we derive the following equivalent abstract weak form
\begin{equation}\left\{
\begin{array}{lll}
v_t+\nu
Av+B(v,v)+B(v,\psi)+B(\psi,v)=P\bar{f}-B(\psi,\psi),\ &\\
\mbox{div}\  v=0\ & \end{array}\right.\label{c1a}\end{equation}
and its global weak solution can be defined as the following.

\begin{Definition}\label{de4.2} Let $\Omega$ be a Lipschitz-like domain, $u_{\tau}\in H,\ f(x,t)\in L^2_{loc}(\mathbb{R};V'),\ \varphi\in L^{\infty}(\partial\Omega)$ and $\varphi\cdot n=0$ on $\partial\Omega$, u is called a weak solution of problem \eqref{a1}
	provided that\\
	(i) $u\in C([\tau,T];H),u(\cdot,\tau)=u_{\tau}$,and $du/dt\in L^{2}(\tau,T;V')$;\\
	(ii) for all $v\in C^{\infty}_{0}(\Omega)$ with $ \mbox{div}\ v=0$, we get
	\begin{eqnarray}&&\frac{d}{dt}<u,v>-\nu<u,\Delta v>-\int_{\Omega}\sum^{2}_{i,j=1}u^{i}u^{j}\frac{\partial v^{i}}{x_{j}}dx=<f,v>\nonumber\end{eqnarray} in the distribution sense on $[\tau,T]$;\\
	(iii) there exists functions $\psi\in C^{2}(\Omega)\cap L^{\infty}(\partial\Omega), q\in C^{1}(\Omega)$ and $g\in L^{2}(\Omega)$ such that
	\begin{equation}\left\{
	\begin{array}{lll}
	\triangle \psi=\nabla q+\hat{g},\ &in& \Omega,\nonumber\\
	\mbox{div}\ \psi=0,\ &in& \Omega,\nonumber\\
	\psi=\varphi\ &on&\partial\Omega,\nonumber\end{array}\right.\end{equation}
	where $\psi$ can reach its boundary values in the sense of non-tangential convergence and $u-\psi\in L^{2}(\tau,T;V)$.
\end{Definition}

The existence of global weak solutions of problem \eqref{c1a} (which is equivalent to \eqref{c1}) can be stated as the following theorem.
\begin{Theorem}\label{th4.1}
	Let $v_{\tau}(x)\in H$, $f(x,t)\in L^2_{loc}(\mathbb{R};V')$, then there exists a
	unique solution $v(t,x)$ of non-autonomous problem \eqref{c1a} satisfying  $$v(t,x)\in L^{\infty}(\tau,T;H)\bigcap L^{2}(\tau,T;V),$$ and $\frac{dv}{dt}$ is uniformly  bounded in  $ L^{2}(\tau,T;V')$.
\end{Theorem}
{\bf Proof.} Applying the Galerkin approximated method and compact argument as in \cite{bps} to our non-autonomous problem, this result could be established. We give the outline of the proof in Section \ref{sece}.$\hfill$$\Box$\\

The continuous dependence of global solution is stated below:

\begin{Theorem}\label{th4.3}
	Let $u_{\tau}\in H,\ f\in L^2_{loc}(\mathbb{R};V'), \varphi\in L^{\infty}(\partial\Omega)$, and $\varphi\cdot n=0$ on $\partial\Omega$. Then, we obtain that\\
	(1) based on the existence of weak solutions, it follows $u(t,x)\in L^4([\tau,T]\times\Omega)$;\\
	(2) the problem \eqref{a1} possesses a
	unique weak solution $u(t,x)\in L^{\infty}(\tau,T;H)\cap L^2(\tau,T;V)$ which is continuously dependent on the initial data.
\end{Theorem}
{\bf Proof.} See Section \ref{sece}.$\hfill$$\Box$

\begin{Remark}\label{re4.4}
	For the original problem \eqref{a1} and the initial data $v_0=u_0-\psi$ of the equation \eqref{c1a}  with initial data $u_0$ and boundary data $\varphi$, based on the wellposedness of problems \eqref{c1a}, we have the generation of process $\{U(t,\tau)\} (t\geq \tau):H\rightarrow H$ which implies $\psi+U(t,\tau)(u_{\tau}-\psi)=\psi+U(t,\tau)v_{\tau}$ is the solution to problem \eqref{a1}.
\end{Remark}

Next, we will define the regularity of solution for problem \eqref{a1}.
\begin{Definition} \label{de7.5} Let $\Omega$ be a Lipschitz-like domain, $u_{\tau}\in  D(A^{\frac{\sigma}{2}}),\ f(x,t)\in L^2_{loc}(\mathbb{R};H),\ \varphi\in L^{\infty}(\partial\Omega)$, u is called a regular weak solution of the problem \eqref{a1}
	provided\\
	(i) $u\in C_w([\tau,T];D(A^{\frac{\sigma}{2}})),u(\cdot,\tau)=u_{\tau}$,and $du/dt\in L^{2}(\tau,T;V')$,\\
	(ii) for all $v\in C^{\infty}_{0}(\Omega)$ with $ \mbox{div}\ v=0$, we get
	\begin{eqnarray}&&\frac{d}{dt}<u,A^{\sigma}v>-\nu<u,A^{\sigma+1} v>+\int_{\Omega}\sum^{2}_{i,j=1}u^{i}u^{j}\frac{\partial (A^{\sigma}v^{i})}{x_{j}}dx=<f,A^{\sigma}v>\nonumber\end{eqnarray} in the distributed sense,\\
	(iii) there exists functions $\psi\in C^{2}(\Omega)\cap L^{\infty}(\partial\Omega), q\in C^{1}(\Omega)$ and $g\in L^{2}(\Omega)$ such that
	\begin{equation}\left\{
	\begin{array}{lll}
	\triangle \psi=\nabla q+\hat{g},\ &in& \Omega,\nonumber\\
	 \mbox{div}\ \psi=0,\ &in& \Omega,\nonumber\\
	\psi=\varphi\ &on&\partial\Omega,\nonumber\end{array}\right.\end{equation}
	where $\psi$ can reach its boundary values in non-tangential sense and $u-\psi\in L^{2}(\tau,T;D(A^{\frac{\sigma+1}{2}}))$.
\end{Definition}

By similar technique in Theorem \ref{th4.1} and uniform estimates in $D(A^{\frac{\sigma}{2}})$ and $V$, the regularity result of solution can be established.
\begin{Theorem}\label{th7.6}
	(1) Assume $f\in L^2_{loc}(\mathbb{R};H)$, let $v(t,x)$ be the solution of \eqref{c1a} with the initial data $v_{\tau}\in  D(A^{\frac{\sigma}{2}})$, then the regular solution $v(t,x)\in L^2(\tau,T;D(A^{\frac{\sigma}{2}}))\cap L^2(\tau,T;D(A^{\frac{\sigma+1}{2}}))$ exists for $\sigma\in[0,\frac{1}{2}]$.

(2) Let $u_{\tau}\in D(A^{\frac{\sigma}{2}}),\ f(x,t)\in L^2_{loc}(\mathbb{R};H),\ \varphi\in L^{\infty}(\partial\Omega)$. Then, the problem \eqref{a1} possesses regular solution $u(t,x)\in C_w([\tau,T];D(A^{\frac{\sigma}{2}}))\cap L^2(\tau,T;D(A^{\frac{\sigma+1}{2}}))$ which is norm-to-weak continuous.
	
	(3) Assume $f\in L^2_{loc}(\mathbb{R};H)$, let $v(t,x)$ be the solution of problem \eqref{c1a} with the initial data $v_{\tau}\in D(A^{\frac{\sigma}{2}})$, then we obtain that the
	regular solution $v(t,x)$ is norm-to-weak continuous, i.e., $v(t,x)\in C_w([\tau,T];D(A^{\frac{\sigma}{2}}))$.
\end{Theorem}
{\bf Proof.} See Section \ref{sec7}.$\hfill$$\Box$

\begin{Remark}
For the autonomous problem, Brown, Perry and Shen \cite{bps} has presented the existence of global unique weak solution in $H$ which can be derived from our result when $f(t)$ reduces to time-independent.
\end{Remark}
\subsection{Main Result II: Existence of minimal and unique family of pullback $\mathcal{D}_{\mu}$-attractors in $H$}\label{sub5.4}
\ \ \ \ In this part, we give a new definition and condition of external force-pullback translation boundedness, and then present the minimal and unique family of pullback attractors.

	\begin{Definition}
		We call the external force $f(t,x)\in L^2_{loc}(\mathbb{R};X)$ is pullback translation bounded if it satisfies
		\begin{eqnarray}\label{13-a}
		\sup_{s\leq \tau,\ h\leq 0}\int^s_{s+h}\|f(r,x)\|^2_{X}dr<+\infty
		\ \ \mbox{or}\ \
		\sup_{s\leq \tau}\int^{s}_{s-1}\|f(r,x)\|^2_{X}dr<+\infty
		\end{eqnarray}
		for the cases: (a-1) some $t \in \mathbb{R}$ and some $h\leq 0$; (a-2) some  $t \in \mathbb{R}$ and all $h\leq 0$; (a-3)  all $t \in \mathbb{R}$ and all $h\leq 0$. In fact, (a-1), (a-2) and (a-3) are equivalent, so we only write \eqref{13-a} with all $t \in \mathbb{R}$ and all $h\leq 0$ for convenient.
Moreover, since $f\in L^2_{loc}(\mathbb{R};X)$, if we use $s\leq t$, the pullback translation bounded also well defined. We denotes all the pullback translation bounded functions as $L^2_{pb}(\mathbb{R};X)$.
	\end{Definition}
	
	\begin{Definition}
			The function $f(t,x)\in L^2_{loc}(\mathbb{R};X)$ is called pullback tempered, if
			\begin{eqnarray}
			\sup_{s\leq \tau}\int^s_{-\infty}e^{\mu(r-s)}\|f(r,x)\|^2_{X}dr<+\infty
			\ \ \mbox{or}\ \
				\lim_{s\rightarrow -\infty}\int^s_{-\infty}e^{\mu(r-s)}\|f(r,x)\|^2_{X}dr<+\infty
				\end{eqnarray}			
			holds for all $t\in\mathbb{R}$. In addition, for $s\leq t$, this is also well-defined. We denotes all the pullback tempered functions as $L^2_{pt}(\mathbb{R};X)$.
		\end{Definition}

\begin{Definition}
	We say $f(t,x)$ is uniformly pullback tempered if
		\begin{eqnarray}
		\int^t_{-\infty}e^{\mu s}\|f(s)\|^2_{X}d s<+\infty\ \
		\mbox{or}\ \
		\int^{\tau}_{-\infty}e^{\mu s}\|f(s)\|^2_{X}d s<+\infty\label{gg1}
		\end{eqnarray}
		thanks to $f\in L^2_{loc}(\mathbb{R};X)$.
		We denotes all the uniform pullback tempered functions as $L^2_{upt}(\mathbb{R};X)$.
\end{Definition}

	\begin{Remark}
		The pullback translation boundedness is equivalent to pullback tempered condition, which is weaker than uniformly pullback tempered condition. The proof of equivalence and more details can be found in Appendix.
	\end{Remark}

	Let $\mathcal{P}(H)$ be the collection of all nonempty subsets in $H$, $\hat{D}=\{D(t)\}\subset \mathcal{P}(H)$ be a subset in $\mathcal{P}(H)$ which is not necessarily bounded, $B(0,\rho_{\hat{D}}(t))$ is a family of balls with center $0$ and radius $\rho_{\hat{D}}(t)$, where $\rho_{\hat{D}}(t)$ satisfies that
\begin{equation}
	e^{-\mu t}\Big(\frac{4}{\nu\mu}\Big[\frac{C\nu^2}{\varepsilon}|\partial\Omega|\|\varphi\|^2_{L^{\infty}(\partial\Omega)}+C\varepsilon\ |\partial\Omega|\|\varphi\|_{L^{\infty}(\partial\Omega)}^2\Big]+|\rho_{\hat{D}}(\tau)|^2\Big)e^{\mu\tau}\rightarrow 0\label{gg6}
\end{equation}
as $\tau\rightarrow -\infty$, which means that there exists a pullback time $\tau_1\leq t$, such that for any $\tau<\tau_1$,
\begin{equation}
	|\rho_{\hat{D}}(\tau)|^2e^{\mu\tau}\rightarrow 0.\label{gg7}
\end{equation}

We define the universe $\mathcal{D}_{\mu}=\{\hat{D}=\{D(t)\}\}$ as
\begin{equation}
	\mathcal{D}_{\mu}=\{\hat{D}|D(t)\subset B(0,\rho_{\hat{D}}(t))\ with\ \rho_{\hat{D}}\ satisfying\ \eqref{gg7}\}.
\end{equation}

Let $\hat{B}_0=\{B_0(t)\}_{t\in\mathbb{R}}$ be a family of balls, where $B_0(t)=\overline{B(0,\rho_0(t))}$ is a ball with center $0$ and radius $\rho_0(t)$, and $\rho_0(t)$ is defined as
\begin{eqnarray}
&&\rho^2_0(t)=1+\frac{4}{\nu\mu}\Big[\frac{C\nu^2}{\varepsilon}|\partial\Omega|\|\varphi\|^2_{L^{\infty}(\partial\Omega)}+C\varepsilon\ |\partial\Omega|\|\varphi\|_{L^{\infty}(\partial\Omega)}^2\Big]\nonumber\\
&&\hspace{2cm}+\|f\|^2_{L^2_{loc}(\mathbb{R};V')}+\frac{4e^{-\mu n_0h}}{\nu\lambda_1(1-e^{-\mu n_0h})}\|f\|^2_{L^2_{pb}(\mathbb{R};V')},
\end{eqnarray}
here $0<\mu\leq \mu_0=\frac{\nu\lambda_1}{2}$, $n_0>0$ is a fixed integer and $h>0$ is a constant.
Using the tempered pullback attractor theory in \cite{gmr} and pullback translation bounded external force, the minimal $\mathcal{D}_{\mu}$-family of pullback attractors in $H$ can be stated as follows.
	\begin{Theorem}\label{th5.6}
		Assume $f \in  L^2_{loc}(\mathbb{R};V')$ is pullback translation bounded ($f\in L^2_{pb}(\mathbb{R};V')$) (or uniformly pullback tempered in $ L^2_{upt}(\mathbb{R};V')$),  let $v_{\tau}\in H$, then the norm-to-weak continuous process $\{U(t,\tau)\}$ possesses a minimal $\mathcal{D}_{\mu}$-family of pullback attractors $\mathcal{A}^H_{\mathcal{D}_{\mu}}=\{\mathcal{A}^H_{\mathcal{D}_{\mu}}(t)\}$ in $H$ for the system \eqref{c1a} which is equivalent to problem \eqref{a1}.
	\end{Theorem}
	{\bf Proof.} By the definition of universes to achieve pullback absorbing set which is no need to be bounded, combining Condition-(MWZ) to achieve asymptotic compactness, the theorem can be proved, see Section \ref{sec5}. $\hfill$$\Box$\\

	Next, if we can show that the universe $\mathcal{D}_{\mu}$ is inclusion closed and $\hat{B}_0\in\mathcal{D}_{\mu}$, the family of pullback $\mathcal{D}_{\mu}$-attractors is unique by Theorem \ref{th5.7}.
	\begin{Theorem}\label{th5.7}
		Assume $f \in  L^2_{loc}(\mathbb{R};V')$ is pullback translation bounded ($f\in L^2_{pb}(\mathbb{R};V')$) (or uniformly pullback tempered in $ L^2_{upt}(\mathbb{R};V')$ with $0<\mu<\mu_0$), let $v_{\tau}\in H$, then the norm-to-weak continuous process $\{U(t,\tau)\}$ possesses a unique $\mathcal{D}_{\mu}$-family of pullback attractors $\mathcal{A}'_{\mathcal{D}_{\mu}}=\{\mathcal{A}'_{\mathcal{D}_{\mu}}(t)\}$ in $H$ for the system \eqref{c1a} which is equivalent to \eqref{a1}.
	\end{Theorem}
	{\bf Proof.} See Section \ref{sec5}.$\hfill$$\Box$

\begin{Remark}
Denoting $\mathcal{D}^H_F$ as the universe of fixed nonempty bounded subset of $H$ with same property as \eqref{gg6}. We see that $\mathcal{D}^H_F\subset \mathcal{D}^H_{\mu}\subseteq \mathcal{D}^H_{\mu_0}$
			which are inclusion closed. The minimal families of pullback attractors
			$\mathcal{A}^H_F(t)$, $\mathcal{A}^H_{\mu}(t)$ and $\mathcal{A}^H_{\mu_0}(t)$
			corresponding to the above universe exists respectively from the above theorems.
If we assume the absorbing set $\mathcal{B}$ belongs to the different subsets of compact set, bounded set $B$, fixed universe $\mathcal{D}^H_F$, universes $\mathcal{D}^H_{\mu}$ and $\mathcal{D}^H_{\mu_0}$ in $H$ with corresponding pullback attractors $\mathcal{A}^H_{cdf}(t),\ \mathcal{A}^H_B(t),\ \mathcal{A}^H_F(t),\ \mathcal{A}^H_{\mathcal{D}_{\mu}}(t)$ and $\mathcal{A}^H_{\mathcal{D}_{\mu_0}}(t)$ respectively, by the structure of pullback attractor as $\omega$-limit set, using the same technique in \cite{yf}, then it follows that
$\mathcal{A}^H_{cdf}(t)\subset\mathcal{A}^H_B(t)\subseteq\mathcal{A}^H_F(t)\subseteq \mathcal{A}^H_{\mathcal{D}_{\mu}}(t)\subseteq \mathcal{A}^H_{\mathcal{D}_{\mu_0}}(t)$. Here $\mathcal{A}^H_{cdf}(t)$ and $\mathcal{A}^H_B(t)$ coincide with the pullback attractors in \cite{cf97} and \cite{wzz} respectively.

\end{Remark}
\subsection{Main Result III: Finite fractal dimension of pullback attractors}\label{sub2.1}
Let $X$ be a separable real Hilbert space, $K\subset X$ be a non-empty compact subset and $\varepsilon>0$, we denote $N_{\varepsilon}(K)$ be the minimum number of open balls in $X$ with radius $\varepsilon$ which are necessary to cover $K$. The fractal dimension of $K$ is defined as $\mbox{d}_F(K)=\displaystyle{\limsup_{\varepsilon\rightarrow 0^+}}\frac{\ln(N_{\varepsilon}(K))}{\ln(\frac{1}{\varepsilon})}$.
Since the Hausdorff dimension $\mbox{d}_H(K)\leq \mbox{d}_F(K)$ which reveals the complexity of nonlinear systems such as the hydrodynamical models, we only need to present the fractal dimension of minimal family of pullback attractors for problem \eqref{a1}.

\begin{Theorem}\label{th5-f}
	Assume $f \in  L^2_{loc}(\mathbb{R};V')$ and satisfies uniformly pullback tempered condition,  let $v_{\tau}\in H$, then the $\mathcal{D}_{\mu}$-family of pullback attractors $\mathcal{A}^H_{\mathcal{D}_{\mu}}=\{\mathcal{A}^H_{\mathcal{D}_{\mu}}(t)\}$ in $H$ for the system \eqref{c1a} has bounded fractal dimension and Hausdorff dimension:
	
	(I) If $\frac{\pi \nu^2 n^2}{2|\Omega|}\geq\frac{C}{\nu^2}\Big\{\frac{\nu^2|\Omega|}{\epsilon}+\epsilon|\partial\Omega|\Big\}\|\varphi\|_{L^{\infty}(\partial\Omega)}+\frac{C}{\nu}M$, here $M=:\displaystyle{\lim_{T\rightarrow +\infty}}\frac{1}{T}\int^t_{t-T}\|f(r)\|^2_{V'}dr$, then $\mbox{dim}_F\mathcal{A}^H_{\mathcal{D}_{\mu}}\leq n$.
	
	(II) If $\frac{\pi \nu^2 n^2}{2|\Omega|}<\frac{C}{\nu^2}\Big\{\frac{\nu^2|\Omega|}{\epsilon}+\epsilon|\partial\Omega|\Big\}\|\varphi\|_{L^{\infty}(\partial\Omega)}+\frac{C}{\nu}M$ and $f\in L^{\infty}(-\infty,T^*;H)$, then the fractal and Hausdorff dimension of pullback attractors has bounded as $\mbox{dim}_F\mathcal{A}^H_{\mathcal{D}_{\mu}}\leq \hat{C}_1 Re+\hat{C}_2 G+\hat{C}_3$, here $\mbox{Re}=\frac{\|\varphi\|_{L^{\infty}(\partial\Omega)}}{\nu\lambda_1}$ is the generalized Reynold number, $\mbox{G}=\frac{\|f\|^2_{L^{\infty}(-\infty,T;H)}}{\nu^2\lambda_1}$ be the generalized Grashof number for the non-autonomous
system \eqref{a1}.
\end{Theorem}
{\bf Proof.} Based on the uniform differentiability of the process $U(\cdot,\cdot)$, we will apply the trace formula (see Lemmas 4.19 and 4.20 in \cite{clr}) to find the bound on the dimension by the generalized Reynold and Grashof numbers. The proof can be found in Section \ref{5.5}.$\hfill$$\Box$

\begin{Remark}
{\bf Consider the 2D Navier-Stokes equation in smooth domain:}

(1) If the Grashof and Reynolds numbers are defined as $G=\frac{|f|}{\nu^2\lambda_1}$ and $\mbox{Re}=\frac{|f|^{1/2}}{\nu\lambda^{1/2}_1}$ respectively, the fractal dimension of global attractor $\mathcal{A}$ in $H$ has finite dimension (see \cite{fmrt}, \cite{te88})
\begin{align}
	&\mbox{dim}_F\mathcal{A}\leq CG,
	\ C\leq (\frac{2}{\pi})^{1/2}(\lambda_1|\Omega|)^{1/2},\label{dm}\\
	&\mbox{dim}_F\mathcal{A}\leq C G^{2/3}(1+\log G)^{1/3}\label{dm1}
	\end{align}
	for Dirichlet boundary and periodic boundary conditions respectively. More delicate results about fractal dimension, we can see \cite{ci}, \cite{cf85}, \cite{cf1}, \cite{cfmt}, \cite{cft}, \cite{dg}, \cite{fmtt}, \cite{ily}, \cite{l1}, \cite{liu}, \cite{mir}, \cite{ro11}, \cite{yf}.

(2) From \cite{clr}, we can see that the fractal dimension for the fibre of pullback attractors $\mathcal{A}(\cdot)$ in $H$ can be estimated similarly as \eqref{dm}-\eqref{dm1} with the non-autonomous Grashof number $\mbox{G}(t)=\frac{\|f\|_{L^2(-\infty,t;V')}}{\nu^2\lambda_1}$ with Dirichlet boundary and periodic boundary conditions on bounded domain respectively, but their union $\displaystyle{\bigcup_{t\in\mathbb{R}}}\mathcal{A}(t)$ can be infinite dimension.
If $\Omega$ is an unbounded smooth domain where the Poincar\'e inequality holds, the fractal dimension of pullback attractors $\mathcal{A}(t)$ if also finite in $H$, i.e., $\mbox{dim}_F\mathcal{A}(t)\leq \max\{1, \frac{2\|f\|^2_{L^2(-\infty,t;V')}}{\nu^4\lambda_1}\}$
	for all $t\in \mathbb{R}$, see \cite{llr}.

{\bf For our problem in Lipschitz-like domain:}	The fibre of pullback attractor has finite fractal dimension which is dependent on the area of domain
as $\hat{C}_1 Re+\hat{C}_2 G+\hat{C}_3$, this is quite different from the finite fractal dimension on smooth domain of the 2D NSE. Comparing with autonomous case on non-smooth domain, the finite fractal dimension of global attractor is described as $C_1G+C_2Re^{3/2}$ in \cite{bps}, which has different exponent index of Reynold's number. Hence, we can conclude that Theorem \ref{th5-f} is an extension result of \cite{bps} to non-autonomous case.
\end{Remark}
\subsection{Main result IV: Continuity of pullback attractors}\label{sub2.2}
The upper semi-continuity of attractors is investigated originally by Hale and Raugel \cite{jkh}, Hale, Lin and Raugel \cite{jlh} in 1988, then the theory has been developed by many mathematicians which contains the situations:
(I) the external force has some perturbation for instance $\varepsilon \sigma(t,x)$ when $\varepsilon\rightarrow 0$ (see \cite{clr01}, \cite{jlh}, \cite{wq});
(II) the system has a damping perturbation, when the damping disappear, the corresponding continuity of attractors (see \cite{ctv}, \cite{ctv07}, \cite{cg});
(III) the external force can be written as singular oscillation $f_0(t,x)+\varepsilon^{-\rho}f_1(t/\varepsilon,x)$, the continuity of attractors when $\varepsilon\rightarrow 0$;
(IV) the coupled system converges to its uncoupled model (see \cite{mm}) and the continuous of attractors for different systems (see \cite{zd});
(V) the convergence of attractors with respect to perturbed domain (see \cite{zlz}) and so on.

In this paper, we consider an external force $f(t,x)=\sigma_0(x)+\delta \sigma(t,x)$, the sum of a stationary forces and a non-autonomous perturbation, which is more natural to physical applications. Then, for our problem, considering the special case $f(t,x)=\sigma_0(x)+\delta \sigma(t,x)$ of \eqref{a1}:
\begin{eqnarray}
	\left\{\begin{array}{ll} u_t-\nu\Delta u+(u\cdot\nabla)u+\nabla p=\sigma_0(x)+\delta \sigma(t,x),\ \ \  (x,t)\in \Omega_{\tau},\\
		 \mbox{div}\  u=0,\ \ \ \ (x,t)\in \Omega_{\tau},\\
		u(t,x)|_{\partial\Omega}=\varphi,\ \ \varphi\cdot n=0,\ \ \ \ (x,t)\in \partial\Omega,\\
		u(\tau,x)=u_{\tau}(x),\ \ \ \ x\in \Omega,\end{array}\right.\label{a1-1}\end{eqnarray}
and its abstract equivalent form in Lipschitz-like domain
\begin{eqnarray}
	\left\{\begin{array}{ll} 	v_t+\nu
		Av+B(v,v)+B(v,\psi)+B(\psi,v)+B(\psi,\psi)=P{(\sigma_0(x)+\nu F(x))+\delta P\sigma(t,x)},\\
		 \mbox{div}\  v=0,\\
		v|_{\partial\Omega}=0,\\
		 v(\tau,x)=v_{\tau}(x)=u_{\tau}(x)+\psi(x).\end{array}\right.\label{a1-2}\end{eqnarray}
The upper semi-continuity of attractors for above problem as the perturbation vanishes is given in the following theorem.
\begin{Theorem}\label{th6.1}
	(1) Assume $v_{\tau}\in H$, the external forces $f(t,x)=\sigma_0(x)+\delta \sigma(t,x)$, here $\sigma_0(x)\in H$ and $\sigma(t,x)\in L^2_{loc}(\mathbb{R};H)$ satisfies \begin{equation}\int^{\tau}_{-\infty}e^{\eta s}|\sigma(s)|^2ds<+\infty,\label{ya-1}\end{equation} then the process $\{U_{\delta}(t,\tau)\}$
	generated by the solution of \eqref{a1-2} possesses a family of pullback attractors
	$\mathcal{A}_{\delta}(t)$ with $\delta> 0$ in $H$.
	
	(2) The pullback
	attractors ${\mathcal A}_{\delta}=\{{\mathcal
		{A}}_{\delta}(t)\}_{t\in {\mathbb R}}$ of the processes for non-autonomous perturbed equation \eqref{a1-2}  with $\delta>0$ and the global attractor $\mathcal {A}$ for
	autonomous case	\eqref{a1-2} with $\delta=0$ satisfy the property of upper semi-continuity in $H$, i.e.,
	\begin{equation}\lim\limits_{\delta\rightarrow
		0^{+}}dist_{H}({\mathcal {A}}_{\delta}(t), {\mathcal {A}})=0
	\quad for~ any~ t\in {\mathbb R}.\label{f14}\end{equation}
\end{Theorem}
{\bf Proof.} See Section \ref{sec6}.$\hfill$$\Box$

\subsection{Main Result V: Regularity of pullback attractor}
\label{sub2.3}
\ \ \ \ Suppose the external force $f(t,x)$ is pullback translation bounded ($f \in L^2_{pb}(\mathbb{R};H)$) or uniformly pullback tempered ($f\in L^2_{upt}(\mathbb{R};H)$),
then we shall prove the regularity of pullback attractors in $D(A^{\frac{\sigma}{2}}) (\sigma\in [0,\frac{1}{2}])$, i.e., the $\mathcal{D}'_{\tilde{\mu}}$-family of pullback attractors in $D(A^{\frac{\sigma}{2}})$ can be stated as follows.

Since $D(A^{\frac{\sigma}{2}})$ is a Hilbert space, let $\mathcal{P}(D(A^{\frac{\sigma}{2}}))$ be the collection of all nonempty subsets in $D(A^{\frac{\sigma}{2}})$, $\hat{D}=\{D(t)\}\subset \mathcal{P}(D(A^{\frac{\sigma}{2}}))$ be a subset in $\mathcal{P}(D(A^{\frac{\sigma}{2}}))$, $B(0,\rho'_{\hat{D}}(t))$ be a family of balls with center $0$ and radius $\rho'_{\hat{D}}(t)$, where $\rho'_{\hat{D}}(t)$ satisfies that
\begin{equation}
e^{-\tilde{\mu} t}\Big(K^2_3+|\rho'_{\hat{D}}(\tau)|^2\Big)e^{\tilde{\mu}\tau}\rightarrow 0\label{gg6-1}
\end{equation}
as $\tau\rightarrow -\infty$.  Since $K^2_3$ is bounded, it suffices to verify that there exists a pullback time $\hat{\tau}_1\leq t$, such that for any $\tau<\hat{\tau}_1$,
\begin{equation}
|\rho'_{\hat{D}}(\tau)|^2e^{\tilde{\mu}\tau}\rightarrow 0.\label{gg7-1}
\end{equation}

We define the universe $\mathcal{D}'_{\tilde{\mu}}=\{\hat{D}=\{D(t)\}\}$ as
\begin{equation}
\mathcal{D}'_{\tilde{\mu}}=\{\hat{D}|D(t)\subset B(0,\rho'_{\hat{D}}(t))\ with\ \rho'_{\hat{D}}\ satisfying\ \eqref{gg7-1}\}.\label{gg7-1-1}
\end{equation}

Let $\hat{B}'_0=\{B'_0(t)\}_{t\in\mathbb{R}}$ be a family of balls, where $B'_0(t)=\overline{B'(0,\rho'_0(t))}$ is a ball with center $0$ and radius $\rho'_0(t)$, where $\rho'_0(t)$ is defined as
\begin{eqnarray}
&&(\rho'_0(t))^2=1+e^{-\tilde{\mu} t}K^2_3e^{\tilde{\mu}\tau}+\frac{2C}{\nu}\|f\|^2_{L^2_{loc}(\mathbb{R};H)}+\frac{2Ce^{-\tilde{\mu} n_0h}}{\nu(1-e^{-\tilde{\mu} n_0h})}\|f\|^2_{L^2_{pb}(\mathbb{R};H)}
\end{eqnarray}
here $0<\tilde{\mu}\leq \tilde{\mu}_0=\nu\lambda_1$.

\begin{Theorem}\label{th7.11}
	Assume $f(t,x)$ is uniformly pullback tempered in $L^2_{upt}(\mathbb{R};H)$ or pullback translation bounded in $L^2_{pb}(\mathbb{R};H)$,  $v_{\tau}\in D(A^{\frac{\sigma}{2}})$, then the norm-to-weak continuous processes $\{U(t,\tau)\}$ possesses a minimal $\mathcal{D}'_{\tilde{\mu}}$-family of pullback attractors $\mathcal{A}_{\mathcal{D}'_{\tilde{\mu}}}=\{\mathcal{A}_{\mathcal{D}'_{\tilde{\mu}}}(t)\}$ in $H$ for the system \eqref{c1a} which is equivalent to our original problem \eqref{a1}.
\end{Theorem}
{\bf Proof.} This is similar to the proof in Theorem \ref{th5.6}, see Section \ref{sub7.1}. $\hfill$$\Box$\\

Next, if we can show that the universe $\mathcal{D}'_{\tilde{\mu}}$ is inclusion closed and $\hat{B}'_0\in\mathcal{D}'_{\tilde{\mu}}$, by Theorem \ref{th8.35}, the family of pullback $\mathcal{D}'_{\tilde{\mu}}$-attractors is unique.
\begin{Theorem}\label{th7.12}
		Assume $f(t,x)$ is uniformly pullback tempered in $L^2_{upt}(\mathbb{R};H)$ or pullback translation bounded in $L^2_{pb}(\mathbb{R};H)$,  $v_{\tau}\in  D(A^{\frac{\sigma}{2}})$, let $0<\tilde{\mu}<\tilde{\mu}_0$, then the norm-to-weak continuous processes $\{U(t,\tau)\}$ possesses a unique $\mathcal{D}'_{\tilde{\mu}}$-family of pullback attractors $\mathcal{A}'_{\mathcal{D}'_{\tilde{\mu}}}=\{\mathcal{A}'_{\mathcal{D}'_{\tilde{\mu}}}(t)\}$ in $D(A^{\frac{\sigma}{2}})$ for the system \eqref{c1a} which is equivalent to \eqref{a1}.
\end{Theorem}
{\bf Proof.} See Section \ref{sub7.1}.$\hfill$$\Box$

\subsection{Coclusions}
From our main results and proof, we can see that the domain is important to the dynamical behavior for Navier-Stokes equation especially turbulence.
However, if the hydrodynamical systems defined on non-cylinder domain or complex thin domain, the dynamics and its continuity is still unknown.

\section{Proof of Main Results}\label{sec4}
\setcounter{equation}{0}
In this section, we will prove our main results by some delicate estimates.
\subsection{Proof of well-poseness}\label{sece}
In this section, we shall first prove the global existence of solutions for the equivalent abstract equation of problem \eqref{a1}. Then the proof of continuous dependence on the initial data and the solution processes will be presented.

\noindent $\bullet$ {\bf Proof of Theorem \ref{th4.1}}: {\bf Step 1:} We shall use the standard Faedo-Galerkin method to establish the existence of approximate solution to problem \eqref{c1a}.
Fix $n\geq1$, $w_j\ (j\geq 1)$ be the normalized eigenfunctions basis for the Stokes operators in the space $H$ with its increasing eigenvalues $\lambda_j\ (j\geq 0)$ being $0<\lambda_1\leq \lambda_2\leq \cdots $ and $\displaystyle{\lim_{j\rightarrow\infty}}\lambda_j=\infty$.
Let $V_n=span\{w_1, w_2, \cdots, w_n\}$, we define an approximate solution $v_{n}$ to problem \eqref{c1a} as $v_{n}(t)=\displaystyle{\sum_{j=1}^{n}}a_{nj}(t)w_{j}\in V_n$ which satisfies the following initial value problem of ordinary differential equation with respect to unknown variables $\{a_{nj}\}^{n}_{j=1}$,
\begin{equation}\left\{
	\begin{array}{ll}
		\frac{d}{dt}(v_{n}, w_j)+\nu
		\langle Av_{n}, w_j\rangle+b(v_{n},v_n, w_j)+b(v_n,\psi,w_j)+b(\psi,v_n,w_j)=(P_{n}\bar{f},w_j)-b(\psi,\psi,w_j),\\
		v_{n}(\tau)=v_{n\tau}.\end{array}\right.\label{c3}\end{equation}

By the local existence theory of solutions for ordinary differential equations, there exists a solution in local interval $(\tau,T)$ for problem \eqref{c3}.

{\bf Step 2:} The uniformly priori $L^{\infty}$-estimates.

Multiplying \eqref{c3} by $a_{nj}$, summing the resulting equations from $j=1$ to $n$, noting $b(v_{n},v_n, v_n)=0$ and $b(\psi,v_n,v_n)=0$ from \eqref{q7}, we have
\begin{eqnarray}
	&&\frac{1}{2}\frac{d}{dt}|v_{n}|^{2}+\nu\|v_{n}\|^{2}\leq |b(v_{n},\psi,v_{n})|+|\langle P_{n}\bar{f},v_{n}\rangle|+|b(\psi, \psi,v_{n})|.\label{c4}\end{eqnarray}

Next, we shall estimate every term on the right-hand side of \eqref{c4}. Similarly to the technique to \cite{bps}, we omit some details here.

(a) Using Hardy's inequality \eqref{7},\eqref{q7} and  \eqref{c5}, choosing suitable $\varepsilon>0$ such that \begin{equation}C_2C_3C_4\varepsilon\|\varphi\|_{L^{\infty}(\partial\Omega)}\leq\frac{\nu}{4},\label{p1}\end{equation}  we obtain
\begin{eqnarray}
	|b(v_{n},\psi,v_{n})|
	&\leq& C_4\|\varphi\|_{L^{\infty}(\partial\Omega)}\int_{dist(x,\partial\Omega)\leq C_{2}\varepsilon}dist(x,\partial\Omega)\frac{|v_{n}|^{2}}{[dist(x,\partial\Omega)]^{2}}dx\nonumber\\
	 &\leq&C_2C_3C_4\varepsilon\|\varphi\|_{L^{\infty}(\partial\Omega)}\|v_{n}\|^{2}\leq\frac{\nu}{4}\|v_{n}\|^{2}.\label{c6}\end{eqnarray}
Similarly, we have
\begin{eqnarray}
	&&|b(\psi,\psi,v_{n})|\leq\frac{\nu}{4}\|v_{n}\|^{2}+\frac{C\varepsilon\ |\partial\Omega|}{\nu}\|\varphi\|_{L^{\infty}(\partial\Omega)}^2,\label{c7}\\
	 &&|<P_{n}\bar{f},v_{n}>|\leq\frac{\nu}{4}\|v_{n}\|^{2}+\frac{C}{\nu}\Big[\frac{|f|^2_{V'}}{\lambda_{1}}+\frac{C\nu^2}{\varepsilon}\|\varphi\|^2_{L^{2}(\partial\Omega)}\Big].\label{c8}\end{eqnarray}
Combining \eqref{c4}--\eqref{c8}, by $\|\varphi\|_{L^{2}(\partial\Omega)}^2\leq C|\partial\Omega|\|\varphi\|_{L^{\infty}(\partial\Omega)}^2$ and the Poincar\'{e} inequality, we get
\begin{eqnarray} \frac{d}{dt}|v_{n}|^{2}+\frac{\nu\lambda_1}{2}|v_{n}|^{2}\leq \frac{C}{\nu}\Big[\frac{|f|^2_{V'}}{\lambda_{1}}+\frac{C\nu^2}{\varepsilon}|\partial\Omega|\|\varphi\|^2_{L^{\infty}(\partial\Omega)}+C\varepsilon\ |\partial\Omega|\|\varphi\|_{L^{\infty}(\partial\Omega)}^2\Big] \equiv K^{2}_{0}.\label{c9}\end{eqnarray}
By the Gronwall inequality, we derive
\begin{equation}
	|v_n|^2\leq |v_{n\tau}|^2e^{-\frac{\nu\lambda_1}{2}(t-\tau)}+\int^t_{\tau}e^{-\frac{\nu\lambda_1}{2}(t-s)}K^2_0ds.\label{nsv1}
\end{equation}
Moreover, since $\|\varphi\|_{L^{2}(\partial\Omega)}^2\leq C|\partial\Omega|\|\varphi\|_{L^{\infty}(\partial\Omega)}^2$ and $f\in L^2_{loc}(\tau,T;V')$, then we have, for an arbitrary $t\in (\tau,T)$,
\begin{eqnarray} &&\hspace{-1cm}\int^t_{\tau}e^{-\frac{\nu\lambda_1}{2}(t-s)}K^2_0ds=\int^t_{\tau}e^{-\frac{\nu\lambda_1}{2}(t-s)}\frac{C}{\nu}\Big[\frac{|f|^2_{V'}}{\lambda_{1}}+\frac{\nu^2}{\varepsilon}\|\varphi\|^2_{L^{2}(\partial\Omega)}+\varepsilon\ |\partial\Omega|\|\varphi\|_{L^{\infty}(\partial\Omega)}^2\Big]ds\nonumber\\
	 &&\hspace{2cm}=\int^t_{\tau}e^{-\frac{\nu\lambda_1}{2}(t-s)}\frac{C}{\nu}\frac{|f|^2_{V'}}{\lambda_1}ds\nonumber\\
	 &&\hspace{2.5cm}+\frac{C}{\nu}\Big[\frac{\nu^2}{\varepsilon}|\partial\Omega|\|\varphi\|^2_{L^{\infty}(\partial\Omega)}+\varepsilon\ |\partial\Omega|\|\varphi\|_{L^{\infty}(\partial\Omega)}^2\Big]\int^t_{\tau}e^{-\frac{\nu\lambda_1}{2}(t-s)}ds\nonumber
\end{eqnarray}
and
\begin{equation}
	 \frac{C}{\nu}\Big[\frac{\nu^2}{\varepsilon}|\partial\Omega|\|\varphi\|_{L^{\infty}(\partial\Omega)}^2+\varepsilon\ |\partial\Omega|\|\varphi\|_{L^{\infty}(\partial\Omega)}^2\Big]\int^t_{\tau}e^{-\frac{\nu\lambda_1}{2}(t-s)}ds\leq C, \label{nsv2}
\end{equation}
\begin{equation}
	 \frac{C}{\nu\lambda_1}\int^t_0e^{-\frac{\nu\lambda_1}{2}(t-s)}|f(s)|^2_{V'}ds\leq \frac{C}{\nu\lambda_1}\int^t_{\tau}|f(s)|^2_{V'}ds\leq C.\label{nsv3}
\end{equation}
 This implies that $v(t,x)\in L^{\infty}(\tau,T;H)$.

{\bf Step 3:} The priori $L^{2}$-estimate.

Integrating \eqref{c4} over $(s,t)$, using \eqref{p1}--\eqref{c8}, we have
\begin{eqnarray}
	&&|v_n|^2+ \frac{\nu}{2}\int^t_s\|v_{n}(r)\|^{2}dr\nonumber\\
	&&\leq |v_n(s)|^2+C\|f\|_{L^2(\tau,T;V')}+\frac{C}{\nu}\Big[\frac{\nu^2}{\varepsilon}|\Omega|\|\varphi\|^2_{L^{\infty}(\partial\Omega)}+\varepsilon\ |\partial\Omega|\|\varphi\|_{L^{\infty}(\partial\Omega)}^2\Big](t-s),\label{nsv4}
\end{eqnarray}
and using the Poincar\'{e} inequality and Gronwall's inequality, we derive that
\begin{eqnarray}
	|v_n|^2&\leq& |v_n(s)|^2e^{-\frac{\nu\lambda_1}{2}(t-s)}+C\|f\|_{L^2_{loc}(\tau,T;V')}e^{-\frac{\nu\lambda_1}{2}(t-s)}\nonumber\\
	 &&+\frac{C}{\nu}\Big[\frac{\nu^2}{\varepsilon}|\Omega|\|\varphi\|^2_{L^{\infty}(\partial\Omega)}+\varepsilon\ |\partial\Omega|\|\varphi\|_{L^{\infty}(\partial\Omega)}^2\Big]\frac{t-s}{e^{\frac{\nu\lambda_1}{2}(t-s)}},\label{yang-5}
\end{eqnarray}
which shows that $v_n(t,x)\in L^{\infty}(\tau,T;H)\cap L^2(\tau,T;V)$ if let $s=\tau$ in \eqref{yang-5}.

{\bf Step 4:} The priori $L^{2}$-estimate of $\frac{dv_n}{dt}$.

Using the inequality
$|b(u,v,w)|\leq C|u|^{\frac{1}{2}}\|u\|^{\frac{1}{2}}\|v\|\|w\|$,
we have
\begin{eqnarray}
	&&\|B(v_n,v_n)\|_{V'}=\sup_{\|w\|=1}\|b(v_n,v_n,w)\|\leq C|v|^{\frac{1}{2}}\|v_n\|^{\frac{1}{2}}\|v_n\|,\nonumber\\
&&\int^T_{\tau}\|B(v_n,v_n)\|_{V'}^2ds
	\leq C\Big(\int^T_{\tau}|v_n|^2ds+\int^T_{\tau}\|v_n\|^2ds\Big),\end{eqnarray}
which implies $B(v_n,v_n)\in L^2(\tau,T;V')$.

Considering $\langle \frac{dv_n}{dt}, w\rangle$ with $w\in V$ in the equation
\begin{eqnarray}
	&&\langle \frac{dv_n}{dt}, w\rangle+(\nu Av_n, w)+(B(v_n,v_n),w)+(B(v_n,\psi),w)+(B(\psi,v_n),w)\nonumber\\
	&&\hspace{5cm}=(P\overline{f},w)-(B(\psi,\psi),w),\label{gu1}
\end{eqnarray}
 and using the Hardy inequality, H\"{o}lder's inequality and the similar technique in \eqref{c6}, we deduce
\begin{eqnarray}
	|b(v_n,\psi,w|
	&\leq& C_4 \|\varphi\|_{L^{\infty}(\partial\Omega)}\int_{dist(x,\partial\Omega)\leq C_2\varepsilon}\frac{|v_n|}{dist(x,\partial\Omega)}|w|dist(x,\partial\Omega)dx\nonumber\\
	&\leq& C \|\varphi\|_{L^{\infty}(\partial\Omega)} \|v_n\| \|w\|
\end{eqnarray}
and
\begin{eqnarray}
	\int^T_{\tau}|B(v_n,\psi)|_{V'}^2d s
	\leq C \|\varphi\|_{L^{\infty}(\partial\Omega)}^2 \|v_n\|^2_{L^2(\tau,T;V)}
\end{eqnarray}
which means that $B(v_n,\psi)\in L^2(\tau,T;V')$.
Similarly, we have $B(\psi,v_n)\in L^2(\tau,T;V')$ and $B(\psi,\psi)\in L^2(\tau,T;V')$.
Since $f\in L^2_{loc}(\tau,T;V')$ and $v_n\in L^{2}(\tau,T;V)$, then $P\overline{f}\in L^2(\tau,T;V')$ and $\nu Av_n\in L^2(\tau,T;V')$. From the equation \eqref{gu1} in the weak sense, we have $\frac{dv_n}{dt}\in L^2(\tau,T;V')$ and $\{\frac{dv_n}{dt}\}$ is bounded.

{\bf Step 5:} The compact argument and existence of weak solutions.

From Steps 3 and 4, using the Lions-Aubin compact argument and the dominated convergence theorem, we can extract a subsequence (relabeled as $v_n$) and derive the existence of function $v\in L^2(\tau,T;V)\cap L^{\infty}(\tau,T;H)$ with $\frac{dv}{dt}\in L^2(\tau,T;V')$ such that
\begin{eqnarray}
	&&v_n\rightarrow v \ \ strongly \ \ in\ \ \ L^2(\tau,T;H),\label{gu2}\\
	&&v_n\rightharpoonup v \ \ weakly \ \ in\ \ \ L^2(\tau,T;V),\\
	&&v_n\rightharpoonup v \ \ weakly\ * \ \ in\ \ \ L^{\infty}(\tau,T;H),\\
	&&\frac{dv_n}{dt}\rightharpoonup \frac{dv}{dt}\  \ weakly\ \ in\ \ L^2(\tau,T;V').
\end{eqnarray}
Next, we shall deal with the convergence of trilinear operators.
Using the H\"{o}lder inequality and the property of trilinear operators, we obtain
\begin{eqnarray}
	\int^T_{\tau}|b(v_n,v_n,\omega_j)-b(v,v,\omega_j)|dt
	&=&
	\int^T_{\tau}|b(v_n,v_n-v,\omega_j)|dt\leq C\int^T_{\tau}\|v_n\|\|\omega_j\||v_n-v|dt\nonumber\\
	&\leq& C\|v_n\|_{L^2(\tau,T;V)}\|v_n-v\|_{L^2(\tau,T;H)}\rightarrow 0.
\end{eqnarray}
Similarly, we have
\begin{eqnarray}
	\int^T_{\tau}|b(v_n-v,v,\omega_j)|dt
	&\leq& \frac{C}{\lambda_1^{\frac{1}{2}}}\|v\|_{L^2(\tau,T;V)}\|v_n-v\|_{L^2(\tau,T;H)}\rightarrow 0,\\
	\int^T_{\tau}|b(v_n,\psi,\omega_j)-b(u,\psi,\omega_j)|dt&\leq& C\|\varphi\|_{L^{\infty}(\partial\Omega)}\|v_n-v\|_{L^2(\tau,T;H)}\rightarrow 0,\\
	\int^T_{\tau}|b(\psi,v_n,\omega_j)-b(\psi,u,\omega_j)|dt&\leq& C\|\varphi\|_{L^{\infty}(\partial\Omega)}\|v_n-v\|_{L^2(\tau,T;H)}\rightarrow 0.\label{gu3}
\end{eqnarray}

Combining \eqref{gu2}--\eqref{gu3}, passing to the limit of \eqref{gu1}, we conclude that $v(t,x)$ is a weak solution to problem \eqref{c1} in the interval $(\tau,T)$, i.e., there exist at least one global in time Hadamard weak solution to problem \eqref{c1}. From the property of background flows class $\psi_{\varepsilon}=\psi\in C^{\infty}(\Omega)$ satisfying \eqref{8} and $v=u-\psi$ and the solution $v$ for problem \eqref{c1} is obtained in Theorem \ref{th4.1} with initial data $v_{\tau}=u_{\tau}-\psi$, it is easy to check that $u$ satisfies the conditions (i) (ii) and (iii) in Definition 3.1 and $u(t,x)\in L^{\infty}(0,T;H)\cap L^2(0,T;V)$.
We then complete the proof of our Theorem \ref{th4.1}.$\hfill$$\Box$\\

\noindent $\bullet$ {\bf Proof of Theorem \ref{th4.3}}\\
Let $u_1(\cdot)$ and $u_2(\cdot)$ be two solutions to problem \eqref{a1} with corresponding initial data
$u^1_{\tau}$ and $u^2_{\tau}$ respectively and background flow functions $\psi_1$ and $\psi_2$, if we take $w=u_1-u_2$, then $w$ satisfies the problem:
\begin{eqnarray}\left\{\begin{array}{ll}\frac{dw}{dt}-\nu\Delta
		w+(u_1\cdot\nabla)u_1-(u_2\cdot\nabla)u_2=0,\\
		 \mbox{div}\  w=0,\ (x,t)\in \Omega_{\tau},\\
		w(t,x)|_{\partial\Omega}=0,\ (x,t)\in \partial\Omega_{\tau},\\
		 w(\tau,x)=w_{\tau}=u^1_{\tau}(x)-u^2_{\tau}(x),\end{array}\right.\label{q5}\end{eqnarray}
which can be written as
\begin{eqnarray}\left\{\begin{array}{ll} \frac{dw}{dt}+\nu A
		w+B(u_1,u_1)-B(u_2,u_2)=0,\\
		 \mbox{div}\  w=0.\end{array}\right.\label{q6}\end{eqnarray}
Let $\omega\in C^{\infty}_{0}(\Omega)$, $ \mbox{div}\ \omega=0$, from the condition (ii) in Definition \ref{de4.2}, we can derive
\begin{eqnarray}\hspace{-0.3cm}\frac{d}{dt}<u_1-u_2,\omega>-\nu<u_1-u_2,\Delta \omega>=\int_{\Omega}\sum^{2}_{i,j=1}(u^{i}_1u^{j}_1-u_2^{i}u_2^{j})\frac{\partial \omega^{i}}{x_{j}}dx.\label{y1-1}\end{eqnarray}
Obviously, \eqref{y1-1} holds for any $\omega\in V$.
In fact, from the condition (ii) and $<u_1-u_2,\Delta \omega>=-((u_1-u_2,\omega))$, we have
\begin{eqnarray} &&u_1-u_2=(u_1-\psi_1)-(u_2-\psi_2)+(\psi_1-\psi_2)\in L^{2}([0,T];V),\nonumber\\
&&\frac{d}{dt}(u_1-u_2)\in L^{2}([0,T];V'),\ \mbox{for}\ \omega\in V.\nonumber\end{eqnarray}
Let $\omega=w=u_1-u_2$ in \eqref{y1-1}, we have
\begin{eqnarray}\frac{1}{2}\frac{d}{dt}|\omega|^2+\nu\|\omega\|^2&\leq& C\int_{\Omega}|u_2||\omega||\nabla \omega|dx
	\leq\|u_2\|_{L^4}|\nabla \omega|^{1/2}|w|^{1/2}|\nabla \omega|\nonumber\\
	 &\leq&\nu\|\omega\|^{2}+C_{\nu}\|u_2\|^{4}_{L^4}|\omega|^{2}.\label{gu5}\end{eqnarray}
Since
\begin{eqnarray}
	\|u_2\|_{L^4(\Omega)}&\leq& \|u_2-\psi\|_{L^4(\Omega)}+\|\psi\|_{L^4(\Omega)}\nonumber\\
	&\leq& C\|\nabla(u_2-\psi)\|^{\frac{1}{2}}_{L^2(\Omega)}\|u_2-\psi\|^{\frac{1}{2}}_{L^2(\Omega)}+\|\psi\|_{L^4(\Omega)},
\end{eqnarray}
we have $u\in L^{4}(\Omega\times(0,T))$ and $\omega(\cdot,0)=0$, hence $\omega=0$, i.e., the solution is unique.
Moreover, we have the continuous dependence on the initial data
\begin{eqnarray}
	|u_1(s)-u_2(s)|^2\leq \|u^1_0-u^2_0\|^2_H\times e^{C_{\nu}\int^t_{\tau}\|u_2(s)\|^4_{L^4(\Omega)}ds}.
	\label{q15}\end{eqnarray}

Note that \eqref{gu5} can be written as
\begin{equation}
	\frac{d}{dt}|\omega|^2+\nu\|\omega\|^2\leq C'_{\nu}\|u_2\|^4_{L^4}|\omega|^2.\end{equation}
Integrating from $\tau$ to $t$, we have
\begin{equation}
	|\omega(t)|^2+\nu\int^t_{\tau}\|\omega(s)\|^2ds\leq |\omega_{\tau}|^2+C'_{\nu}\int^t_{\tau}\|u_2(s)\|^4_{L^4}|\omega(s)|^2ds.\label{gu6}\end{equation}
Neglecting the first term on the left-hand side of \eqref{gu6}, and using \eqref{q15}, we derive
\begin{eqnarray}
	\int^t_{\tau}\|u_1(s)-u_2(s)\|^2ds\leq \frac{1}{\nu}|u^1_{\tau}-u^2_{\tau}|^2\times \Big(C'_{\nu}\int^t_{\tau}\|u_2(s)\|^4_{L^4}dse^{C_{\nu}\int^t_{\tau}\|u_2(s)\|^4_{L^4(\Omega)}ds}+1\Big).\label{q17}
\end{eqnarray}

Hence, \eqref{q15} and \eqref{q17} imply the continuous dependent on the initial data for the global weak solutions, and hence $u(t,x)\in C(0,T;H)\cap L^2(0,T;V)$, which completes the proof. $\hfill$$\Box$

\subsection{Theory of pullback attractors}\label{ap}
In this section, we shall recall the theory of pullback attractors which can be found in \cite{clr2006}, \cite{clr}, \cite{gmr}.

\noindent $\bullet$ {\bf Continuity of Processes}\label{sub8.3}
\begin{Definition}\label{de8.20}
	Let $X$ be a Banach space and the bi-parameters operators $\{U(t,\tau)|t\geq\tau\}$ be a family of processes on $X$.
	
	(1) (Pullback strong continuous process) We say that the process is pullback strong continuous if
	for given $t\in \mathbb{R}$ and fixed, $x_n\rightarrow x$ in $X$, we have that $U(t,\tau_n)x_n\rightarrow U(t,\tau)x$ for any $\tau_n\rightarrow \tau\in (-\infty,t]$ in $X$.

	(2) (Pullback weak continuous process) The process is pullback weak continuous if for given $t\in \mathbb{R}$ and fixed, $x_n\rightarrow x$ in $X$, we have that $U(t,\tau_n)x_n\rightharpoonup U(t,\tau)x$ for any $\tau_n\rightarrow \tau\in (-\infty,t]$ in $X$.

	(3) (Pullback norm-to-weak continuous process-Wang and Zhong \cite{wz}) The bi-parameters operators $\{U(t,\tau)|t\geq\tau\}$ be family of processes on  $X$, we say that $\{U(t,\tau)|t\geq\tau\}$ is pullback norm-to-weak continuous on $X$ if
	(i) $U(\tau,\tau)=Id$,
	(ii) $U(t,s)U(s,\tau)=U(t,\tau)$,
	(iii) for given $t\in \mathbb{R}$ and fixed, $x_n\rightarrow x$ in $X$, $\tau_n\rightarrow \tau$, it yields $U(t,\tau_n)x_n\rightharpoonup U(t,\tau)x$.

	(4)	(Pullback closed process-Garci\'{a}-Luengo, Mar\'{i}n-Rubio and Real-\cite{gmr}) A process $U(\cdot,\cdot)$ on $X$ is said to be closed if for any $\tau\leq t$ and any sequence $\{x_n\}\subset X$ with $x_n\rightarrow x\in X$ and $U(t,\tau)x_n\rightarrow y\in X$, then $y=U(t,\tau)x$.
\end{Definition}

\begin{Remark}\label{re8.23}
	(1) The pullback norm-to-weak continuous process is weaker than pullback strong and weak continuous processes;
	
	(2) The pullback closed process $U(t,\tau)$ is weaker than norm-to-weak continuous process.
	
	(3) The pullback closed is most generalized pullback continuity for the process.
\end{Remark}

\noindent $\bullet$ {\bf Pullback attractors for dissipative systems}\label{sub8.4}\\
Let us denote by $\mathcal{P}(X)$ be the family of all nonempty subsets of $X$, and consider a family of nonempty sets: $\hat{D}_0=\{D_0(t)|t\in\mathbb{R}\}\subset \mathcal{P}(X)$, here we do not require any additional condition on these sets such as the boundedness or compactness.

\begin{Definition}\label{de8.24}
	(Pullback asymptotically compact) A family of process $U(\cdot,\cdot)$ on $X$ is called pullback $\hat{D}_0$-asymptotically compact if for any $t\in\mathbb{R}$, and any sequence $\{\tau_n\}\subset (-\infty,t]$ and $\{x_n\}\subset X$, the sequence $\{U(t,\tau_n)x_n\}$ is relatively compact in $X$ for $\tau_n\rightarrow -\infty$ and $x_n\in D_0(\tau_n)$ as $n\rightarrow +\infty$.
\end{Definition}

\begin{Definition}\label{de8.25} (Universe) Let $\mathcal{D}$ be a nonempty class of families with parameters in time $\hat{D}=\{D(t)|t\in\mathbb{R}\}\subset \mathcal{P}(X)$, the class $\mathcal{D}$ is called a universe in $\mathcal{P}(X)$.\end{Definition}

\begin{Definition}\label{de8.26}
	(Inclusion closed) We say the family $\mathcal{D}$ is inclusion closed , if for any $\hat{D}\in\mathcal{D}$, and $\hat{C}=\{C(t)\}_{t\in\mathbb{R}}\in 2^X$, such that $C(t)\subset X,\ C(t)\neq \Phi,\ C(t)\subset D(t)$ for all $t\in\mathbb{R}$, then one has $\hat{C}\in\mathcal{D}$.
\end{Definition}

\begin{Definition}\label{de8.27}
	(Dissipation-pullback $\mathcal{D}$-absorbing) The set class $\hat{D}_0=\{D_0(t)|t\in\mathbb{R}\}\subset \mathcal{P}(X)$ is called pullback $\mathcal{D}$-absorbing for the process $U(\cdot,\cdot)$ on $X$ if for any $\hat{D}\in \mathcal{D}$ and fixed $t\in \mathbb{R}$, there exists a pullback time $\tau_0(t,\hat{D})\leq t$ such that $U(t,\tau)D(\tau)\subset D_0(t)$  for all $\tau \leq \tau_0(t,\hat{D})$.
\end{Definition}

\begin{Definition}\label{de8.28}
	(Pullback $\mathcal{D}$-asymptotic compactness) A family of process $U(\cdot,\cdot)$ on $X$ is called pullback $\mathcal{D}$-asymptotically compact if the process is $\hat{D}$-asymptotically compact for any $ \hat{D}\in\mathcal{D}$, i.e., if for any $t\in\mathbb{R}$, any $\hat{D}\in\mathcal{D}$, and any sequence $\{\tau_n\}\subset (-\infty,t]$ and $\{x_n\}\subset X$ satisfying $\tau_n\rightarrow -\infty$ and $x_n\in D(\tau_n)$ for all $n$, the sequence $\{U(t,\tau_n)x_n\}$ is relatively compact in $X$.
\end{Definition}

\begin{Definition}\label{de8.29}
	(Pullback $\mathcal{D}$-condition (MWZ)-see Wang and Zhong \cite{wz})  Let $X$ be a Banach space or complete metric space, we call the family of processes $U(t,\tau):X\rightarrow X$ satisfies pullback $\mathcal{D}$-condition (MWZ) if for any $\mathcal{B}=\{B(t)\}_{t\in\mathbb{R}}\in \mathcal{D}$,  for any $t\in\mathbb{R}$ and fixed, any $\varepsilon>0$, there exists a pullback time $\tau_{\varepsilon}=\tau(t,\varepsilon,\mathcal{B})\leq t$ and a finite dimensional subspace $X_1\subset X$ such that
	\begin{eqnarray}
	&&(i)\ P(\cup_{s\leq \tau_{\varepsilon}}U(t,s)B(s))\ is\ bounded,\nonumber\\
	&&(ii)\ \|(I-P)(\cup_{s\leq \tau_{\varepsilon}}U(t,s)B(s))\|_X\leq \varepsilon,\nonumber
	\end{eqnarray}
	where $P:X\rightarrow X_1$ is bounded projector.
\end{Definition}

\begin{Theorem}\label{re8.30}
	(Wang and Zhong \cite{wz}) For the process $U(t,\tau)$, the pullback $\mathcal{D}$-condition (MWZ) implies pullback $\mathcal{D}$-asymptotically compact, but the converse is not true. Moreover, if $X$ is a uniformly convex Banach space, such as $X$ is a Banach space, then they are equivalent.
\end{Theorem}

\begin{Theorem}\label{th8.31}
	If $\hat{D}_0$ is pullback $\mathcal{D}$-absorbing of $U(\cdot,\cdot)$, $U(\cdot,\cdot)$ is pullback $\hat{D}_0$-asymptotically compact (pullback $\mathcal{D}_0$-condition (MWZ)), then the process $U(\cdot,\cdot)$ is pullback $\mathcal{D}$-asymptotically compact.
\end{Theorem}
Let $\Lambda(\hat{D}_0,t)=\displaystyle{\cap_{s\leq t}}\overline{\displaystyle{\cup_{\tau\leq s}}U(t,\tau)D_0(\tau)}^X$ for all $t\in\mathbb{R}$, similarly, we can define  $\Lambda(\hat{D},t)=\displaystyle{\cap_{s\leq t}}\overline{\displaystyle{\cup_{\tau\leq s}}U(t,\tau)D(\tau)}^X$.
\begin{Theorem}\label{th8.35}
	Let $X$ be a complete metric space or Banach space, $\{U(t,\tau)\}: \mathbb{R}^2_d\times X\rightarrow X$ be a process, $\mathcal{D}=\{\hat{D}=\{D(t)\}_{t\in\mathbb{R}}\}\subset \mathcal{P}(X)$ be the universe. Suppose that the process satisfies:
	
	(1) $\{U(t,\tau)\}$ is norm-to-weak or closed;
	
	(2) $\{U(t,\tau)\}$ is pullback $\mathcal{D}_0$-asymptotically compact ((pullback $\mathcal{D}_0$-condition (MWZ)));
	
	(3) $\{U(t,\tau)\}$ admits a pullback $\mathcal{D}$-absorbing family $\hat{D}_0=\{D_0(t)\}_{t\in\mathbb{R}}$ which are not necessary in the universe.
	
	Then we have $U(t,\tau)$ possesses a minimal pullback $\mathcal{D}$-attractor $\mathcal{A}_{\mathcal{D}}=\{A_{\mathcal{D}}(t)\}$, here
	 $A_{\mathcal{D}}(t)=\overline{\displaystyle{\bigcup_{\hat{D}\in\mathcal{D}}}\Lambda(\hat{D},t)}^X$, $\mathcal{A}_{\mathcal{D}}=\{A_{\mathcal{D}}(t)\}$ is a family satisfying:
	
	(a) for any $t\in\mathbb{R}$, $A_{\mathcal{D}}(t)$ is a nonempty compact subset in $X$, and $A_{\mathcal{D}}(t)\subset \Lambda(\hat{D},t)$;
	
	(b) $A_{\mathcal{D}}(t)$ is pullback $\mathcal{D}$-attracting, i.e., $\displaystyle{\lim_{\tau\rightarrow-\infty}}dist_X(U(t,\tau)D(\tau),A_{\mathcal{D}}(t))=0$
	for all $\hat{D}\in\mathcal{D}$ and $t\in\mathbb{R}$;
	
	(c) $A_{\mathcal{D}}(t)$ is invariant, i.e., $U(t,\tau)A_{\mathcal{D}}(\tau)=A_{\mathcal{D}}(t)$ for all $\tau\leq t$;
	
	(d) if $\hat{D}_0\in\mathcal{D}$, then $A_{\mathcal{D}}(t)=\Lambda(\hat{D}_0,t)\subset \overline{D_0(t)}^X$ for all $t\in\mathbb{R}$.

	The family $\mathcal{A}_{\mathcal{D}}$ is minimal in the sense that if $\hat{C}=\{C(t)|t\in\mathbb{R}\}\subset \mathcal{P}(X)$ is a family of closed sets such that for any $\hat{D}=\{D(t)|t\in\mathbb{R}\}\subset \mathcal{D}$, $$\displaystyle{\lim_{\tau\rightarrow-\infty}}dist_X(U(t,\tau)D(\tau),C(t))=0,$$ then $A_{\mathcal{D}}(t)\subset C(t)$.
	Moreover, if $\hat{D}_0\in\mathcal{D}$ and $\mathcal{D}$ is inclusion closed, then $\mathcal{A}$ is unique and $\mathcal{A}\in\mathcal{D}$.
\end{Theorem}
{\bf Proof.} See,  e.g. \cite{gmr}. $\hfill$$\Box$

\subsection{Proof of Theorems \ref{th5.6} and \ref{th5.7}: minimal and unique family of pullback attractors in $H$}\label{sec5}
 \noindent $\bullet$ {\bf Step 1: Processes and its continuity}\label{sub5.1}\\
In \eqref{c1}, assume that $f\in  L^2_{loc}(\mathbb{R};V')$ and $v_{\tau}\in H$, then the problem admits a unique weak solution denoted as $v=v(\cdot,\tau,v_{\tau})$.
Hence, we can define a family of processes $U(\cdot,\cdot):\mathbb{R}^2_d\times H\rightarrow H$ as $U(t,\tau)v_{\tau}=v(t,\tau,v_{\tau})$.

\begin{Lemma}\label{th5.1}
	(Continuity of processes) If $v_{\tau}\in H$, $f\in L^2_{loc}(\mathbb{R};V')$, then $U(t,\tau)$ is pullback norm-to-weak continuous from $H$ to $H$ for problem \eqref{c1}.
\end{Lemma}
{\bf Proof.} The pullback norm-to-weak continuity of processes in $H$ is easily verified from Theorem \ref{th4.3}. For the process $v(t,\tau,v_{\tau})=U(t,\tau)v_{\tau}$ with $\tau_n\rightarrow \tau$ and $(v_{\tau})_n \rightarrow v_{\tau} $ in $H$ as $n\rightarrow \infty$, since $v\in L^2(0,T;V)\cap C(0,T,H)$ and $v$ is the weak solution of problem \eqref{c1a},
we can easily get the weak convergence subsequence of $U(t,\tau)(v_{\tau})_n$ from the energy equality
\begin{eqnarray}
	 &&|v_n(t)|^2+2\nu\int^t_{\tau}\|v_n(s)\|^2ds=|(v_{\tau})_n|^2-2[(B(v_n,\psi),v_n)+(B(\psi,v_n),v_n)\nonumber\\
	&&\hspace{5.4cm}+(B(\psi,\psi),v_n)]+2(P\bar{f}, v_n),
\end{eqnarray}
the technique of passing to the limit is similar to the proof of existence of weak solutions, here we omit the detail.
\\

\noindent $\bullet$ {\bf Step 2: Pullback $\mathcal{D}_{\mu}$-absorbing family of sets in $H$} \label{sub5.2a}
\begin{Lemma}\label{le5.3}
	Assume $f \in  L^2_{loc}(\mathbb{R};V')$ is pullback translation bounded in $L^2_{pb}(\mathbb{R};V')$ (or uniformly pullback tempered), $v_{\tau}\in H$, if we choose parameter $\mu\in (0,\mu_0=\frac{\nu\lambda_1}{2}]$ and fixed, then the solution $v$ to the problem \eqref{c1} satisfies that for any $\tau\leq t$,
	\begin{equation}
		|v|^2\leq |v_{\tau}|^2e^{-\mu(t-\tau)}+e^{-\mu t}\int^t_{\tau}e^{\mu s}K^2_0ds,
	\end{equation}
	here \begin{eqnarray}  K^{2}_{0}=\frac{4}{\nu}\Big[\frac{\|f\|^2_{V'}}{\lambda_{1}}+\frac{C\nu^2}{\varepsilon}|\partial\Omega|\|\varphi\|^2_{L^{\infty}(\partial\Omega)}+C\varepsilon\ |\partial\Omega|\|\varphi\|_{L^{\infty}(\partial\Omega)}^2\Big].\end{eqnarray}
	Moreover,
	\begin{eqnarray}
		&&\hspace{-1.5cm}|v(t)|^2\leq \Big(\frac{4}{\nu\mu}\Big[\frac{C\nu^2}{\varepsilon}|\partial\Omega|\|\varphi\|^2_{L^{\infty}(\partial\Omega)}+C\varepsilon\ |\partial\Omega|\|\varphi\|_{L^{\infty}(\partial\Omega)}^2\Big]+|v_{\tau}|^2\Big)e^{-\mu(t-\tau)}\nonumber\\
		 &&+\frac{4}{\nu\mu}\Big[\frac{C\nu^2}{\varepsilon}|\partial\Omega|\|\varphi\|^2_{L^{\infty}(\partial\Omega)}+C\varepsilon\ |\partial\Omega|\|\varphi\|_{L^{\infty}(\partial\Omega)}^2\Big]\nonumber\\
&&+\|f\|^2_{L^2_{loc}(\mathbb{R};V')}+\frac{4e^{-\mu n_0h}}{\nu\lambda_1(1-e^{-\mu n_0h})}\|f\|^2_{L^2_{pb}(\mathbb{R};V')}.\label{gg5}\end{eqnarray}
\end{Lemma}
{\bf Proof.} Using the same technique in the proof of Theorem \ref{th4.1} and noting that $\tau=t-n_0 h$, it follows
\begin{eqnarray}
&&\int^t_{-\infty}e^{-\mu(t-s)}\|f\|^2_{V'}ds\nonumber\\
&=&\int^t_{\tau}e^{-\mu(t-s)}\|f\|^2_{V'}ds+\int^{t-n_0h}_{t-(n_0+1)h}e^{-\mu(t-s)}\|f\|^2_{V'}ds+\int^{t-(n_0+1)h}_{t-(n_0+2)h}e^{-\mu(t-s)}\|f\|^2_{V'}ds+\cdots\nonumber\\
&\leq& \int^t_{\tau}\|f\|^2_{V'}ds+(e^{-\mu n_0h}+e^{-2\mu n_0 h}+\cdots)\|f\|^2_{L^2_{pb}(\mathbb{R};V')}\nonumber\\
&\leq&\int^t_{\tau}\|f\|^2_{V'}ds+\frac{e^{-\mu n_0h}}{1-e^{-\mu n_0h}}\|f\|^2_{L^2_{pb}(\mathbb{R};V')},
\end{eqnarray} we can easily get the result. $\hfill$$\Box$\\

\begin{Lemma}\label{th5.4}
	Assume $f \in  L^2_{loc}(\mathbb{R};V')$ is pullback translation bounded in $L^2_{pb}(\mathbb{R};V')$ (or uniformly pullback tempered), for any small enough $\varepsilon_1>0$, there exists a pullback time $\tau(t,\varepsilon_1)$, such that for any $\tau< \tau(t,\varepsilon_1)\leq t$, $\hat{B}_0(t)$ is a $\mathcal{D}_{\mu}$-family of pullback absorbing sets for the process $U(t,\tau)$.
\end{Lemma}
{\bf Proof.} Noting that \begin{eqnarray}
	&&\hspace{-1.5cm}|U(t,\tau)v_{\tau}|^2\leq \Big(\frac{4}{\nu\mu}\Big[\frac{C\nu^2}{\varepsilon}|\partial\Omega|\|\varphi\|^2_{L^{\infty}(\partial\Omega)}+C\varepsilon\ |\partial\Omega|\|\varphi\|_{L^{\infty}(\partial\Omega)}^2\Big]+|\rho_{\hat{D}}(t)|^2\Big)e^{-\mu(t-\tau)}\nonumber\\
	 &&\hspace{1cm}+\frac{4}{\nu\mu}\Big[\frac{C\nu^2}{\varepsilon}|\partial\Omega|\|\varphi\|^2_{L^{\infty}(\partial\Omega)}+C\varepsilon\ |\partial\Omega|\|\varphi\|_{L^{\infty}(\partial\Omega)}^2\Big]\nonumber\\
&&\hspace{1cm}+\|f\|^2_{L^2_{loc}(\mathbb{R};V')}+\frac{4e^{-\mu n_0h}}{\nu\lambda_1(1-e^{-\mu n_0h})}\|f\|^2_{L^2_{pb}(\mathbb{R};V')}
\end{eqnarray}
and
there exists a pullback time $\tau(t,\varepsilon_1)$, such that for any $\tau<\tau(t,\varepsilon_1)\leq t$, it follows
\begin{eqnarray}
	\hspace{-0.35cm}e^{-\mu t}\Big(\frac{4}{\nu\mu}\Big[\frac{C\nu^2}{\varepsilon}|\partial\Omega|\|\varphi\|^2_{L^{\infty}(\partial\Omega)}+C\varepsilon\ |\partial\Omega|\|\varphi\|_{L^{\infty}(\partial\Omega)}^2\Big]+|\rho_{\hat{D}}(\tau)|^2\Big)e^{\mu\tau}\leq \varepsilon_1
\end{eqnarray}
since
$e^{-\mu_0 (t-s)}\leq e^{-\mu(t-s)}$
holds for any $\mu\leq \mu_0$,
 we have \begin{equation}
	|U(t,\tau)v_{\tau}|^2\leq \varepsilon_1+\rho^2_0(t)-\frac{1}{2}\leq \rho^2_0(t),
\end{equation}
which implies that $U(t,\tau)D(\tau)\subset B_0(t)$, i.e., $\hat{B}_0(t)$ is a family of pullback $\mathcal{D}_{\mu}$-absorbing sets. \\

\noindent $\bullet$ {\bf Step 3: Pullback $\mathcal{D}_{\mu}$-asymptotic compactness for the pullback norm-to-weak continuous processes in $H$}\label{sub5.3}\\
 In this section, we shall prove the pullback $\mathcal{D}_{\mu}$-condition-(MWZ) (see Definition \ref{de8.29} and Remark \ref{re8.30}) to achieve the pullback $\mathcal{D}_{\mu}$-asymptotic compactness for the processes.

\begin{Lemma}\label{th5.5}
	Assume $f \in  L^2_{loc}(\mathbb{R};V')$ is pullback translation bounded in $L^2_{pb}(\mathbb{R};V')$ (or uniformly pullback tempered),  $v_{\tau}\in H$, then the processes $U(t,\tau)$ satisfies pullback $\mathcal{D}_{\mu}$-condition (MWZ) which implies pullback $\mathcal{D}_{\mu}$-asymptotically compact in $H$ for the system \eqref{c1a} which is equivalent to problem \eqref{a1}.
\end{Lemma}
{\bf Proof.} {\bf Step 1:} Let $\mathcal{B}=\{B_0(t)\}_{t\in\mathbb{R}}$ be the pullback $\mathcal{D}_{\mu}$-absorbing family given in Theorem \ref{th5.4}, then there exists a pullback time $\tau_{t,\varepsilon_1}$ such that
$|v|^2=|U(t,\tau)v_{\tau}|^2\leq \rho_{0}(t)$.

{\bf Step 2:} Since $H$ is a Hilbert space, the processes are defined as $U(t,\tau): H\rightarrow H$. Let $H=H_1\bigoplus H_2$, here $H_1=span\{\omega_1,\omega_2,\cdots, \omega_m\}$ and $H_1\bot H_2$ which means $H_2$ is the orthogonal complement space of $H_1$.

Let $\tilde{P}$ be an orthonormal projector from $H$ to $H_1$, then we have the decomposition
\begin{equation}
	v=\tilde{P}v+(I-\tilde{P})v:=v_1+v_2,
\end{equation} for $v_1\in H_1$, $v_2\in H_2$.

From the existence of global solutions and pullback $\mathcal{D}_{\mu}$-absorbing family of sets, we know $\|v_1\|^2\leq \rho_{0}(t)$, and next we only need to prove the $H$-norm of $v_2$ is small enough as $\tau\rightarrow -\infty$.

{\bf Step 3:} Taking inner product of \eqref{c1a} with $v_2$, since $(v_1, v_2)=0$, we have
\begin{eqnarray}&&\hspace{-1cm}(v_t,v_2)+(\nu A v, v_2)+(B(v),v_2)+(B(v,\psi),v_2)+(B(\psi,v),v_2)\nonumber\\
	&&\hspace{2.5cm}=-(B(\psi),v_2)+<\tilde{P} f,v_2>-<\nu \tilde{P} F,v_2>,\label{g2-a}\end{eqnarray}
that is,
\begin{eqnarray}&&\frac{1}{2}\frac{d}{dt}|v_2|^2+\nu\|v_2\|^2\leq |(B(v),v_2)|+|(B(v,\psi),v_2)|+|(B(\psi,v),v_2)|\nonumber\\
	&&\hspace{4cm}+|(B(\psi),v_2)|+|<\tilde{P} f,v_2>|+|<\nu \tilde{P} F,v_2>|.\label{g2}\end{eqnarray}

Next, we shall estimate every term on the right-hand side of \eqref{g2}.\\
(a) Using the property of trilinear operators, \eqref{q7} and the $\varepsilon$-Young inequality, we have
\begin{eqnarray}
	|(B(v),v_2)|&=&|b(v,v_2,v)|\leq\int_{\Omega}|v||\nabla v_2||v|dx\nonumber\\
	&\leq&C\|v\|_H^{\frac{1}{2}}\|\nabla v_2\|_H^{\frac{1}{2}}\|v\|_H\leq\frac{\nu}{12}\|v_2\|^{2}+\frac{C}{\nu}\|v\|_H^2\nonumber\\
	&\leq& \frac{\nu}{12}\|v_2\|^{2}+\frac{C}{\nu}\rho_0(t).\end{eqnarray}
(b) By the similar technique in (b), we derive
\begin{eqnarray}
	|(B(\psi,v),v_2)|&\leq&\int_{\Omega}|\psi||\nabla v_2||v|dx\leq C_4\|\varphi\|_{L^{\infty}(\partial\Omega)}\int_{\Omega}|\nabla v_2||v|dx\nonumber\\
	&\leq& \frac{\nu}{12}\|v_2\|^{2}+\frac{C\|\varphi\|^2_{L^{\infty}(\partial\Omega)} }{\nu}\rho_0(t).\end{eqnarray}
(c) By virtue of the property of trilinear operators, \eqref{4} and the Hardy inequality, we get
\begin{eqnarray}
	|(B(v,\psi),v_2)|&\leq&\int_{\Omega}|v||\nabla \psi||v_2|dx\nonumber\\
	&\leq& \int_{dist(x,\partial\Omega)\leq C'_2\varepsilon}|v||\nabla \psi|\cdot dist(x,\partial\Omega)\frac{|v_2|}{dist(x,\partial\Omega)}dx\nonumber\\
	&\leq& C_4\|\varphi\|_{L^{\infty}(\partial\Omega)}\int_{\Omega}\frac{|v_2|}{dist(x,\partial\Omega)}|v|dx\nonumber\\
	&\leq& C_4\|\varphi\|_{L^{\infty}(\partial\Omega)}\Big(\int_{\Omega}\frac{|v_2|^2}{dist^2(x,\partial\Omega)}\Big)^{\frac{1}{2}}\Big(\int_{\Omega}|v|^2dx\Big)^{\frac{1}{2}}\nonumber\\
	&\leq& C_4\|\varphi\|_{L^{\infty}(\partial\Omega)} \|A^{\frac{1}{2}}v_2\|_{H}\|v\|_H\nonumber\\
	&\leq& C\|\varphi\|^2_{L^{\infty}(\partial\Omega)} \frac{1}{\nu}\|v\|^2_H+\frac{\nu}{12}\|v_2\|^{2}\nonumber\\
	&\leq& \frac{\nu}{12}\|v_2\|^{2}+\frac{C\|\varphi\|^2_{L^{\infty}(\partial\Omega)} }{\nu}\rho_0(t).\end{eqnarray}
(d) By the Hardy inequality and the property of trilinear operators, we obtain
\begin{eqnarray}
	|(B(\psi,\psi),v_2)|&\leq&\int_{\Omega}|\psi||\nabla \psi||v_2|dx\nonumber\\
	&=& \int_{\Omega}|\psi||\nabla \psi|\cdot dist(x,\partial\Omega)\frac{|v_2|}{dist(x,\partial\Omega)}dx\nonumber\\
	&\leq& C_4\|\varphi\|_{L^{\infty}(\partial\Omega)}\int_{dist(x,\partial\Omega)\leq C'_2\varepsilon}|\psi|\frac{|v_2|}{dist(x,\partial\Omega)}dx\nonumber\\
	&\leq& C_4\|\varphi\|_{L^{\infty}(\partial\Omega)}\Big(\int_{\Omega}\frac{|v_2|^2}{dist^2(x,\partial\Omega)}\Big)^{\frac{1}{2}}\Big(\int_{dist(x,\partial\Omega)\leq C'_2\varepsilon}|\psi|^2dx\Big)^{\frac{1}{2}}\nonumber\\
	&\leq& \frac{\nu}{12}\|v_2\|^{2}+\frac{C\|\varphi\|^4_{L^{\infty}(\partial\Omega)}|\partial\Omega|\varepsilon}{\nu}.\label{ya-10}\end{eqnarray}
(e) By the Hardy inequality and Lemma \ref{le3.4}, since $\|\varphi\|_{L^{2}(\partial\Omega)}^2\leq |\partial\Omega|\|\varphi\|_{L^{\infty}(\partial\Omega)}^2$, we derive
\begin{eqnarray}
	\nu|<F,v_2>|
	&\leq& \nu\int_{\Omega}|F||v_2|dx\leq \nu\int_{\Omega}|F|\frac{|v_2|}{dist(x,\partial\Omega)}\cdot dist(x,\partial\Omega)dx\nonumber\\
	&\leq& \nu \Big(\int_{\Omega}\frac{|v_2|^2}{dist^2(x,\partial\Omega)}dx\Big)^{\frac{1}{2}}\Big(\int_{C'_1\varepsilon\leq dist(x,\partial\Omega)\leq C'_2\varepsilon}|F|^2dist^2(x,\partial\Omega)dx\Big)^{\frac{1}{2}}\nonumber\\
	&\leq& \nu\|v_2\|\times C'_2\varepsilon\|F\|_{L^2(\Omega)}\leq \nu\|v_2\|\times \frac{C}{\varepsilon^{\frac{1}{2}}}\|\varphi\|_{L^2(\partial\Omega)}\nonumber\\
	&\leq& \frac{\nu}{12}\|v_2\|^{2}+\frac{C\nu\|\varphi\|^2_{L^{\infty}(\partial\Omega)}|\partial\Omega|}{\varepsilon}.\label{g6}
\end{eqnarray}
(f) Using the Cauchy inequality and Young's inequality, we have
\begin{eqnarray}
	|<f,v_2>|&\leq& \|f\|_{V'}\|v_2\|\leq\frac{\nu}{12}\|v_2\|^{2}+\frac{C}{\nu}\|f\|^2_{V'}.\label{g6'}\end{eqnarray}
Combining \eqref{g2}--\eqref{g6'}, we conclude
\begin{eqnarray}&&\hspace{-1.6cm}\frac{d}{dt}|v_2|^2+\nu\|v_2\|^2\leq\Big(\frac{2C}{\nu}+\frac{4C\|\varphi\|^2_{L^{\infty}(\partial\Omega)} }{\nu}\Big)\rho_0(t) +\frac{4C\|\varphi\|^4_{L^{\infty}(\partial\Omega)}|\partial\Omega|\varepsilon}{\nu}\nonumber\\
	&&\hspace{2.6cm}
	 +\frac{2C\nu\|\varphi\|^2_{L^{\infty}(\partial\Omega)}|\partial\Omega|}{\varepsilon}+\frac{2C}{\nu}\|f\|^2_{V'}.\label{5.23}\end{eqnarray}
By Poincar\'{e}'s inequality, using the definition of fractal power of operators $A$, since $\lambda_j$ is increasing with respect to $j=m+1,\cdots$, for $v_2=\displaystyle{\sum_{j=m+1}}a_j\omega_j$ and $A^sv_2=\displaystyle{\sum_{j=m+1}}\lambda^sa_j\omega_j$, we have
\begin{eqnarray}
	 \|v_2\|^2&=&\|A^{\frac{1}{2}}v_2\|^2_H=\|\sum_{j=m+1}\lambda^{\frac{1}{2}}_ja_j\omega_j\|^2_H\geq \|\sum_{j=m+1}\lambda^{\frac{1}{2}}_{m+1}a_j\omega_j\|^2_H\nonumber\\
	&\geq& \lambda_{m+1}\|\sum_{j=m+1}a_j\omega_j\|^2_H
	=\lambda_{m+1}|v_2|^2.\end{eqnarray}
Thus, we derive from \eqref{5.23}
\begin{eqnarray}&&\hspace{-1.5cm}\frac{d}{dt}|v_2|^2+\nu\lambda_{m+1}|v_2|^2\nonumber\\
	 &&\hspace{-1.5cm}\leq\Big(\frac{2C}{\nu}+\frac{4C\|\varphi\|^2_{L^{\infty}(\partial\Omega)} }{\nu}\Big)\rho_0(t) +\frac{4C\|\varphi\|^4_{L^{\infty}(\partial\Omega)}|\partial\Omega|\varepsilon}{\nu}
	 +\frac{2C\nu\|\varphi\|^2_{L^{\infty}(\partial\Omega)}|\partial\Omega|}{\varepsilon}+\frac{2C}{\nu}\|g\|^2_{V'}\nonumber\\
	&&\hspace{-1.5cm}=K_1+\frac{2C}{\nu}\|f\|^2_{V'},\label{gg10}\end{eqnarray}
where $K_1=\Big(\frac{2C}{\nu}+\frac{4C\|\varphi\|^2_{L^{\infty}(\partial\Omega)} }{\nu}\Big)\rho_0(t) +\frac{4C\|\varphi\|^4_{L^{\infty}(\partial\Omega)}|\partial\Omega|\varepsilon}{\nu}
+\frac{2C\nu\|\varphi\|^2_{L^{\infty}(\partial\Omega)}|\partial\Omega|}{\varepsilon}$.\\
Applying the Gronwall inequality on $[\tau,t]$ to \eqref{gg10}, we conclude
\begin{eqnarray}
	\|v_2(t)\|_H^2&\leq& \|v_2(\tau)\|_H^2e^{-\nu\lambda_{m+1}(t-\tau)}+\int^t_{\tau}e^{-\nu\lambda_{m+1}(t-s)}|f(s)|^2_{V'}ds\nonumber\\
	 &&+(\frac{K_1}{\nu\lambda_{m+1}}-\frac{K_1e^{\nu\lambda_{m+1}\tau}}{\nu\lambda_{m+1}})\nonumber\\
	&\leq& e^{-\nu\lambda_{m+1}t}\rho^2_{\hat{D}}(t)e^{\nu\lambda_{m+1}\tau}+e^{-\nu\lambda_{m+1}t}e^{\nu\lambda_{m+1}\tau}\int^t_{\tau}|f(s)|^2_{V'}ds\nonumber\\
	 &&+(\frac{K_1}{\nu\lambda_{m+1}}-\frac{K_1e^{\nu\lambda_{m+1}\tau}}{\nu\lambda_{m+1}})\nonumber\\
	&=&I_1+I_2+I_3.\label{g7}
\end{eqnarray}
By the definition of the universe $\mathcal{D}_{\mu}$, using the technique in \eqref{gg7}, there exists a pullback time $\tau'_1<\tau_1$, then for any $\tau<\tau'_1\leq t$, it follows
\begin{eqnarray}
	I_1=e^{-\nu\lambda_{m+1}t}\rho^2_{\hat{D}}(t)e^{\nu\lambda_{m+1}\tau}\leq e^{-\nu\lambda_{m+1}t}\rho^2_{\hat{D}}(t)e^{\mu\tau}e^{(\nu\lambda_{m+1}-\mu)\tau}\leq \frac{\varepsilon}{3}\end{eqnarray}
as $m\rightarrow +\infty$ and $\tau\rightarrow -\infty$.\\
Since $f\in L^2_{loc}(\mathbb{R};V')$, using $\displaystyle{\lim_{m\rightarrow\infty}}\lambda_{m+1}=+\infty$ if we can choose $m$ large enough and letting $\tau\rightarrow -\infty$, i.e., there exists a pullback time $\tau_2\leq t$ and positive integer $m_0>0$, for any $\tau<\tau_2$ and $m>m_0$, we obtain
\begin{eqnarray}
	I_2&=&e^{-\nu\lambda_{m+1}t}e^{\nu\lambda_{m+1}\tau}\|f(s)\|^2_{L^2_{loc}(\mathbb{R};V')}< \frac{\varepsilon}{3}.\end{eqnarray}
Moreover, since $K_1$ is bounded, using $\displaystyle{\lim_{m\rightarrow\infty}}\lambda_{m+1}=+\infty$ if we can choose $m$ large enough, i.e., $m>m_0$ as above, letting $\tau\rightarrow -\infty$ (there exists a pullback time $\tau_3$, for any $\tau<\tau_3$), then we can achieve
\begin{eqnarray}
	 I_3=(\frac{K_1}{\nu\lambda_{m+1}}-\frac{K_1e^{\nu\lambda_{m+1}\tau}}{\nu\lambda_{m+1}})< \frac{\varepsilon}{3}.\label{g8}\end{eqnarray}
Combining  \eqref{g7}--\eqref{g8}, choosing a pullback time $$\hat{\tau}\leq min\{\tau_1,\tau_2,\tau_3,\tau'_1,\tau(t,\varepsilon_1)\},$$ such that for any $\tau<\hat{\tau}\leq t$, we conclude
\begin{eqnarray}
	&&\|(I-P)U(t,\tau)v_{\tau}\|^2_H=|v_2|^2\leq \varepsilon\end{eqnarray}
which implies the pullback $\mathcal{D}_{\mu}$-condition-(MWZ) holds for the norm-to-weak process $\{U(t,\tau)\}$, i.e., $U(t,\tau)$ is pullback $\mathcal{D}_{\mu}$-asymptotically compact. Hence, we complete the proof.$\hfill$$\Box$\\

 \noindent $\bullet$ {\bf Proof of Theorem \ref{th5.6}.}\\ Combining the norm-to-weak continuity of process, the existence of pullback $\mathcal{D}_{\mu}$-absorbing sets $\mathcal{B}$ and the pullback $\mathcal{D}_{\mu}$-asymptotic compactness of $\{U(t,\tau)\}$, using Theorem \ref{th8.35}, it is easy to prove the main result here. $\hfill$$\Box$\\

 \noindent $\bullet$ {\bf Proof of Theorem \ref{th5.7}.}\\ {\bf Inclusion closed:}\\
We begin by defining our tempered universe $\mathcal{D}_{\mu}$.
Let $R_{0}$ be the family
\begin{eqnarray}
	R_{0} = \big\{ \rho':\mathbb{R}\to [0,\infty) \, | \; e^{\mu\tau}  | \rho'(\tau)|^{2}  \to 0
	\;\; \mbox{as} \;\; \tau \to - \infty \big\} ,
\end{eqnarray}
where $\mu>0$ is defined in Lemma \ref{le5.3}, and the closed balls
$$
\overline{B_{t}(0,\rho'(t))}= \{ v_{\tau} \in H \, | \; \Vert v_{\tau} \Vert_{H} \le \rho'(t) \}.
$$
Then we can define a family of sets of as the form
\begin{eqnarray} \hspace{-0.355cm}\label{D-hat}
	\hat{D} = \big\{ \{D(t)\}_{t \in \mathbb{R}} \, | \;
	D(t) \neq \emptyset \;\; \mbox{and} \;\;
	D(t)\subset \overline{B_{t}(0,\rho_{\hat{D}}(t))} \;\;
	\mbox{with} \;\; \rho_{\hat{D}} \in R_0 \big\}.\nonumber
	\\
\end{eqnarray}
The class of all families of the form $\hat{D}$ is denoted by
\begin{equation} \label{universe}
	\mathcal{D}_{\mu}= \{\hat{D} \, | \; \hat{D} \;\; \mbox{satisfies} \;\; (\ref{D-hat}) \}.
\end{equation}
Clearly $\mathcal{D}_{\mu}$ is inclusion closed.

{\bf  We want to show $\hat{B}_0\in\mathcal{D}_{\mu}$:}\\
{\bf Case 1: $f\in L^2_{pb}(\mathbb{R};V')$}\\
Since $f\in L^2_{loc}(\mathbb{R},V')$ and satisfies $\int^t_{-\infty}e^{-\mu(t-s) }\|f\|^2_{V'}d s=\frac{1}{(1-e^{-\mu})}\|f\|^2_{L^2_{pb}(\mathbb{R};V')}$. If we want to show that $\rho^2_0(\tau)e^{\mu_0\tau}\rightarrow 0$ as $\tau\rightarrow -\infty$,
i.e.,
\begin{eqnarray}
	 &&\hspace{-2cm}\Big(1+\frac{4}{\nu\mu}\Big[\frac{C\nu^2}{\varepsilon}|\partial\Omega|\|\varphi\|^2_{L^{\infty}(\partial\Omega)}+C\varepsilon\ |\partial\Omega|\|\varphi\|_{L^{\infty}(\partial\Omega)}^2\Big]+\frac{4}{\nu\lambda_1}\int^t_{-\infty}e^{-\mu(t-s) }\|f\|^2_{V'}ds\Big)e^{\mu_0\tau}\rightarrow 0,
\end{eqnarray}
it is suffices to prove that
$\Big[\|f\|^2_{L^2_{loc}(\mathbb{R};V')}+\frac{4e^{-\mu n_0h}}{\nu\lambda_0(1-e^{-\mu n_0h})}\|f\|^2_{L^2_{pb}(\mathbb{R};V')}\Big]e^{\mu_0\tau}\rightarrow 0$ for every $\mu$, which is obviously.\\
{\bf Case 2: if $f(t,x)\in L^2_{upt}(\mathbb{R};V')$}\\
Since $f\in L^2_{loc}(\mathbb{R},V')$ and satisfies \begin{eqnarray}
\int^t_{-\infty}e^{-\mu(t-s) }\|f\|^2_{V'}d s&=&\int^t_{\tau}e^{-\mu(t-s) }\|f\|^2_{V'}d s+\int_{\tau}^{-\infty}e^{-\mu(t-s) }\|f\|^2_{V'}d s,\nonumber\\
\int^t_{\tau}e^{-\mu(t-s) }\|f\|^2_{V'}d s&\leq& C,\nonumber
\end{eqnarray} where $C>0$ is a constant, if we want to show that $\rho^2_0(\tau)e^{\mu_0\tau}\rightarrow 0$ as $\tau\rightarrow -\infty$,
i.e.,
\begin{eqnarray}
&&\hspace{-2cm}\Big(1+\frac{4}{\nu\mu}\Big[\frac{C\nu^2}{\varepsilon}|\partial\Omega|\|\varphi\|^2_{L^{\infty}(\partial\Omega)}+C\varepsilon\ |\partial\Omega|\|\varphi\|_{L^{\infty}(\partial\Omega)}^2\Big]\nonumber\\
&&\hspace{3cm}+\frac{4}{\nu\lambda_1}\int^t_{-\infty}e^{-\mu(\tau-s) }\|f\|^2_{V'}d s\Big)e^{\mu_0\tau}\rightarrow 0,
\end{eqnarray}
it is suffices to prove that
\begin{eqnarray}
\Big(\frac{4}{\nu\lambda_1}e^{-\mu\tau }\int^{\tau}_{-\infty}e^{\mu s }\|g\|^2_{V'}d s\Big)e^{\mu_0\tau}=\frac{4}{\nu\lambda_1}e^{(\mu_0-\mu)\tau }\int^{\tau}_{-\infty}e^{\mu s }\|f\|^2_{V'}d s\rightarrow 0,\nonumber
\end{eqnarray}
which needs that $\mu_0-\mu>0$, i.e., $0<\mu<\mu_0$.

By Theorem \ref{th8.35}, from the inclusion closed and $\hat{B}_0\in\mathcal{D}_{\mu}$, since the process is continuous, it is easy to achieve the result.$\hfill$$\Box$

\subsection{Proof of Theorem \ref{th5-f}: The fractal dimension of the pullback attractors}\label{5.5}
\ \ \ \ \ \ In this subsection, we shall use the abstract theory in \cite{clr} to estimate the dimension of pullback attractor $\mathcal{A}^H_{\mathcal{D}_{\mu}}$ in $H$. The first variation equation of \eqref{c1a} can be written as
\begin{eqnarray}
\left\{\begin{array}{ll} \frac{dU}{dt}+\nu AU+B(U,\psi)+B(\psi,U)+B(v,U)+B(U,v)=0,\\
U(\tau)=\xi\in H.\end{array}\right.\label{5-5-1}\end{eqnarray}
Multiplying \eqref{5-5-1} with $U$ and integrating over $\Omega$, using the property of trilinear operator, we can derive the estimate
\begin{equation}
|(B(\psi,U)+B(U,\psi),U)|\leq \frac{\nu}{4}\|U\|^2
\end{equation}
and the unique solution
$U\in L^2(\tau,T;V)\cap C(\tau,T;H)$ for any $T>\tau$.

\begin{Lemma}
	(1) There exists a bounded linear operator $\Lambda(t,s;v_0): H\rightarrow H$ such that
	\begin{equation}
	\sup_{\hat{v}_0, v_0\in \mathcal{A}_{\mathcal{D}_{\mu}}(s)}\sup_{\|\hat{v}_0-v_0\|\leq \varepsilon}\frac{\|U(t,s)\hat{v}_0-U(t,s)v_0-\Lambda(t,s;\hat{v}_0)(\hat{v}_0-v_0)\|}{\|\hat{v}_0-v_0\|}\rightarrow 0\label{5-5-9}
	\end{equation}
	 as $\varepsilon\rightarrow 0$.
	
	 (2) Moreover, $\Lambda(t,s;v_0)\xi=U(t)$ is the solution of the equation \eqref{5-5-1} and the operator $\Lambda(t,s;v_0)$ is compact for all $t\geq s$.
\end{Lemma}
{\bf Proof.} (1) Let $v(t)$ and $\hat{v}(t)$ be two solutions of the equation
\begin{equation}
\frac{dv}{dt}+\nu Av+B(v,v)+B(v,\psi)+B(\psi,v)+B(\psi,\psi)-\nu F(x)=Pf(t)
\end{equation}
with $v(s)=v_0$ and $\hat{v}(s)=\hat{v}_0$ respectively, and consider the solution $U(t)$ of problem \eqref{5-5-1} with $U(s)=v_0-\hat{v}_0$, hence, denoting $w=v-\hat{v}$, $\theta=v-\hat{v}-U$ satisfies the following problem
\begin{eqnarray}
\left\{\begin{array}{ll} \frac{d\theta}{dt}+\nu A\theta+B(\theta,\psi)+B(\psi,\theta)+B(v,\theta)+B(\theta,v)-B(w,w)=0,\\
\theta(s)=0.\end{array}\right.\label{5-5-2}\end{eqnarray}
Taking inner product of \eqref{5-5-2} by $\theta$, and using \eqref{q7}, we can derive
\begin{eqnarray}
\frac{1}{2}\frac{d}{dt}|\theta|^2+\nu\|\theta\|^2=|b(\theta,\psi,\theta)|+|b(\theta,v,\theta)|+|b(w,w,\theta)|. \label{5-5-3}
\end{eqnarray}
Using the similar technique in \eqref{c6}, choosing appropriate $\varepsilon$, we conclude
\begin{eqnarray}
|b(\theta,\psi,\theta)|\leq \frac{\nu}{3}\|\theta\|^2. \label{5-5-4}
\end{eqnarray}
By Ladyzhenskaya's and Young's inequalities, we obtain
\begin{eqnarray}
|b(\theta,v,\theta)|\leq C|\theta|\|\theta\|\|v\|\leq \frac{\nu}{3}\|\theta\|^2+\frac{C}{\nu}|\theta|^2\|v\|^2\label{5-5-5}
\end{eqnarray}
and
\begin{eqnarray}
|b(w,w,\theta)|\leq C|w|\|w\|\|\theta\|\leq \frac{\nu}{3}\|\theta\|^2+\frac{C}{\nu}|w|^2\|w\|^2.\label{5-5-6}
\end{eqnarray}
Combining \eqref{5-5-3}-\eqref{5-5-6}, it follows
\begin{eqnarray}
\frac{d}{dt}|\theta|^2\leq \frac{C}{\nu}|\theta|^2\|v\|^2+\frac{C}{\nu}|w|^2\|w\|^2.\label{5-5-7}
\end{eqnarray}
Using the same argument in Theorem 11.6 (see Carvalho, Langa and Robinson \cite{clr}), we conclude
\begin{eqnarray}
|\theta|^2\leq \Big[\frac{C}{\nu}e^{D(t,s)}\Big(1+D(t,s)e^{D(t,s)}\Big)\Big]|w(s)|^4\label{5-5-8}
\end{eqnarray}
with $D(t,s)=\frac{C}{\nu}\int^t_s\|v(r)\|dr$ and $w(s)=v_0-\hat{v}_0$, here we use the initial data $\theta(s)=0$. Since $D(t,s)$ is bounded, \eqref{5-5-8} implies that $|\theta|^2\leq C\epsilon^4$ and hence \eqref{5-5-9} is proved.

(2) {\bf Step 1:} Taking inner product of \eqref{5-5-1} with $U$, integrating over $\Omega$, using Lemma 2.3 and Hardy's inequality, choosing $C\epsilon\|\phi(x)\|_{L^{\infty}(\partial\Omega)}\leq \frac{\nu}{4}$, we derive
\begin{eqnarray}
\frac{1}{2}\frac{d}{dt}|U|^2+\nu\|U\|^2&\leq& |b(U,\psi,U)|+|b(U,v,U)|\nonumber\\
&\leq& \int_{\Omega}\frac{|U|^2}{dist^2(x,\partial\Omega)}|\nabla \psi|dist^2(x,\partial\Omega)dx+C|U|^{\frac{1}{2}}|\nabla v||U|^{\frac{1}{2}}\nonumber\\
&\leq& C\epsilon\|\phi(x)\|_{L^{\infty}(\partial\Omega)} \|U\|^2+\frac{\nu}{4}\|U\|^2+\frac{C}{\lambda_1}\|v\|^2\nonumber\\
&\leq& \frac{\nu}{2}\|U\|^2+\frac{C}{\lambda_1}\|v\|^2,\label{yan-2}
\end{eqnarray}
i.e.,
\begin{eqnarray}
\frac{d}{dt}|U|^2+\nu\lambda_1|U|^2&\leq& \frac{C}{\lambda_1}\|v\|^2.
\end{eqnarray}
Using Gronwall's inequality, we derive
\begin{eqnarray}
|U(t)|^2\leq e^{-\nu\lambda_1(t-s)}|\xi|^2+\frac{C}{\lambda_1}\|v\|^2_{L^2(\tau,T;V)},
\end{eqnarray}
since $v\in L^2(\tau,T;V)$, there exists an time $\tau(\|v\|^2_{L^2(\tau,T;V)},\xi)$, such that
$
|U(t,s;\xi)|^2\leq \rho^2_0
$ for all $\tau\leq \tau(\|v\|^2_{L^2(\tau,T;V)},\xi)$.\\
Integrating \eqref{yan-2} from $t$ to $t+1$, we derive
\begin{eqnarray}
|U(t+1)|^2+\nu\int^{t+1}_t\|U(s)\|^2ds&\leq& |U(t)|^2+\frac{C}{\lambda_1}\|v\|_{L^2(0,T;V)}\nonumber\\
&\leq& e^{-\nu\lambda_1(t-s)}|\xi|^2+\frac{2C}{\lambda_1}\|v\|^2_{L^2(\tau,T;V)}.
\end{eqnarray}

{\bf Step 2:} Multiplying \eqref{5-5-1} with $A^{\sigma}U$, integrating by parts over $\Omega$, we have
\begin{eqnarray}&&\hspace{-0.5cm}
\frac{1}{2}\frac{d}{dt}\int_{\Omega}|A^{\frac{\sigma}{2}}U|^2dx
+\nu\int_{\Omega}|A^{\frac{\sigma+1}{2}}U|dx\nonumber\\
&&\hspace{-0.5cm}\leq|b(v,U,A^{\sigma}U)|+|b(U,v,A^{\sigma}U)|+|b(U,\psi,A^{\sigma}U)|+|b(\psi,U,A^{\sigma}U)|.\label{yan28}\end{eqnarray}		
Next, we shall estimate every term on the right-hand side of \eqref{yan28} for $\sigma\in (0,\frac{1}{2}]$.\\
(a-1): By Lemmas \ref{le3.4} and \ref{le6.3}, using Cauchy's inequality and the Poincar\'{e} inequality, we derive
\begin{eqnarray}
| b(v,U,A^{\sigma}U)|&\leq& C|v|^{\frac{1}{2}}|A^{\frac{1}{2}} U|^{\frac{1}{2}}|A^{\sigma}U|^{\frac{1}{2}}\leq C|v|^{\frac{1}{2}}\frac{1}{\lambda_1^{\frac{\sigma+1}{4}}}|A^{\frac{\sigma+1}{2}} U|^{\frac{3}{2}}\nonumber\\
&\leq&  \frac{\nu}{8}|A^{\frac{\sigma+1}{2}}U|^2+\frac{C}{\nu\lambda^{2+\sigma}_1}\|v\|^2\label{yan-6}
\end{eqnarray}
and
\begin{eqnarray}
|b(U,v,A^{\sigma}U)|&\leq& C|U||\nabla v|^{\frac{1}{2}}|A^{\sigma}U|^{\frac{1}{2}}
\leq \frac{\nu}{8}|A^{\frac{\sigma+1}{2}}U|^2+\frac{C}{\nu\lambda^{3+\sigma}_1}\|v\|^2.\label{yan-7}
\end{eqnarray}

(a-2): By the Hardy inequality and Lemma \ref{le7.2}, we derive
\begin{eqnarray}
|b(U,\psi,A^{\sigma}U)|
&\leq&\int_{\Omega}|U||\nabla \psi||A^{\sigma}U|dx\nonumber\\
&\leq& C\|\varphi\|_{L^{\infty}(\partial\Omega)} |A^{\frac{1}{2}}U||A^{\sigma}U|\nonumber\\
&\leq& \frac{\nu}{8} |A^{\frac{\sigma+1}{2}}U|^2+C\|\varphi\|_{L^{\infty}(\partial\Omega)} \|U\|^2\label{yan-2a}
\end{eqnarray}
and
\begin{eqnarray}|b(\psi,U,A^{\sigma}U)|
\leq C_4\|\varphi\|_{L^{\infty}(\partial\Omega)}|A^{\frac{1}{2}} U| |A^{\sigma}U|\leq \frac{\nu}{8} |A^{\frac{\sigma+1}{2}}U|^2+C\|\varphi\|_{L^{\infty}(\partial\Omega)} \|U\|^2.\label{yan-3}
\end{eqnarray}

Hence combining \eqref{yan28}--\eqref{yan-3}, we have \begin{eqnarray}
&&\hspace{-1cm}\frac{d}{dt}|
A^{\frac{\sigma}{2}}U|^{2}+\nu\int_{\Omega}|A^{\frac{\sigma+1}{2}}U|dx\leq
\Big(\frac{C}{\nu\lambda^{2+\sigma}_1}+\frac{C}{\nu\lambda^{3+\sigma}_1}\Big)\|v\|^2+2C\|\varphi\|_{L^{\infty}(\partial\Omega)} \|U\|^2.\label{yan31}\end{eqnarray}
Integrating \eqref{yan31} from $r$ to $t+1$ with $t\leq r\leq t+1$, we obtain
\begin{eqnarray}
|A^{\frac{\sigma}{2}}U(t+1)|^2&\leq& \Big(\frac{C}{\nu\lambda^{2+\sigma}_1}+\frac{C}{\nu\lambda^{3+\sigma}_1}\Big) \|v\|^2_{L^2(\tau,T;V)}+2C\|\varphi\|_{L^{\infty}(\partial\Omega)} \int^{t+1}_t\|U(r)\|^2dr\nonumber\\
&&+|A^{\frac{\sigma}{2}}U(r)|^2\nonumber\\
&\leq&
\Big(\frac{C}{\nu\lambda^{2+\sigma}_1}+\frac{C}{\nu\lambda^{3+\sigma}_1}+\frac{2C\|\varphi\|_{L^{\infty}(\partial\Omega)}}{\nu\lambda_1}\Big)\|v\|_{L^2(\tau,T;V)}\nonumber\\
&&+2C\|\varphi\|_{L^{\infty}(\partial\Omega)}e^{-\nu\lambda_1(t-s)}|\xi|^2+|A^{\frac{\sigma}{2}}U(r)|^2.\label{yan-9}
\end{eqnarray}
Then integrating with respect to $r$ from $t$ to $t+1$, since  $D(A^{\frac{\sigma}{2}})$ with $\sigma\in [0,\frac{1}{2}]$ is compact in $V$, we conclude
\begin{eqnarray}
|A^{\frac{\sigma}{2}}U(t+1)|^2
&\leq&
\Big(\frac{C}{\nu\lambda^{2+\sigma}_1}+\frac{C}{\nu\lambda^{3+\sigma}_1}+\frac{2C\|\varphi\|_{L^{\infty}(\partial\Omega)}}{\nu\lambda_1}\Big)\|v\|^2_{L^2(\tau,T;V)}\nonumber\\
&&+2C\|\varphi\|_{L^{\infty}(\partial\Omega)}e^{-\nu\lambda_1(t-s)}|\xi|^2+\int^{t+1}_t|A^{\frac{1}{2}}U(r)|^2dr\nonumber\\
&\leq& e^{-\nu\lambda_1(t-s)}|\xi|^2+\frac{2C}{\lambda_1}\|v\|^2_{L^2(\tau,T;V)}\nonumber\\
&&+\Big(\frac{C}{\nu\lambda^{2+\sigma}_1}+\frac{C}{\nu\lambda^{3+\sigma}_1}+\frac{2C\|\varphi\|_{L^{\infty}(\partial\Omega)}}{\nu\lambda_1}\Big)\|v\|^2_{L^2(\tau,T;V)}\nonumber\\
&&+2C\|\varphi\|_{L^{\infty}(\partial\Omega)}e^{-\nu\lambda_1(t-s)}|\xi|^2\nonumber\\
&=:&\hat{\rho}^2_0\label{yan-10}
\end{eqnarray}
for all $\tau\leq \tau(\|v\|^2_{L^2(\tau,T;V)},\xi)$.

{\bf Step 3:} From Step 2, there exists a bounded set that pullback absorbs all bounded sets with initial data in $H$. Since $D(A^{\frac{\sigma}{2}})$ is compact in $H$, then we prove that the operator $\Lambda(t,s;u_0)$ is compact.   \hfill $\Box$\\

{\bf Proof of Theorem \ref{th5-f}.} {\bf Step 1:} From the existence of family of pullback attractors, for a fixed $\tau^*$, $\bigcup{\mathcal{A}}(t)$ is precompact in $H$.

{\bf Step 2:} For each $t\geq \tau$ and $v_0\in H$, we can define the linear operator $\Lambda(t,s;v_0)\cdot \xi=U(t)$, where $U(t)$ is the solution of \eqref{5-5-1}.

For each $\xi_1,\ \xi_2,\ \cdots, \xi_n\in H$, let $U_i(t)=\Lambda(t,s;v_0)\cdot \xi_i$ with $v_0\in H$, we define $F(U(t,\tau)v_0,t)=-\nu A-B(\psi,\cdot)-B(\cdot,\psi)-B(\cdot,\cdot)-B(\psi,\psi)+\bar{f}$, then $F(\cdot,t)$ is Gateaux differential in $V$ at $U(t,\tau)v_0$ which satisfies
\begin{eqnarray}
F'(U(t,\tau)v_0)U=-\nu A U-B(\psi,U)-B(U,\psi)-B(U(t,\tau)v_0,U)-B(U,U(t,\tau)v_0).
\end{eqnarray}
Hence $F'(U(t,\tau)v_0,t)\in \mathcal{L}(V,V')$ is a continuous linear operator for $(v,t)\in V\times\mathbb{R}$, which also means the problem
\begin{eqnarray}
	\left\{\begin{array}{ll}\frac{dU}{dt}=F'(U(t,\tau)v_0,t)U,\ \ v_0\in H,\\
U(\tau,x)=\xi\end{array}\right.\label{yan-11}
\end{eqnarray}
possesses a unique solution $U(t)=U(t,\tau;v_0,\xi)\in L^2(\tau,T;V)\cap C(\tau,T;H)$.

Let $U_1(s)=U_1(s,\tau;v_0,\xi_1), \ U_2(s)=U_2(s,\tau;v_0,\xi_2),\ \cdots,\ U_n(s)=U_n(s,\tau;v_0,\xi_n)$ be the solution of problem \eqref{yan-11} with initial data $U_i(\tau)=\xi_i (i=1,2,\cdots,n)$ respectively, $Q_n(s)$ denote the projection from $H$ to the space $span\{U_1(s),\ U_2(s),\ \cdots,\ U_n(s)\}$, then by Lemma 4.19 in \cite{clr}, it follows
\begin{eqnarray}
&&\|U_1(t)\wedge U_2(t)\wedge\cdots \wedge U_n(t)\|_{\Lambda^n(H)}\nonumber\\
&=&\|\xi_1\wedge \xi_2\wedge\cdots \wedge \xi_n\|_{\Lambda^n(H)}\exp\Big(\int^t_{\tau}Tr_n(F'(U(s,\tau)v_0,s)\circ Q_n(s)ds)\Big)\label{yan-12}
.\end{eqnarray}

Let $\{e_1(s),\ e_2(s),\ \cdots,\ e_n(s)\}$ be an orthonormal basis for $span\{U_1(s),\ U_2(s),\ \cdots,\ U_n(s)\}$, then
\begin{eqnarray}
Tr_n(F'(U(s,\tau)v_0,s)=\sup_{\xi_i\in H, |\xi_i|\leq 1, i\leq n}\Big(\sum^n_{i=1}\langle F'(U(s,\tau)v_0,s)e_i,e_i\rangle\Big).\label{yan-13}
\end{eqnarray}
Since $U_i(s)\in L^2(\tau,T;V)$, then $U_i(s)\in V$ for a.e. $s\geq \tau$, hence
$e_i(s)\in V$ for a.e. $s\geq \tau$ and $i=1,2,\cdots,n$.

Noting that $b(U(t,\tau)v_0,e_i(s),e_i(s))=0$, $b(\psi,e_i(s),e_i(s))=0$ and
\begin{eqnarray}
|b(e_i(s),\psi,e_i(s))|&\leq& \int_{\Omega}|e_i(s)||\nabla\psi||e_i(s)|ds\nonumber\\
&\leq& C\epsilon \|\phi\|_{L^{\infty}(\partial\Omega)}\|e_i(s)\|^2\leq \frac{\nu}{4}\|e_i(s)\|^2,
\end{eqnarray}
here we choose appropriate $\epsilon$ such that $C\epsilon \|\phi\|_{L^{\infty}(\partial\Omega)}\leq \frac{\nu}{4}$, then we derive
\begin{eqnarray}
Tr_n(F'(U(s,\tau)v_0,s)\circ Q_n(s)&=&\sum^n_{i=1}\langle F'(U(s,\tau)v_0,s)e_i(s),e_i(s)\rangle\nonumber\\
&=&\sum^n_{i=1}\Big(-\nu\|e_i(s)\|^2-b(\psi,e_i(s),e_i(s))-b(e_i(s),\psi,e_i(s))\nonumber\\
&&\hspace{1cm}-b(U(s,\tau)v_0,e_i(s),e_i(s))-b(e_i(s),U(s,\tau)v_0,e_i(s))\Big)\nonumber\\
&\leq& -\frac{3}{4}\nu\sum^n_{i=1}\|e_i(s)\|^2+\sum^n_{i=1}|b(e_i(s),U(s,\tau)v_0,e_i(s))|.\label{yan-14}
\end{eqnarray}
For the second term in \eqref{yan-14},  by the Lied-Thirring inequality $\int_{\Omega}\Big(\displaystyle{\sum^n_{i=1}}|e_i(s)|^2\Big)^2ds\leq C_1\displaystyle{\sum^n_{i=1}}\|e_i(s)\|^2$, we could proceed
\begin{eqnarray}
\sum^n_{i=1}|b(e_i(s),U(s,\tau)v_0,e_i(s))|&\leq& C\sum^n_{i=1}|e_i(s)|^{\frac{1}{2}}\|U(s,\tau)v_0\||e_i(s)|^{\frac{1}{2}}\nonumber\\
&\leq& \|U(s,\tau)v_0\| \Big[\int_{\Omega}\Big(\sum^n_{i=1}|e_i(s,x)|\Big)^2ds\Big]^{\frac{1}{2}} \nonumber\\
&\leq& \frac{C}{\nu}\|U(s,\tau)v_0\|^2+\frac{\nu}{4}\sum^n_{i=1}\|e_i(s)\|^2.
\end{eqnarray}
Using the variational principle and $\displaystyle{\sum^n_{i=1}\lambda_i}\geq \frac{\pi n^2}{|\Omega|}$ from \cite{ily}, taking the average, we obtain
\begin{eqnarray}
Tr_n(F'(U(s,\tau)v_0,s)\circ Q_n(s)
&\leq& -\frac{\nu}{2}\sum^n_{i=1}\|e_i(s)\|^2+\frac{C}{\nu} \|U(s,\tau)v_0\|^2\nonumber\\
&\leq& -\frac{\nu}{2}\sum^n_{i=1}\lambda_i+\frac{C}{\nu} \|U(s,\tau)v_0\|^2\nonumber\\
&\leq& -\frac{\pi\nu n^2}{2|\Omega|}+\frac{C}{\nu} \|U(s,\tau)v_0\|^2.
\end{eqnarray}
Defining
\begin{eqnarray}
&&q_n=\sup_{t\in\mathbb{R}}\sup_{v_0\in\mathcal{A}(t)}\Big(\frac{1}{T}\int^t_{t-T}Tr_n(F'(U(s,\tau)v_0,s)\Big)ds,\\
&&\hat{q}_n=\limsup_{T\rightarrow +\infty}q_n,
\end{eqnarray}
we derive
\begin{eqnarray}
q_n\leq -\frac{\pi\nu n^2}{2|\Omega|}+\frac{C}{\nu}\sup_{t\in\mathbb{R}}\sup_{v_0\in\mathcal{A}(t)}\Big(\frac{1}{T}\int^t_{t-T}\|U(s,\tau)v_0\|\Big)ds
\end{eqnarray}
and
\begin{eqnarray}
\hat{q}_n\leq -\frac{\pi\nu n^2}{2|\Omega|}+\frac{C}{\nu}q,
\end{eqnarray}
where $q=\displaystyle{\limsup_{T\rightarrow +\infty}}\sup_{t\in\mathbb{R}}\sup_{v_0\in\mathcal{A}(t)}\frac{1}{T}\int^t_{t-T}\|U(s,\tau)v_0\|^2ds$.\\
Using the technique as \eqref{nsv4}, we derive
\begin{eqnarray}
\frac{\nu}{2}\int^t_s\|v(r)\|^2dr\leq |v_0(s)|^2+C\|f(r)\|_{L^2(s,t;V')}+\frac{C}{\nu}\Big[\frac{\nu^2|\Omega|}{\epsilon}+\epsilon|\partial\Omega|\Big]\|\varphi\|^2_{L^{\infty}(\partial\Omega)}(t-s).\label{yan-15}
\end{eqnarray}
Setting $s=t-T$ in \eqref{yan-15}, it follows
\begin{eqnarray}
q\leq \frac{C}{\nu^2}\Big[\frac{\nu^2|\Omega|}{\epsilon}+\epsilon|\partial\Omega|\Big]\|\varphi\|^2_{L^{\infty}(\partial\Omega)}+\frac{C}{\nu}\lim_{T\rightarrow +\infty}\frac{1}{T}\int^t_{t-T}\|f(r)\|^2_{V'}dr.
\end{eqnarray}
Defining $M=\displaystyle{\lim_{T\rightarrow +\infty}}\frac{1}{T}\int^t_{t-T}\|f(r)\|^2_{V'}dr$, then
\begin{eqnarray}
\hat{q}_n\leq -\frac{\pi\nu n^2}{2|\Omega|}+\frac{C}{\nu^3}\Big[\frac{\nu^2|\Omega|}{\epsilon}+\epsilon|\partial\Omega|\Big]\|\varphi\|^2_{L^{\infty}(\partial\Omega)}+\frac{C}{\nu^2}M.
\end{eqnarray}

{\bf Case 1:} If $\frac{\pi\nu n^2}{2|\Omega|}>\frac{C}{\nu^3}\Big[\frac{\nu^2|\Omega|}{\epsilon}+\epsilon|\partial\Omega|\Big]\|\varphi\|^2_{L^{\infty}(\partial\Omega)}+\frac{C}{\nu^2}M$, then by Lemma 4.19 in \cite{clr}, we have $\mbox{dim}_B(\mathcal{A}(t))\leq n$.

{\bf Case 2:} Otherwise, by the theory in Temam \cite{te}, Carvalho, Langa and Robinson \cite{clr}, we see that the fractal and Hausdorff dimension of pullback attractors proceed as  $\mbox{dim}(\mathcal{A}(t))\leq \frac{C|\Omega|^{\frac{1}{2}}}{\nu} q^{\frac{1}{2}}+\hat{C}$.

Denoting $\hat{M}=\displaystyle{\limsup_{T\rightarrow +\infty}}\frac{1}{T}\int^t_{t-T}|f(r)|^2dr\rightarrow \|f\|^2_{L^{\infty}(-\infty,T^*;H)}$, then we derive
\begin{eqnarray}
\mbox{dim}(\mathcal{A})&\leq& \frac{C|\Omega|^{\frac{1}{2}}}{\nu^4}\Big[\frac{\nu|\Omega|^{\frac{1}{2}}}{\epsilon}+\sqrt{\epsilon}|\partial\Omega|^{\frac{1}{2}}\Big]\|\varphi\|_{L^{\infty}(\partial\Omega)}+\frac{C|\Omega|^{\frac{1}{2}}}{\nu^3}M+\hat{C}_3\nonumber\\
&\leq& \frac{C|\Omega|^{\frac{1}{2}}}{\nu^4}\Big[\frac{\nu|\Omega|^{\frac{1}{2}}}{\epsilon}+\sqrt{\epsilon}|\partial\Omega|^{\frac{1}{2}}\Big]\|\varphi\|_{L^{\infty}(\partial\Omega)}+\frac{C|\Omega|^{\frac{1}{2}}}{\nu^3}\|f\|_{L^{\infty}(-\infty,T^*;H)}+\hat{C}_3\nonumber\\
&\leq& \hat{C}_1\mbox{Re}+\hat{C}_2\mbox{G}+\hat{C}_3,
\end{eqnarray}
here $\mbox{G}=\frac{\|f\|_{L^{\infty}(-\infty,T^*;H)}}{\nu^2\lambda_1}$, $\mbox{Re}=\frac{\|\varphi\|_{L^{\infty}(\partial\Omega)}}{\nu \lambda_1^{1/2}}$. This means the proof is complete.   \hfill $\Box$

\subsection{\bf Upper semi-continuity theory of pullback attractors for non-autonomous systems }\label{sub8.5}
Considering the non-autonomous systems with perturbed external forces
\begin{equation}
\frac{\partial u}{\partial t}=\hat{A}_{\sigma}u+\varepsilon \hat{\sigma}(t,x),\label{yang-1}\end{equation} this section is to show the upper semi-continuity theory between pullback attractors ${\mathcal
	A}_{\varepsilon}=\{A_{\varepsilon}(t)\}_{t\in {\mathbb R}}$ with respect to the parameter $\varepsilon\in (0,\varepsilon_{0}]$  for the evolutionary process
$U_{\varepsilon}(\cdot,\cdot)$ and global
attractor $\mathcal {A}$ for \eqref{yang-1} with the cases $\varepsilon>0$ and
$\varepsilon=0$ respectively.

\begin{Definition}
	(See Carvalho, Langa and Robinson \cite{clr}) Let $X$ be a Banach space and $\Lambda$ be metric spaces, ${A_{\lambda}} (\lambda\in\Lambda)$ be a family of subsets of
	$X$. We say that the family of pullback attractors $A_{\lambda}$ is upper semi-continuous as $\lambda\rightarrow \lambda_0$ if
	\begin{eqnarray}
	\lim_{\lambda\rightarrow \lambda_0}dist_X(A_{\lambda},A_{\lambda_0})=0.
	\end{eqnarray}
\end{Definition}

\begin{Theorem}\label{th3.13}
	(See Carvalho, Langa and Robinson \cite{clr}) Suppose that for each $\varepsilon\in [0,1)$, $U_{\varepsilon}(\cdot,\cdot)$ is a family of processes such that
	
	(i) $U_{\varepsilon}(\cdot,\cdot)$ has a pullback attractor $\mathcal{A}_{\varepsilon}(\cdot)$ for each $\varepsilon\in [0,1)$.
	
	(ii) For some $t\in\mathbb{R}$ and any $T\geq 0$, for each bounded set $D$ in $X$, we have
	\begin{eqnarray}
	\sup_{s\in[0,T]}\sup_{u_0\in D}dist_X(U_{\varepsilon}(t+s,t)u_0,U_{0}(t+s,t)u_0)\rightarrow 0\ \ as\ \ \varepsilon\rightarrow 0.
	\end{eqnarray}
	
	(iii) There exist $\delta>0$ and $t_0\in\mathbb{R}$ such that  $\displaystyle{\bigcup_{\varepsilon\in(0,\delta)}}\bigcup_{s\leq t_0}\mathcal{A}_{\varepsilon}(s)$ is bounded.
	
	Then  ${\mathcal A}_{\varepsilon}$ are upper semi-continuous to $\mathcal {A}_0$ as $\varepsilon$ goes to $0$, i.e., \begin{equation}\lim\limits_{\varepsilon\rightarrow
		0}dis_{X}(A_{\varepsilon}(t),{\mathcal {A}_0})=0 \quad for~ any~ t\in
	{\mathbb R}.\label{a6-1}\end{equation}	
\end{Theorem}

Next, using non-compact measure and decomposition method respectively, we shall present some equivalent results to Theorem \ref{th3.13} for the non-autonomous
dissipative system to obtain the upper semi-continuity between
pullback attractors ${\mathcal A}_{\varepsilon}$ and global
attractor $\mathcal {A}_0$.
\begin{Theorem}\label{th3}
	(See Caraballo, Langa and Robinson \cite{clr01}) For each $\tau\leq t\in {\mathbb R}$ and $x\in X$, we
	assume \begin{equation}(H_{1})\quad\quad\quad\quad\quad
	~\lim\limits_{\varepsilon\rightarrow
		0}d_{X}(U_{\varepsilon}(t,t-\tau)x, S(\tau)x)=0,\label{a4a}\end{equation} holds uniformly on
	bounded sets of $X$.
	
	Assume that $(H_{1})$ holds, there exist pullback attractors ${\mathcal
		A}_{\varepsilon}=\{A_{\varepsilon}(t)\}_{t\in {\mathbb R}}$ for all $\varepsilon\in
	(0,\varepsilon_{0}]$. If
	there exists a compact set $K\subset X$, such that
	 \begin{equation}(H_{2})\quad\quad\quad\quad\quad~\lim\limits_{\varepsilon\rightarrow
		0}dis_{X}(A_{\varepsilon}(t), K)=0 \quad for~ any~ t\in {\mathbb
		R}.\label{a5}\end{equation} Then ${\mathcal A}_{\varepsilon}$ are upper semi-continuous to $\mathcal {A}_0$.
\end{Theorem}

\begin{Theorem}\label{th1}
	(See Wang and Qin \cite{wq}) Assume the family of sets ${\mathcal B}=\{B(t)\}_{t\in{\mathbb R}}$ is
	pullback absorbing for the process $U(\cdot,\cdot)$, ${\mathcal
		K}_{\varepsilon}=\{K_{\varepsilon}(t)\}_{t\in {\mathbb R}}$ is a
	family of compact sets in $X$ for each
	$\varepsilon\in (0,\varepsilon_{0}]$. Suppose the decomposition
	 $U_{\varepsilon}(\cdot,\cdot)=U_{1,\varepsilon}(\cdot,\cdot)+U_{2,\varepsilon}(\cdot,\cdot):{\mathbb
		R}\times{\mathbb R}\times X\rightarrow X$ satisfies
	
	(i) for any $t\in {\mathbb R}$ and $\varepsilon\in
	(0,\varepsilon_{0}]$,
	\begin{eqnarray} \parallel
	U_{1,\varepsilon}(t,t-\tau)x_{t-\tau}\parallel_{X}\leq\Phi(t,\tau),\quad
	\forall~ x_{t-\tau}\in B(t-\tau),\quad \tau>0,\end{eqnarray} where
	$\Phi(\cdot,\cdot):{\mathbb R}\times{\mathbb R}\rightarrow {\mathbb
		R}^{+}$ satisfies $\lim\limits_{\tau\rightarrow
		+\infty}\Phi(t,\tau)=0$ for each $t\in{\mathbb R}$,
	
	(ii)  for any $t\in{\mathbb R}$ and any $T\geq 0$, $\displaystyle{\bigcup_{0\leq
			\tau\leq T}}U_{2,\varepsilon}(t,t-\tau)B(t-\tau)$ is bounded, and for
	any $t\in {\mathbb R}$, there exists a time $T_{{\mathcal B}}(t)>0$,
	which is independent of $\varepsilon$, such that
	\begin{equation}U_{2,\varepsilon}(t,t-\tau)B(t-\tau)\subset
	K_{\varepsilon}(t),\quad \forall~\tau\geq T_{{\mathcal B}}(t),\quad
	\varepsilon\in (0,\varepsilon_{0}]\label{a7}\end{equation} and there exists a compact set
	$K\subset X$, such that \begin{equation}(H'_{2})\quad\quad\quad\quad\quad
	\lim\limits_{\varepsilon\rightarrow
		0}dist_{X}(K_{\varepsilon}(t),K)=0, \quad for ~ any~ t\in {\mathbb
		R}.\label{a8}\end{equation}
	
	Then (a) for each $\varepsilon\in (0,\varepsilon_{0}]$, the system \eqref{yang-1}
	possesses a family of pullback attractors ${\mathcal
		A}_{\varepsilon}=\{A_{\varepsilon}(t)\}_{t\in {\mathbb R}}$, (b)
	$(H_{2})$ holds and hence $\mathcal{A}_{\varepsilon}$ has the upper semi-continuity to $\mathcal{A}_0$.
\end{Theorem}

\begin{Remark}
	For the upper semi-continuity of attractors between the system \eqref{yang-1} with $\varepsilon>0$ and $\varepsilon=0$, the same initial data are necessary, i.e., every trajectory should begin at the same point.
\end{Remark}

\subsection{Proof of Theorem \ref{th6.1}: Upper semi-continuity of pullback attractors in $H$  for perturbed non-autonomous external forces}\label{sec6}
 First, let us give some inequalities and recall the preliminary theory of continuity of pullback attractors.
\begin{Lemma}\label{le7.2}
	For the fractal operator $A^{\sigma}$, we have the property
	
	(a) the embedding inequality:
	\begin{eqnarray}
	|A^{\frac{\sigma}{2}}u|^2\leq \frac{C}{\lambda_1}|A^{\frac{\sigma+1}{2}}u|^2, \ \ \ \forall\ u\in D(A^{\frac{\sigma+1}{2}}),\label{2-10}
	\end{eqnarray}

	(b) the generalized Hardy's inequality:
	 \begin{eqnarray}&&\int_{\Omega}\frac{|A^{\sigma}u|^{2}}{dist(x,\partial\Omega)}dx\leq C\int_{\Omega}|A^{\sigma+\frac{1}{4}}u|^{2}dx\leq \frac{C}{\lambda^{1-2\sigma}_1}\|A^{\frac{\sigma+1}{2}}u\|^2_H,\ \ 0<\sigma<\frac{1}{2},\label{1-1}\end{eqnarray}
	for all $u\in D(A^{\frac{\sigma+1}{2}})$,
	
	(c) the Gagliardo-Nirenberg inequality:
	\begin{eqnarray} &&|A^{\sigma/2}u|^2\leq
	C_{2}|A^{\sigma/4}u||A^{3\sigma/4}u|,\ \ \forall \ u\in D(A^{3\sigma/4}).\label{6-1}\end{eqnarray}
\end{Lemma}
{\bf Proof.} (a) Using \eqref{2}, applying the Poincar\'{e} inequality to the fractal operator $A^{\frac{\sigma}{2}}$, we can prove the result \eqref{2-10}.

(b) Since $0\leq \sigma< \frac{1}{2}$, noting the definition \eqref{2}, using the classical Hardy's inequality and  the Poincar\'{e} inequality, \eqref{1-1} holds.

(c) By \eqref{6}, substituting the derivative as $\sigma/2$, $\sigma/4$ and $3\sigma/4$, the inequality follows from $\frac{1}{2}-\frac{\sigma}{2}=\frac{1}{2}(\frac{1}{2}-\frac{\frac{\sigma}{2}}{2})+\frac{1}{2}(\frac{1}{2}-\frac{\frac{3\sigma}{2}}{2})$.         $\hfill$$\Box$\\

Next, we shall prove the upper semi-continuity of pullback attractors $\mathcal{A}_{\delta}(t)$ to the corresponding global attractor $\mathcal{A}$ for the case $\delta=0$ in Brown, Perry and Shen \cite{bps} as $\delta\rightarrow 0$ in $H$ of \eqref{a1-2} which is equivalent to \eqref{a1-2}.
\begin{Lemma} \label{le6.2} Denote
	$R_{\mu}=\{r(\xi)|\displaystyle{\lim_{\xi\rightarrow-\infty}}e^{\mu
		\xi}r^2(\xi)=0\}$ and ${\mathcal{D}}_{\eta}$ be the class of
	families $\hat{D}=\{D(t):t\in \mathbb{R}\}\subset {\mathcal{D}}(H)$ as the universe such that
	$D(t)\subset \bar{B}(0,r_{\hat{D}}(t))$,
	where $\bar{B}(0,r_{\hat{D}}(t))$ is a closed ball in $H$
	with radius $r_{\hat{D}}(t)$ at center zero.
	
	Assume $v_{\tau}\in H$, the external forces $\sigma_0(x)\in H$ and $\sigma(t)\in L^2(\tau,T;H)$ satisfying \eqref{ya-1}. Then for any fixed $t\in {\mathbb R}$
	and any bounded set $B\subset H$, there exists a time $T(B,t)>0$, such that
	\begin{eqnarray}
		\| U_{\delta}(t,t-\tau)v_{t-\tau}\|^{2}_{H}\leq
		R^2_{\delta}(t), \quad \forall~ \tau\geq T(B,t),
		~u_{t-\tau}\in B,\label{f19}\end{eqnarray} where
	 $R^2_{\delta}(t)=\frac{1}{\nu\lambda_1}\Big(\frac{C\|\varphi\|^4_{L^{\infty}(\partial\Omega)}|\partial\Omega|\varepsilon}{\nu}
+\frac{C\nu\|\varphi\|^2_{L^{\infty}(\partial\Omega)}|\partial\Omega|}{\varepsilon}\Big)+\frac{C}{\nu^2\lambda_1}|\sigma_0(x)|^2+\frac{2C\delta^2}{\nu\lambda_1}
	\int^t_{-\infty}e^{-\nu\lambda_1(t-\xi)}|\sigma(\xi)|^2d\xi$.
	
	Moreover, if we denote $B_{\delta}(t)=\{v_{\delta}\in H:
	|v_{\delta}|^2\leq R^2_{\delta}(t)\},$ then
	${\mathcal
		B}_{\delta}=\{B_{\delta}(t)\}_{t\in {\mathbb{R}}}\in \mathcal{D}_{\mu}$ is
	the family of pullback absorbing sets in $H$, i.e., $
	\lim\limits_{t\rightarrow -\infty}e^{\eta t}R_{\delta}(t)=0$
	for any $\delta>0$.
	
\end{Lemma}
{\bf Proof.}
For any fixed $t\in \mathbb{R}$ and any $\tau\in \mathbb{R}$, we denote \begin{eqnarray} v_{\delta}(r):=v(r;t-\tau,v_{\tau})=v_{\delta}(r-t+\tau,t-\tau,v_{\tau})\
	for\ r\geq t-\tau.\label{f20}\end{eqnarray}

Taking inner product of \eqref{a1-2} with $e^{\eta t}v_{\delta}$ in $H$ ($\eta$ will be determined later), using $(B(v_{\delta},v_{\delta}),v_{\delta})=0$ and $(B(\psi,v_{\delta}),v_{\delta})=0$, we obtain  \begin{eqnarray}
	&&\hspace{-1.5cm}\frac{d}{dt}\Big(e^{\eta t}|v_{\delta}|^2\Big)+2\nu e^{\eta
		t}\|v_{\delta}\|^2\non\\
	&&\hspace{-1.5cm}\leq \eta e^{\eta
		t}|v_{\delta}|^2+2e^{\eta t}|b(v_{\delta},\psi,v_{\delta})|+2e^{\eta t}|b(\psi,\psi,v_{\delta})|+2e^{\eta
		t}(\sigma_0(x)+\nu F(x)+\delta \sigma(t), v_{\delta})
	\label{f21}\end{eqnarray} which holds for all $v_{\delta}\in H$.

Using the same technique in \eqref{c6}, \eqref{ya-10} and \eqref{g6}, we can derive
\begin{eqnarray}
	&&|b(v_{\delta},\psi,v_{\delta})|\leq\int_{\Omega}|v_{\delta}||\nabla \psi||v_{\delta}|dx\leq C_2C_3C_4\varepsilon\|\varphi\|_{L^{\infty}(\partial\Omega)}\|v_{\delta}\|^{2}\leq\frac{\nu}{10}\|v_{\delta}\|^{2},\label{c6-1}\\
	&&|b(\psi,\psi,v_{\delta})|\leq\int_{\Omega}|\psi||\nabla \psi||v_{\delta}|dx\leq \frac{\nu}{10}\|v_{\delta}\|^{2}+\frac{C\|\varphi\|^4_{L^{\infty}(\partial\Omega)}|\partial\Omega|\varepsilon}{\nu},\label{c6-2}\\
	&&\nu|<F,v_{\delta}>|\leq \nu\int_{\Omega}|F||v_{\delta}|dx\leq \frac{\nu}{10}\|v_{\delta}\|^{2}+\frac{C\nu\|\varphi\|^2_{L^{\infty}(\partial\Omega)}|\partial\Omega|}{\varepsilon}.\label{c6-3}
\end{eqnarray}

By the Cauchy-Schwarz inequality and Young's inequality, we obtain
\begin{eqnarray}
	(\sigma_0(x), v_{\delta})\leq \|\sigma_0(x)\|_{V'}\|v_{\delta}\|_V\leq \frac{\nu}{10}\|v_{\delta}\|^{2}+\frac{C}{\nu}|\sigma_0(x)|^2
\end{eqnarray}
and
\begin{eqnarray}
	(\delta\sigma(t), v_{\delta})\leq \|\sigma(t)\|_{H}\|v_{\delta}\|_H\leq \frac{\nu}{10}\|v_{\delta}\|^{2}+\frac{C\delta^2}{\nu\lambda_1}|\sigma(t)|^2.
\end{eqnarray}

Then taking into account the Poincar\'{e} inequality and choosing $\eta=\frac{\nu\lambda_1}{2}$, we have \begin{eqnarray}
	&&\hspace{-1cm}\frac{d}{dt}\Big(e^{\eta t}|
	v_{\delta}(t)|^2\Big)+\nu\lambda_1 e^{\eta
		t}|v_{\delta}(t)|^2\leq e^{\eta
		 t}\Big(\frac{C\|\varphi\|^4_{L^{\infty}(\partial\Omega)}|\partial\Omega|\varepsilon}{\nu}+\frac{C\nu\|\varphi\|^2_{L^{\infty}(\partial\Omega)}|\partial\Omega|}{\varepsilon}\Big)\nonumber\\
	&&\hspace{5cm}+
	\frac{C}{\nu}e^{\eta
		t}|\sigma_0(x)|^2+\frac{Ce^{\eta
			t}\delta^2}{\nu\lambda_1}|\sigma(t)|^2,\label{f22}\end{eqnarray}
which implies for all $\tau\in \mathbb{R}$
\begin{eqnarray} &&\hspace{-1cm}|v_{\delta}(t)|^2\leq
	e^{-\nu\lambda_1(t-\tau)}| v_{\tau}|^2+\frac{1}{\nu\lambda_1}\Big(\frac{C\|\varphi\|^4_{L^{\infty}(\partial\Omega)}|\partial\Omega|\varepsilon}{\nu}+\frac{C\nu\|\varphi\|^2_{L^{\infty}(\partial\Omega)}|\partial\Omega|}{\varepsilon}\Big)\nonumber\\
	 &&\hspace{3cm}+\frac{C}{\nu^2\lambda_1}|\sigma_0(x)|^2+\frac{C\delta^2}{\nu\lambda_1}\int^t_{\tau}e^{-\nu\lambda_1(t-\xi)}|\sigma(\xi)|^2d\xi.
	\label{f23}\end{eqnarray}

Let $\hat{D}\in {\mathcal{D}}_{\eta}$ be given above, then for any $u_{0}\in D(\tau)$ and $t\geq\tau$, it follows \begin{eqnarray}
	&&|U_{\delta}(t,t-\tau)u_{t-\tau}|^2\leq e^{-\eta
		 (t-\tau)}r^2_{\hat{D}}+\frac{1}{\nu\lambda_1}\Big(\frac{C\|\varphi\|^4_{L^{\infty}(\partial\Omega)}|\partial\Omega|\varepsilon}{\nu}+\frac{C\nu\|\varphi\|^2_{L^{\infty}(\partial\Omega)}|\partial\Omega|}{\varepsilon}\Big)\nonumber\\
	 &&\hspace{4cm}+\frac{C}{\nu^2\lambda_1}|\sigma_0(x)|^2+\frac{C\delta^2}{\nu\lambda_1}\int^t_{-\infty}e^{-\nu\lambda_1(t-\xi)}|\sigma(\xi)|^2d\xi.\label{f24}\end{eqnarray}
Setting $e^{-\nu\lambda_1
	(t-\tau)}r^2_{\hat{D}}\leq
\frac{C\delta}{\nu\lambda_1}\int^t_{-\infty}e^{-\nu\lambda_1(t-\xi)}|\sigma(\xi)|^2d\xi$,
then for each fixed
$t\in \mathbb{R}$, we denote $R_{\delta}(t)$ as
\begin{eqnarray}
	 &&(R_{\delta}(t))^2=\frac{1}{\nu\lambda_1}\Big(\frac{C\|\varphi\|^4_{L^{\infty}(\partial\Omega)}|\partial\Omega|\varepsilon}{\nu}+\frac{C\nu\|\varphi\|^2_{L^{\infty}(\partial\Omega)}|\partial\Omega|}{\varepsilon}\Big)+\frac{C}{\nu^2\lambda_1}|\sigma_0(x)|^2\nonumber\\
	&&\hspace{2cm}+\frac{2C\delta^2}{\nu\lambda_1}
	 \int^t_{-\infty}e^{-\nu\lambda_1(t-\xi)}|\sigma(\xi)|^2d\xi.\label{f25}\end{eqnarray}
Considering the family of closed balls $\hat{B}_{\delta}$ for any fixed $t\geq \tau$ in
$H$ defined by \begin{eqnarray} B_{\delta}(t)=\{u_{\delta}\in H||u_{\delta}|^2\leq
	2R^2_{\delta}(t)\},\label{f26}\end{eqnarray} it is easy to check that
$\mathcal{B}_{\delta}(t)\in {\mathcal{D}}_{\eta}$ and hence
$\mathcal{B}_{\eta}(t)$ is the family of ${\mathcal{D}}_{\eta}$-pullback absorbing sets for the
process $\{u_{\delta}(t,t-\tau)\}$.$\hfill$$\Box$

\begin{Lemma} \label{le6.3} Let $R_{\delta}(t)$, $B_{\delta}(t)$ are defined
	in Lemma \ref{le6.2}, then for any $t\geq \tau\in {\mathbb R}$, the solution $v(t)=U_{1,\delta}(t,t-\tau)v(t-\tau)$ of
	\eqref{f17} satisfies, for all $\tau\in \mathbb{R}$ and
	$u_{t-\tau}\in B_{\delta}(t-\tau)$, \begin{eqnarray}
		&&|U_{1,\delta}(t,t-\tau)v_{t-\tau}|^{2}\leq
		e^{-\nu\lambda_1\tau}R^2_{\delta}(t-\tau),\label{a25}\\\
		&&\int^t_{t-\tau}\|v_1(s)\|^2ds\leq J_{\delta}(t),\label{a25-a}\end{eqnarray}  where $J_{\delta}(t)$ is dependent on $\tau, \ R^2_{\delta}(t-\tau),\ \nu$ and $\lambda_1$.
\end{Lemma}
{\bf Proof.}  Taking inner product of \eqref{f17} with $v_1$ in $H$, since $(B(v_1,v_1),v_1)=0$, we derive \begin{eqnarray} \frac{1}{2}\frac{d}{dt}|
	v_1(t)|^2+\nu\|v_1(t)\|^{2}\leq
	|b(v_1,\psi,v_1)|.\label{ya-3}\end{eqnarray}
Using the same techniques as in \eqref{c6}, we can get
\begin{eqnarray}
	&&|b(v_{1},\psi,v_{1})|\leq\int_{\Omega}|v_{1}||\nabla \psi||v_{1}|dx\leq\frac{\nu}{2}\|v_{1}\|^{2}.\end{eqnarray}
By the Poincar\'{e} inequality, it yields \begin{equation} \frac{d}{dt}|
	v_1(t)|^2+\nu\lambda_1|v_1(t)|^{2}\leq
	0.\label{a25-1}\end{equation}
Applying Gronwall's inequality to \eqref{a25-1} from $t-\tau$ to $t$, we get \begin{eqnarray}
	|U_{1,\delta}(t,t-\tau)v_{t-\tau}|^{2}&\leq&| v_{t-\tau}|^2e^{-\nu\lambda_1\tau}\leq e^{-\nu\lambda_1\tau}R^2_{\delta}(t-\tau)\rightarrow 0, \ as\ \tau\rightarrow +\infty,\label{a261}\end{eqnarray}
and \eqref{a25-a} is the direct result of \eqref{a261}, this therefore completes the proof. \hfill$\Box$

\begin{Lemma}\label{le6.4} Let the family of pullback absorbing sets ${\mathcal B}_{\delta}(t)=\{B_{\delta}(t)\}_{t\in {\mathbb R}}$ be given
	by Lemma \ref{le6.2}, then for any fixed $t\geq \tau\in {\mathbb R}$, there exist a time
	$T_{\delta}(t,{\mathcal B})>0$ and a function $I_{\delta}(t)>0$, such that
	the solution $U_{2,\delta}(t,\tau)v_{\tau}=w(t)$ of \eqref{f18} satisfies
	\begin{eqnarray}
		\|U_{2,\delta}(t,t-\tau)v_{t-\tau}\|^{2}_{D(A^{\frac{\sigma}{2}})}\leq
		I_{\delta}(t),\label{a27}\end{eqnarray} for all  $\tau\geq
	T_{\delta}(t,{\mathcal B}) $ and any $v_{t-\tau}\in
	B_{\delta}(t-\tau)$.
\end{Lemma}
{\bf Proof.} Multiplying \eqref{f18} with $A^{\sigma}v_2$ and integrating by parts over $\Omega$, we have
\begin{eqnarray}&&\hspace{-0.5cm}
	\frac{1}{2}\frac{d}{dt}\int_{\Omega}|A^{\frac{\sigma}{2}}v_2|^2dx
	+\nu\int_{\Omega}|A^{\frac{\sigma+1}{2}}v_2|dx\nonumber\\
	 &&\hspace{-0.5cm}\leq|b(v_1,v_2,A^{\sigma}v_2)|+|b(v_2,v_1,A^{\sigma}v_2)|+|b(v_2,v_2,A^{\sigma}v_2)|+|b(\psi,\psi,A^{\sigma}v_2)|\nonumber\\
	&&+|b(v_2,\psi,A^{\sigma}v_2)|+|b(\psi,v_2,A^{\sigma}v_2)|+|(P{(\sigma_0(x)+\nu \sigma_0(x))+\delta P\sigma(t,x)},A^{\sigma}v_2)|.\label{a28}\end{eqnarray}		 
Next, we shall estimate every term on the right-hand side of \eqref{a28} for $\sigma\in (0,\frac{1}{2}]$.\\
(a) By Lemmas \ref{le3.4} and \ref{le6.3}, using Cauchy's inequality and the Poincar\'{e} inequality, we derive
\begin{eqnarray}
	| b(v_1,v_2,A^{\sigma}v_2)|&\leq& C|v_1||A^{\frac{1}{2}} v_2|^{\frac{1}{2}}|A^{\sigma}v_2|_{L^4}\leq C|v_1|\frac{1}{\lambda_1^{\frac{\sigma}{4}}}|A^{\frac{\sigma+1}{2}} v_2|^{\frac{1}{2}}|A^{\sigma+\frac{1}{4}}v_2|\nonumber\\
	&\leq& \frac{C}{\lambda_1^{\frac{1-\sigma}{4}}}|v_1||A^{\frac{\sigma+1}{2}} v_2|^{\frac{3}{2}}\leq \frac{\nu}{9}|A^{\frac{\sigma+1}{2}}v_2|^2+\frac{C}{\nu\lambda^{1-\sigma}_1}|v_1|^4,\label{x-6}
\end{eqnarray}
\begin{eqnarray}
	|b(v_2,v_1,A^{\sigma}v_2)|&\leq& C|v_2||\nabla v_1|^{\frac{1}{2}}|A^{\sigma}v_2|^{\frac{1}{2}}\leq \frac{C}{\lambda_1}\|v_1\|^{\frac{1}{2}}|A^{\frac{\sigma+1}{2}} v_2|^{\frac{3}{2}}\nonumber\\
	&\leq& \frac{\nu}{9}|A^{\frac{\sigma+1}{2}}v_2|^2+\frac{C}{\nu\lambda^{4}_1}\|v_1\|^2\label{x-7}
\end{eqnarray}
and
\begin{eqnarray}
	|b(v_2,v_2,A^{\sigma}v_2)|
	&\leq&|v_2|^{\frac{1}{2}}|A^{\frac{1}{2}}v_2|^{\frac{1}{2}}|A^{\sigma}v_2|\leq \frac{C}{\lambda^{\frac{2-\sigma}{4}}_1}|A^{\frac{\sigma}{2}}v_2|^{\frac{1}{2}}|A^{\frac{\sigma+1}{2}}v_2|^{\frac{3}{2}} \nonumber\\
	&\leq& \frac{\nu}{9} |A^{\frac{\sigma+1}{2}}v_2|^2+\frac{C}{\nu\lambda^{2-\sigma}_1}|A^{\frac{\sigma}{2}}v_2|^2.\label{x-1}
\end{eqnarray}
(b) By the Hardy inequality and Lemma \ref{le7.2}, we derive
\begin{eqnarray}
	|b(v_2,\psi,A^{\sigma}v_2)|
	&\leq&\int_{\Omega}|v_2||\nabla \psi||A^{\sigma}v_2|dx\nonumber\\
	&\leq& C_4\|\varphi\|_{L^{\infty}(\partial\Omega)}\Big(\int_{\Omega}\frac{|v_2|^2}{dist^2(x,\partial\Omega)}dx\Big)^{\frac{1}{2}}
	\Big(\int_{\Omega}|A^{\sigma}v_2|^2dx\Big)^{\frac{1}{2}}\nonumber\\
	&\leq& C\|\varphi\|_{L^{\infty}(\partial\Omega)} |A^{\frac{1}{2}}v_2||A^{\sigma}v_2|\nonumber\\
	&\leq& C\|\varphi\|_{L^{\infty}(\partial\Omega)} |A^{\frac{1}{2}}v_2| |A^{\frac{\sigma}{2}}v_2|^{1-\sigma}|A^{\frac{\sigma+1}{2}}v_2|^{\sigma}\nonumber\\
	&\leq& \frac{C\|\varphi\|_{L^{\infty}(\partial\Omega)}}{\lambda_1^{\frac{\sigma}{2}}} |A^{\frac{\sigma}{2}}v_2|^{1-\sigma}|A^{\frac{\sigma+1}{2}}v_2|^{1+\sigma}\nonumber\\
	&\leq& \frac{\nu}{9} |A^{\frac{\sigma+1}{2}}v_2|^2+\frac{C\|\varphi\|^{\frac{2}{1-\sigma}}_{L^{\infty}(\partial\Omega)}}{\nu\lambda_1^{\frac{\sigma}{1-\sigma}}} |A^{\frac{\sigma}{2}}v_2|^2,\label{x-2}
\end{eqnarray}
\begin{eqnarray}|b(\psi,v_2,A^{\sigma}v_2)|
	\leq C_4\|\varphi\|_{L^{\infty}(\partial\Omega)}|A^{\frac{1}{2}} v_2| |A^{\sigma}v_2|\leq \frac{\nu}{9} |A^{\frac{\sigma+1}{2}}v_2|^2+\frac{C\|\varphi\|^{\frac{2}{1-\sigma}}_{L^{\infty}(\partial\Omega)}}{\nu\lambda_1^{\frac{\sigma}{1-\sigma}}} |A^{\frac{\sigma}{2}}v_2|^2\label{x-3}
\end{eqnarray}
and
\begin{eqnarray}
	|b(\psi,\psi,A^{\sigma}v_2)|
	&\leq&C\|\varphi\|_{L^{\infty}(\partial\Omega)}\int_{\Omega}|\nabla \psi||A^{\sigma}v_2|dx\nonumber\\
	&\leq& C\|\varphi\|_{L^{\infty}(\partial\Omega)}
	\Big(\int_{\Omega}|\nabla \psi|^2dist(x,\partial\Omega)dx\Big)^{\frac{1}{2}}
	 \Big(\int_{\Omega}\frac{|A^{\sigma}v_2|^2}{dist(x,\partial\Omega)}dx\Big)^{\frac{1}{2}}
	\nonumber\\
	&\leq& C\|\varphi\|_{L^{\infty}(\partial\Omega)}\|\varphi\|_{L^2(\partial\Omega)} |A^{\sigma+\frac{1}{4}}v_2|\nonumber\\
	&\leq& \frac{C}{(\lambda_1)^{\frac{2\sigma-1}{4}}}\|\varphi\|_{L^{\infty}(\partial\Omega)} \|\varphi\|_{L^2(\partial\Omega)} |A^{\frac{\sigma+1}{2}}v_2|\nonumber\\
	&\leq& \frac{\nu}{9} |A^{\frac{\sigma+1}{2}}v_2|^2+\frac{C}{\nu\lambda_1^{\sigma-\frac{1}{2}}}|\partial\Omega|\|\varphi\|^4_{L^{\infty}(\partial\Omega)}.\label{x-4}
\end{eqnarray}
(c) By the Young inequality and Lemma \ref{le7.2}, we obtain
\begin{eqnarray}
	|(P\nu F, A^{\sigma}v_2)|&\leq& \nu\int_{\Omega}|F|A^{\sigma}v_2|dx\leq C\nu \sqrt{\varepsilon}\int_{\Omega}|F| \frac{|A^{\sigma}v_2|}{dist^{\frac{1}{2}}(x,\partial\Omega)}dx\nonumber\\
	&\leq& C\nu \sqrt{\varepsilon}|F||A^{\sigma+\frac{1}{4}}v_2|
	\leq \frac{C\nu \sqrt{\varepsilon}}{\lambda^{\frac{1-2\sigma}{4}}_1}|F||A^{\frac{\sigma+1}{2}}v_2|
	\nonumber\\
	&\leq& \frac{C\nu \varepsilon}{\lambda_1^{\frac{1-2\sigma}{2}}}|F|^2+\frac{\nu}{16}|A^{\frac{\sigma+1}{2}}v_2|^2\nonumber\\
	&\leq& \frac{C\nu}{\lambda_1^{\frac{1-2\sigma}{2}} \varepsilon}|\partial\Omega|\|\varphi\|_{L^{\infty}(\partial\Omega)}^2+\frac{\nu}{9}|A^{\frac{\sigma+1}{2}}v_2|^2,\label{x-5}
\end{eqnarray}
\begin{eqnarray}
	(P \sigma_0(x),A^{\sigma}v_2)\leq |\sigma_0(x)| |A^{\sigma} v_2|\leq \frac{1}{\lambda_1^{\frac{1-\sigma}{2}}}|\sigma_0(x)| |A^{\frac{1+\sigma}{2}} v_2|\leq \frac{\nu}{9}|
	A^{\frac{\sigma+1}{2}}w(t)|^2+\frac{C}{\nu\lambda_1^{1-\sigma}}|
	\sigma_0(x)|^{2}\label{a30-a}
\end{eqnarray}
and
\begin{eqnarray}
|(\delta P\sigma(t,x),A^{\sigma}v_2)|\leq  \frac{\nu}{9}|
A^{\frac{\sigma+1}{2}}w(t)|^2+\frac{C\delta^2}{\nu\lambda_1^{1-\sigma}}|
\sigma(t)|^{2}.\label{a30}
\end{eqnarray}
Hence combining \eqref{a28}--\eqref{a30}, we derive \begin{eqnarray}
	&&\hspace{-1cm}\frac{d}{dt}|
	A^{\frac{\sigma}{2}}v_2|^{2}\leq
	 \Big(\frac{2C\|\varphi\|^{\frac{2}{1-\sigma}}_{L^{\infty}(\partial\Omega)}}{\nu\lambda_1^{\frac{\sigma}{1-\sigma}}} +\frac{C}{\nu\lambda^{2-\sigma}_1}\Big)|A^{\frac{\sigma}{2}}v_2|^2\nonumber\\
	&&\hspace{1.5cm}+\frac{C\nu}{\lambda_1^{\frac{1-2\sigma}{2}} \varepsilon}|\partial\Omega|\|\varphi\|_{L^{\infty}(\partial\Omega)}^2+\frac{C}{\nu\lambda_1^{\sigma-\frac{1}{2}}}|\partial\Omega|\|\varphi\|^4_{L^{\infty}(\partial\Omega)}\nonumber\\
	 &&\hspace{1.5cm}+\frac{C}{\nu\lambda^{1-\sigma}_1}|v_1|^4+\frac{C\delta^2}{\nu\lambda_1^{1-\sigma}}|
	\sigma(t)|^{2}+\frac{C}{\nu\lambda_1^{1-\sigma}}|
	\sigma_0(x)|^{2}+\frac{C}{\nu\lambda^{4}_1}\|v_1\|^2.\label{a31}\end{eqnarray}
Applying the Gronwall inequality to \eqref{a31} from $t-\tau$ to $t$, and using Lemma \ref{le6.2}, we conclude that for all $t\geq\tau$ \begin{eqnarray}|
	A^{\frac{\sigma}{2}}v_2(t)|^{2}\leq
	I_{\delta}(t)=I_{\delta}(t,\tau,\nu,\lambda_1,R_{\delta}(t-\tau),J_{\delta}, \|\varphi\|_{L^{\infty}(\partial\Omega)},|f|, \|\sigma\|_{L^2_{loc}(\mathbb{R};H)})\leq C.\label{a32}\end{eqnarray}
	 Which completes the proof.\hfill $\Box$

\begin{Lemma} \label{le6.5}  Assume $v_{\tau}\in H$, the external forces $\sigma_0(x)\in H$ and $\sigma(t)\in L^2(\tau,T;H)$ satisfies \eqref{ya-1} for any fixed $t\geq\tau\in {\mathbb R}$. If $v_{\tau}$ belongs to some bounded set, then the solution
	$v_{\delta}(t)=U_{\delta}(t,t-\tau)v_{\tau}$ of perturbed non-autonomous problem \eqref{a1-1}
	has upper semi-continuity property to the solution
	$v(t)=S(t)v_{\tau}$ of the autonomous problem \eqref{a1-1} with
	$\delta=0$ uniformly in $H$ as $\delta\rightarrow 0^{+}$, i.e., \begin{eqnarray} \lim\limits_{\delta\rightarrow
			0^{+}}\sup\limits_{v_{\tau}\in B}\|v_{\delta}(t)-v(t)\|_{H}=0,\label{a33}\end{eqnarray}
	where $B$ is an arbitrary bounded subset in $H$.
\end{Lemma}
{\bf Proof.} Since $v_{\delta}(t)$ and $v(t)$ satisfy the following non-autonomous and autonomous problem respectively,
\begin{eqnarray} \left\{\begin{array}{ll} (v_{\delta}(t))_t+\nu
		 Av_{\delta}+B(v_{\delta},v_{\delta})+B(v_{\delta},\psi)+B(\psi,v_{\delta})+B(\psi,\psi)\\
		\hspace{7cm}=P{(\sigma_0(x)+\nu \sigma_0(x))+\delta P\sigma(t,x)},\\
		 \mbox{div}\  v_{\delta}=0,\\
		v_{\delta}|_{\partial\Omega}=0, \\
		v_{\delta}|_{t=\tau}=v_{\tau}\\
	\end{array}\right.\label{a35-1}\end{eqnarray} and
\begin{eqnarray} \left\{\begin{array}{ll} v_t+\nu Av+B(v,v)+B(v,\psi)+B(\psi,v)+B(\psi,\psi)=P(\sigma_0(x)+\nu \sigma_0(x)),\\
		 \mbox{div}\  v=0,\\
		v|_{\partial\Omega}=0, \\
		v|_{t=\tau}=v_{\tau},\\
	\end{array}\right.\label{a35-2}\end{eqnarray}
then we can verify that $y^{\delta}(t)=v_{\delta}(t)-v(t)$ satisfies the problem
\begin{eqnarray} \left\{\begin{array}{ll} (y^{\delta})_t+\nu
		 Ay^{\delta}=-B(y^{\delta},\psi)-B(\psi,y^{\delta})-B(v_{\delta},v_{\delta})+B(v,v)+\delta P\sigma(t,x),\\
		 \mbox{div}\  y^{\delta}=0,\\
		y_{\delta}|_{\partial\Omega}=0, \\
		y_{\delta}|_{t=\tau}=0.\\
	\end{array}\right.\label{a35}\end{eqnarray}  Multiplying \eqref{a35} with $y^{\delta}(t)$, integrating by parts and using the property of trilinear operator $b(\cdot,\cdot,\cdot)$, we have
\begin{eqnarray}
	&&\frac{1}{2}\frac{d}{dt}|y^{\delta}|^{2}+\nu \|
	y^{\delta}\|^2\nonumber\\
	&=&\langle B(v,v)-B(y^{\delta},\psi)-B(v_{\delta},v_{\delta}),y^{\delta} \rangle+\langle\delta
	\sigma(t),y^{\delta}\rangle\nonumber\\
	&\leq&|b(y_{\delta},\psi,y_{\delta})|+|\langle
	B(u,u)-B(u_{\delta},u_{\delta}),y^{\delta})\rangle|+\frac{\nu}{3}\|
	y^{\delta}(t)\|^2+\frac{C\delta^2}{\nu}|
	\sigma(t)|^{2}.\label{ya-0}\end{eqnarray}
Using the same technique in \eqref{c6} and choosing $\delta$ small enough, we can derive
\begin{eqnarray}
	&&|b(y_{\delta},\psi,y_{\delta})|\leq\int_{\Omega}|y_{\delta}||\nabla \psi||y_{\delta}|dx\leq C_2C_3C_4\delta\|\varphi\|_{L^{\infty}(\partial\Omega)}\|y_{\delta}\|^{2}\leq\frac{\nu}{3}\|y_{\delta}\|^{2}.\end{eqnarray}
From the Young inequality and $b(v,y^{\delta}, y^{\delta})=0$, we obtain \begin{eqnarray}|\langle
	B(v,v)-B(v_{\delta},v_{\delta}),y^{\delta})\rangle|
	= |
	b(y^{\delta},v_{\delta},y^{\delta})|
	\leq C|y^{\delta}|\|v_{\delta}\|\|y^{\delta}\|\leq\frac{\nu}{3}\|
	y^{\delta}\|^2+\frac{C}{\nu}|y^{\delta}|^2\|v_{\delta}\|^2.\label{ya-01}
\end{eqnarray}
Hence, combining \eqref{ya-0}-\eqref{ya-01}, it yields
\begin{eqnarray}&&\frac{d}{dt}|y^{\delta}|^{2}\leq \frac{C}{\nu}|y^{\delta}|^{2}\|v_{\delta}\|^{2}+\frac{C\delta^2}{\nu}|
	\sigma(t)|^{2}.\label{ya-6}\end{eqnarray}
From Theorem \ref{th4.1}, Lemmas \ref{le6.3}--\ref{le6.4} and \eqref{ya-1}, since $\sigma(t)\in L^{2}_{loc}({\mathbb R},H)$,
by the Gronwall inequality
to \eqref{ya-6}, we conclude
\begin{eqnarray}|y^{\delta}|^{2}
	\leq
	 \frac{C\delta^2}{\nu}e^{\frac{C}{\nu}\|u\|^2_{L^2(\tau,T;V)}}\int_{t-\tau}^{t}|\sigma(s)|^{2}ds\leq \frac{C\delta^2}{\nu_0}e^{\frac{C}{\nu}\|u\|^2_{L^2(\tau,T;V)}} \|\sigma\|^2_{L^{2}_{loc}({\mathbb R},H)}\leq C_0\delta^2\rightarrow
	0\end{eqnarray} as $\delta\rightarrow 0^+$, which means \eqref{a33} holds. The proof is completed.
\hfill$\Box$
\\

\noindent $\bullet$ {\bf Proof of Theorem \ref{th6.1}.} Using the theory in Section \ref{sub8.5}, firstly, we decompose equation \eqref{a1-1} as a linear equation with non-homogeneous initial data and a nonlinear equation with homogeneous initial data, i.e.,
\begin{equation}\left\{
	\begin{array}{lll}
		(v_1)_t+\nu
		Av_1=-B(v_1,v_1)-B(\psi,v_1)-B(v_1,\psi),\ &\\
		 \mbox{div}\  v_1=0,\ &\\
		v_1|_{\partial\Omega}=0,\ &\\
		v_1(\tau,x)=v_{\tau}(x)=u_{\tau}(x)+\psi(x),\ & \end{array}\right.\label{f17}\end{equation}
and note that \begin{equation}\left\{
	\begin{array}{lll}
		(v_2)_t+\nu
		Av_2=-B(v_1,v_2)-B(v_2,v_1)-B(v_2,v_2)-B(v_2,\psi)\\
		\hspace{3.2cm}-B(\psi,v_2)-B(\psi,\psi)+P{(\sigma_0(x)+\nu \sigma_0(x))+\delta P\sigma(t,x)},\ &\\
		 \mbox{div}\  v_2=0,\ &\\
		v_2|_{\partial\Omega}=0,\ &\\
		v_2(\tau,x)=0,\ & \end{array}\right.\label{f18}\end{equation}
the solution
$u_{\delta}(t)=U_{\delta}(t,\tau)u_{\tau}$ of \eqref{a1-2} with
initial data $v_{\tau}\in H$ can be decomposed as  \begin{eqnarray}
	 v_{\delta}=U_{\delta}(t,\tau)v_{\tau}=U_{1,\delta}(t,\tau)v_{\tau}+U_{2,\delta}(t,\tau)v_{\tau},\end{eqnarray}
where \begin{eqnarray} U_{1,\delta}(t,\tau)v_{\tau}=v_1(t),\ \
	U_{2,\delta}(t,\tau)v_{\tau}=v_2(t),\label{f16}\end{eqnarray} solve the
problems \eqref{f17} and \eqref{f18} respectively.

Secondly, under the existence of pullback absorbing set in Lemma \ref{le6.2}, we should estimate the solutions $v_1$ with $H$-norm small enough of \eqref{f17} and $v_2$ with uniform boundedness in more regular space of \eqref{f18} by Lemmas \ref{le6.3} and \ref{le6.4} respectively.

Finally, by the continuous theory of pullback attractors, we shall prove the upper semi-continuity of pullback attractors in $H$ by Lemma \ref{le6.5}.  \hfill $\Box$

\subsection{Proof of Theorems \ref{th7.11} and \ref{th7.12}: Regularity of pullback attractors}\label{sec7}
\noindent $\bullet$ {\bf Norm-to-weak continuity of the process in $D(A^{\frac{\sigma}{2}})$ $(0\leq \frac{\sigma}{2}\leq\frac{1}{4} )$}\label{sub7.1}\\
\ \ \ \ \ \ In this section, we only need to prove the generalization of the norm-to-weak continuous process which will be used for pullback attractors in the next section.

{\bf Proof of Theorem \ref{th7.6}.} Firstly, by the definition of regular solutions, using the Galerkin approximated technique and compact argument again,  noting that $v=u-\psi$, we need to estimate the norms of $v$ in some more regular norm for the non-autonomous problem, i.e., the asymptotic regularity of $v$.
For the case $\sigma=\frac{1}{2}$, for the boundedness of $v$ in $D(A^{\frac{1}{4}})$, we can refer to \cite{bps}. In fact, if $0\leq \frac{\sigma}{2}<\frac{1}{4} $, we can use interpolation inequality to achieve the bounded of $v$ in $D(A^{\frac{\sigma}{2}})$, here we omit the detailed proof.

Next, since the systems \eqref{a1} and \eqref{c1a} are equivalent, we only need to verify the norm-to-weak continuity of $u(t,x)$, and hence $v(t,x)$ yields.
Let $u_1(\cdot)$ and $u_2(\cdot)$ be two solutions to problem \eqref{a1} and with corresponding initial data
$u^1_{\tau}$ and $u^2_{\tau}$ respectively and background flow functions $\psi_1$ and $\psi_2$, take $w=u_1-u_2$. Then $w$ satisfies the abstract problem:
\begin{eqnarray}\left\{\begin{array}{ll} \frac{dw}{dt}+\nu A
w+B(u_1,u_1)-B(u_2,u_2)=0,\\
 \mbox{div}\  w=0,\\
w(t,x)|_{\partial\Omega}=0,\ (x,t)\in \partial\Omega_{\tau},\\
w(\tau,x)=w_{\tau}=u^1_{\tau}(x)-u^2_{\tau}(x).\end{array}\right.\label{q6-1}\end{eqnarray}

Let $\omega\in C^{\infty}_{0}(\Omega)$, $ \mbox{div}\ \omega=0$, from the condition (ii) in Definition \ref{de7.5}, we can derive the weak formulation
\begin{eqnarray}\frac{d}{dt}(w,A^{\sigma}\omega)-\nu<w,A^{\sigma+1} \omega>=b(u_2,A^{\sigma}\omega,u_2)-b(u_1,A^{\sigma}\omega,u_1).\label{y1-1-1}\end{eqnarray}

Let $\omega=u_1-u_2$ in \eqref{y1-1-1}, we have
\begin{eqnarray}\frac{1}{2}\frac{d}{dt}|A^{\frac{\sigma}{2}}w|^2+\nu|A^{\frac{\sigma+1}{2}}w|^2&\leq& |b(w,A^{\sigma}w,u_2)|+|b(u_1,A^{\sigma}w,w)|.\label{gu5-1}\end{eqnarray}
Since
\begin{eqnarray}
|b(w,A^{\sigma}w,u_2)|&\leq& \int_{\Omega}|w||\nabla u_2||A^{\sigma}w|dx\leq C|\nabla u_2|^{\frac{1}{2}}|w|^{\frac{1}{2}}|A^{\frac{\sigma+1}{2}}w|\nonumber\\
&\leq& C\| u_2\|_{V}|A^{\frac{\sigma}{2}}w|+\frac{\nu}{4}|A^{\frac{\sigma+1}{2}}w|^2,\label{yang-3}
\end{eqnarray}
and
\begin{eqnarray}
|b(u_1,A^{\sigma}w,w)|&\leq& \int_{\Omega}|w||\nabla w||A^{\sigma}w|dx\leq C| u_1|\|w\|^{\frac{1}{2}}_{V}|A^{\frac{\sigma+1}{2}}w|^{\frac{1}{2}}\nonumber\\
&\leq& C| u_1||A^{\frac{\sigma+1}{2}}w|^{\frac{1}{2}}|A^{\frac{\sigma+1}{2}}w|^{\frac{1}{2}}\nonumber\\
&\leq& C| u_1|^2+\frac{\nu}{4}|A^{\frac{\sigma+1}{2}}w|^2,\label{yang-4}
\end{eqnarray}
we have
\begin{eqnarray}\frac{d}{dt}|A^{\frac{\sigma}{2}}w|^2+\nu|A^{\frac{\sigma+1}{2}}w|^2\leq C_{\lambda_1}|A^{\frac{\sigma}{2}}w|^2+ C(\| u_2\|^2_{V}+\| u_1\|^2_{V}).\label{gu5-1-1} \end{eqnarray}
Neglecting the second term on the left-hand side, since the index $\sigma\in [0,\frac{1}{2}]$, the trilinear operator can not imply appropriate $D(A^{\frac{\sigma}{2}})$ and $D(A^{\frac{\sigma+1}{2}})$ estimates on the right-hand side of \eqref{yang-3} and \eqref{yang-4}, we can only derive
\begin{eqnarray}
&&\hspace{-1.4cm}\|u_1(s)-u_2(s)\|^2_{D(A^{\frac{\sigma}{2}})}\leq \Big[\|u^1_0-u^2_0\|^2_{ D(A^{\frac{\sigma}{2}})}+C\Big(\| u_2\|^2_{L^2(\tau,T;V)}+\| u_1\|^2_{L^2(\tau,T;V)}\Big)\Big] e^{C_{\lambda_1}(t-\tau)}
\label{q15-1}\end{eqnarray}
and
\begin{eqnarray}
\int^t_{\Omega}\|u_1(s)-u_2(s)\|^2_{D(A^{\frac{\sigma+1}{2}})}ds\leq \hat{C}\Big[\|u^1_0-u^2_0\|^2_{ D(A^{\frac{\sigma}{2}})}+\Big(\| u_2\|^2_{L^2(\tau,T;V)}+\| u_1\|^2_{L^2(\tau,T;V)}\Big)\Big] e^{C_{\lambda_1}(t-\tau)},\nonumber
\label{q15-1-1}\end{eqnarray}	
this means the strong continuity can not be obtained, however, from the compact argument and existence of regular weak solutions, $u(t,x)$ and $v(t,x)$ are norm-to-weak continuous, which complete the proof. $\hfill$$\Box$

\begin{Lemma}\label{th7.7}
	(Norm-to-weak continuity of process) 	
	If $v_{\tau}\in D(A^{\frac{\sigma}{2}})$, $f\in L^2_{loc}(\mathbb{R};H)$, then we can define the process $U(\cdot,\cdot):\mathbb{R}^2_d\times D(A^{\frac{\sigma}{2}})\rightarrow D(A^{\frac{\sigma}{2}})$ as $U(t,\tau)v_{\tau}=v(t,\tau,v_{\tau})$.
	
	Moreover, the process $\{U(t,\tau)\}$ is pullback norm-to-weak continuous from $D(A^{\frac{\sigma}{2}})$ to $D(A^{\frac{\sigma}{2}})$.
\end{Lemma}
{\bf Proof.}
Using Theorem \ref{th7.6}, we can derive the pullback norm-to-weak continuity of process in $D(A^{\frac{\sigma}{2}})$ for the problem \eqref{c1a}.\\

\noindent $\bullet$ {\bf Pullback $\mathcal{D}'_{\mu}$-absorbing sets family in $D(A^{\frac{\sigma}{2}})$}
\begin{Lemma}
	Assume $f \in  L^2_{loc}(\mathbb{R};H)$ is pullback translation bounded in $L^2_{pb}(\mathbb{R};H)$ (or uniformly pullback tempered), $v_{\tau}\in  D(A^{\frac{\sigma}{2}})$, if we choose parameter $\tilde{\mu}\in (0,\lambda_1\nu]$ fixed, then the solution $v$ to the problem \eqref{c1} satisfies that for any $\tau\leq t$,
	\begin{eqnarray}
	|A^{\frac{\sigma}{2}}v(t)|^2
	&\leq& \Big(\Big[\frac{2}{\nu\lambda_1}\|\varphi\|^4_{L^{\infty}(\partial\Omega)}|\partial\Omega|+\frac{2C\nu|\partial\Omega|}{\varepsilon}\|\varphi\|^2_{L^{\infty}(\partial\Omega)}\Big]
	+|A^{\frac{\sigma}{2}}v_{\tau}|^2\Big)e^{-\tilde{\mu}(t-\tau)}\nonumber\\ &&+\Big(\frac{2C}{\nu\lambda^{\frac{2}{3}}}+4C\|\varphi\|^2_{L^{\infty}(\partial\Omega)}\Big)e^{-\tilde{\mu}\lambda_1(t-\tau)}\int^t_{\tau}\|v(s)\|^2ds\nonumber\\ &&+\frac{2C}{\nu}\|f\|^2_{L^2_{loc}(\mathbb{R};H)}+\frac{2Ce^{-\tilde{\mu}n_0h}}{\nu(1-e^{-\tilde{\mu}n_0h})}\|f\|^2_{L^2_{pb}(\mathbb{R};H)}.\label{gg5-1}\end{eqnarray}
	If we denote \begin{eqnarray}K^2_3&=&\Big(\frac{2C}{\nu\lambda^{\frac{2}{3}}}+4C\|\varphi\|^2_{L^{\infty}(\partial\Omega)}\Big)\|v\|^2_{L^2(\tau,t;V)}+\frac{2}{\nu\lambda_1}\|\varphi\|^4_{L^{\infty}(\partial\Omega)}|\partial\Omega|+\frac{2C\nu|\partial\Omega|}{\varepsilon}\|\varphi\|^2_{L^{\infty}(\partial\Omega)},\end{eqnarray}
	then
	\begin{eqnarray}
	\hspace{-1cm}|A^{\frac{\sigma}{2}}v(t)|^2
	&\leq& e^{-\tilde{\mu} t}\Big(K^2_3+|A^{\frac{\sigma}{2}}v_{\tau}|^2\Big)e^{\tilde{\mu}\tau}+\frac{2C}{\nu}\|f\|^2_{L^2_{loc}(\mathbb{R};H)}+\frac{2Ce^{-\tilde{\mu}n_0h}}{\nu(1-e^{-\tilde{\mu}n_0h})}\|f\|^2_{L^2_{pb}(\mathbb{R};H)}+1.\label{gg5-2}\end{eqnarray}
\end{Lemma}
{\bf Proof.} Taking inner product with \eqref{c1a} with $A^{\sigma}v$ in $H$, we obtain that
\begin{eqnarray}
&&\hspace{-0.5cm}(v_t,A^{\sigma}v)+\nu(A v,A^{\sigma}v)+(B(v,v),A^{\sigma}v)+(B(v,\psi),A^{\sigma}v)+(B(\psi,v),A^{\sigma}v)\nonumber\\
&&
=-(B(\psi,\psi),A^{\sigma}v)+\langle P f,A^{\sigma}v\rangle -\langle \nu P F, A^{\sigma}v\rangle,\label{gg11}
\end{eqnarray}
i.e.,
\begin{eqnarray}
&&\frac{1}{2}\frac{d}{dt}\|A^{\frac{\sigma}{2}}v\|^2_H
+\nu\|A^{\frac{\sigma+1}{2}}v\|^2_H\leq |(B(v,v),A^{\sigma}v)|+|(B(v,\psi),A^{\sigma}v)|+|(B(\psi,v),A^{\sigma}v)|\nonumber\\
&&
\hspace{5cm}+|(B(\psi,\psi),A^{\sigma}v)|+|\langle P f,A^{\sigma}v\rangle|+|\langle \nu P F, A^{\sigma}v\rangle|.\label{gg12}
\end{eqnarray}
Using the same techniques as in Lemma \ref{le6.4}, we have
\begin{eqnarray}
|(B(v,v),A^{\sigma}v)|&\leq& \frac{\nu}{12}|A^{\frac{\sigma+1}{2}}v|^2+\frac{C}{\nu\lambda_1^{\frac{2}{3}}}\|v\|^2,\label{gg13}\\
|(B(v,\psi),A^{\sigma}v)|&\leq&\frac{\nu}{12}|A^{\frac{\sigma+1}{2}}v|^2
+C\|\varphi\|^2_{L^{\infty}(\partial\Omega)}\|v\|^2,\label{gg14}\\
|(B(\psi,v),A^{\sigma}v)|&\leq& \frac{\nu}{12}|A^{\frac{\sigma+1}{2}}v|^2+C\|\varphi\|^2_{L^{\infty}(\partial\Omega)}\|v\|^2,\label{gg15}\\
|(B(\psi,\psi),A^{\sigma}v)|&\leq& \frac{\nu}{12}|A^{\frac{\sigma+1}{2}}v|^2
+\frac{1}{\nu\lambda_1}\|\varphi\|^2_{L^{\infty}(\partial\Omega)}\|\varphi\|^2_{L^2(\partial\Omega)}\nonumber\\
&\leq& \frac{\nu}{12}|A^{\frac{\sigma+1}{2}}v|^2
+\frac{1}{\nu\lambda_1}\|\varphi\|^4_{L^{\infty}(\partial\Omega)}|\partial\Omega|,\label{gg16}\\
|\langle P f,A^{\sigma}v\rangle|&\leq& \frac{C}{\nu}|f|^2+\frac{\nu}{12}|A^{\frac{\sigma+1}{2}}v|^2,\label{gg17}
\end{eqnarray}
\begin{eqnarray}
|\langle \nu P F, A^{\sigma}v\rangle|&\leq&\nu\int_{\Omega}|F||A^{\sigma}v|dx\leq \nu\Big(\int_{\Omega}|F|^2dx\Big)^{\frac{1}{2}}\Big(\int_{\Omega}|A^{\sigma}v|^2dx\Big)^{\frac{1}{2}}\nonumber\\
&\leq&\nu\|F\|_{L^2(\Omega)}|A^{\sigma}v|\leq \frac{C\nu}{\varepsilon^{\frac{3}{2}}}\|\varphi\|_{L^2(\partial\Omega)}|A^{\sigma}v|\nonumber\\
&\leq& \frac{C\nu}{\varepsilon^{\frac{1}{2}}}|\partial\Omega|^{\frac{1}{2}}\varepsilon\|\varphi\|_{L^{\infty}(\partial\Omega)}|A^{\frac{\sigma+1}{2}}v|\nonumber\\
&\leq& \frac{C\nu|\partial\Omega|}{\varepsilon}\|\varphi\|^2_{L^{\infty}(\partial\Omega)}+
\frac{\nu}{12}|A^{\frac{\sigma+1}{2}}v|^2,\label{gg18}
\end{eqnarray}
and thus combining \eqref{gg11}--\eqref{gg18}, we can conclude
\begin{eqnarray}
\frac{d}{dt}|A^{\frac{\sigma}{2}}v|^2
+\lambda_1\nu|A^{\frac{\sigma+1}{2}}v|^2
&\leq& \Big(\frac{2C}{\nu\lambda_1^{\frac{2}{3}}}+4C\|\varphi\|^2_{L^{\infty}(\partial\Omega)}\Big)\|v\|^2+\frac{2}{\nu\lambda_1}\|\varphi\|^4_{L^{\infty}(\partial\Omega)}|\partial\Omega|\nonumber\\
&&+\frac{2C\nu|\partial\Omega|}{\varepsilon}\|\varphi\|^2_{L^{\infty}(\partial\Omega)}+\frac{2C}{\nu}|f|^2\nonumber\\
&=& K^2_2+\frac{2C}{\nu}|f|^2,\label{gg19}
\end{eqnarray}
where $K^2_2=\Big(\frac{2C}{\nu\lambda^{\frac{2}{3}}}+4C\|\varphi\|^2_{L^{\infty}(\partial\Omega)}\Big)\|v\|^2+\frac{2}{\nu\lambda_1}\|\varphi\|^4_{L^{\infty}(\partial\Omega)}|\partial\Omega|
+\frac{2C\nu|\partial\Omega|}{\varepsilon}\|\varphi\|^2_{L^{\infty}(\partial\Omega)}$.

Using the Gronwall inequality in $[\tau,t]$ to \eqref{gg19}, we obtain
\begin{eqnarray}
|A^{\frac{\sigma}{2}}v|^2
&\leq& e^{-\lambda_1\nu t}|A^{\frac{\sigma}{2}}v_{\tau}|^2 e^{\lambda_1\nu\tau}+\int^t_{\tau}e^{-\lambda_1\nu(t-s)}(K^2_2+\frac{2C}{\nu}|f(s)|^2)ds.\label{gg20}
\end{eqnarray}
Using \eqref{gg20}, we can easily get the desired result. $\hfill$$\Box$\\

\begin{Lemma}\label{th7.9}
	Assume $f \in  L^2_{loc}(\mathbb{R};H)$ is pullback translation bounded in $L^2_{pb}(\mathbb{R};H)$ (or uniformly pullback tempered), for any small enough $\varepsilon'_1>0$, there exists a pullback time $\hat{\tau}(t,\varepsilon_1)$, such that for any $\tau< \hat{\tau}(t,\varepsilon_1)\leq t$, $\hat{B}'_0(t)$ is a family of pullback $\mathcal{D}'_{\tilde{\mu}}$-absorbing sets for the process $U(t,\tau)$.
\end{Lemma}
{\bf Proof.} Noting that \begin{eqnarray}
&&\hspace{-1cm}\|U(t,\tau)v_{\tau}\|^2_{D(A^{\frac{\sigma}{2}})}=|A^{\frac{\sigma}{2}}v|^2\nonumber\\
&\leq& \Big(K^2_3+|\rho'_{\hat{D}}(t)|^2\Big)e^{-\tilde{\mu}(t-\tau)}+\frac{2C}{\nu}\|f\|^2_{L^2_{loc}(\mathbb{R};H)}+\frac{2Ce^{-\tilde{\mu}n_0h}}{\nu(1-e^{-\tilde{\mu}n_0h})}\|f\|^2_{L^2_{pb}(\mathbb{R};H)}
\end{eqnarray}
and
there exists a pullback time $\hat{\tau}(t,\varepsilon_1)$, such that for any $\tau<\hat{\tau}(t,\varepsilon_1)\leq t$, it follows
\begin{eqnarray}
\hspace{-0.3cm}e^{-\tilde{\mu} t}\Big(K^2_3+|\rho'_{\hat{D}}(\tau)|^2\Big)e^{\tilde{\mu}\tau}\leq \varepsilon_1
\end{eqnarray}
since
\begin{equation}
e^{-\tilde{\mu} (t-s)}\leq e^{-\tilde{\mu}_0(t-s)}
\end{equation}
holds for any $\tilde{\mu}_0\leq \tilde{\mu}$.

Hence, we have \begin{equation}
|U(t,\tau)v_{\tau}|^2_{D(A^{\frac{\sigma}{2}})}\leq \varepsilon_1+(\rho'_0(t))^2-\frac{1}{2}\leq (\rho'_0(t))^2,
\end{equation}
which implies that $U(t,\tau)D(\tau)\subset B'_0(t)$, i.e., $\hat{B}'_0(t)$ is a family of pullback $\mathcal{D}'_{\tilde{\mu}}$-absorbing balls. \\

\noindent $\bullet$ {\bf Pullback $\mathcal{D}'_{\tilde{\mu}}$-asymptotic compactness for the pullback norm-to-weak continuous process in $D(A^{\frac{\sigma}{2}})$}\label{sub7.4}\\
In this section, since $D(A^{\frac{\sigma}{2}})$ is a uniformly convex space, the pullback $\mathcal{D}'_{\tilde{\mu}}$-condition-(MWZ) (see Definition \ref{de8.29} and Remark \ref{re8.30}) is equivalent to the pullback $\mathcal{D}'_{\tilde{\mu}}$-asymptotic compactness for the processes.

\begin{Lemma}\label{th7.10}
	Assume $f \in  L^2_{loc}(\mathbb{R};H)$ is pullback translation bounded in $L^2_{pb}(\mathbb{R};H)$ (or uniformly pullback tempered),  $v_{\tau}\in D(A^{\frac{\sigma}{2}})$, then the norm-to-weak continuous processes $U(t,\tau)$ satisfies pullback $\mathcal{D}'_{\tilde{\mu}}$-condition (MWZ) which is equivalent to pullback $\mathcal{D}'_{\tilde{\mu}}$-asymptotically compact in $D(A^{\frac{\sigma}{2}})$ for the system \eqref{c1a} which is equivalent to \eqref{a1}.
\end{Lemma}
{\bf Proof.} {\bf Step 1:} Let $\hat{B}'_0=\{B'_0(t)\}_{t\in\mathbb{R}}$ be the pullback $\mathcal{D}'_{\tilde{\mu}}$-absorbing family given in Theorem \ref{th7.9}, then there exists a pullback time $\tau'_{t,\varepsilon_1}$ such that
$|A^{\frac{\sigma}{2}}v|^2=\|U(t,\tau)v_{\tau}\|_{D(A^{\frac{\sigma}{2}})}^2\leq \rho'_{0}(t)$.

Since $V^{\sigma}=D(A^{\frac{\sigma}{2}})\subset H$ is a Hilbert space, the norm-to-weak continuous process is defined as $U(t,\tau): V^{\sigma}\rightarrow V^{\sigma}$, where $V^{\sigma}=V^{\sigma}_1\bigoplus V^{\sigma}_2$, $V^{\sigma}_1=span\{\omega'_1,\omega'_2,\cdots, \omega'_m\}$ and $V^{\sigma}_1\bot V^{\sigma}_2$.
$\bar{P}$ be a orthonormal projector from $V^{\sigma}$ to $V^{\sigma}_1$, then we have the decomposition
\begin{equation}
v=\bar{P}v+(I-\bar{P})v:=v_1+v_2,
\end{equation} for $v_1\in V^{\sigma}_1$, $v_2\in V^{\sigma}_2$.

{\bf Step 2:} From the existence of global solutions and pullback $\mathcal{D}'_{\tilde{\mu}}$-absorbing family of sets $\hat{B}'_0$, we know $\|A^{\frac{\sigma}{2}}v_1\|^2_H\leq \rho'_{0}(t)$.

{\bf Step 3:}  Next we only need to prove the $D(A^{\frac{\sigma}{2}})$-norm of $v_2$ is small enough as $\tau\rightarrow -\infty$.

Taking inner product of \eqref{c1a} with $A^{\sigma}v_2$, since $(A^{\frac{\sigma}{2}}v_1, A^{\frac{\sigma}{2}}v_2)=0$ and $(v_1,v_2)=0$, we have
\begin{eqnarray}&&\frac{1}{2}\frac{d}{dt}|A^{\frac{\sigma}{2}}v_2|^2+\nu|A^{\frac{\sigma+1}{2}}v_2|^2\nonumber\\
&\leq& |(B(v,v),A^{\sigma}v_2)|+|(B(v,\psi),A^{\sigma}v_2)|
+|(B(\psi,v),A^{\sigma}v_2)|\nonumber\\
&&+|(B(\psi),A^{\sigma}v_2)|+|<\bar{P} f,A^{\sigma}v_2>|+|<\nu \bar{P} F,A^{\sigma}v_2>|.\label{g2-1}\end{eqnarray}

Next, we shall estimate every term on the right-hand side of \eqref{g2-1}.

(a) By the property of $b(\cdot,\cdot,\cdot)$, using the $\varepsilon$-Young inequality, we have
\begin{eqnarray}
|(B(v,v),A^{\sigma}v_2)|&=&|b(v,v,A^{\sigma}v_2)|\leq\int_{\Omega}|v||\nabla v||A^{\sigma}v_2|dx\nonumber\\
&\leq&C |v|^{\frac{1}{2}}|A^{\sigma}v_2|^{\frac{1}{2}}|\nabla v|\leq \frac{C}{\lambda_1^{\frac{1}{4}}}|v|^{\frac{3}{2}}|A^{\sigma}v_2|^{\frac{1}{2}}\nonumber\\
&\leq&\frac{\nu}{12}|A^{\frac{\sigma+1}{2}}v_2|^{2}+\frac{C}{\nu\lambda_1^{\frac{2}{3}}}\|v\|^2_V.\label{g2-2}\end{eqnarray}

(b) Using the same technique in Lemma \ref{le6.4}, we get
\begin{eqnarray}
|(B(v,\psi),A^{\sigma}v_2)|\leq\int_{\Omega}|v||\nabla \psi||v_2|dx
&\leq& \frac{\nu}{12}|A^{\frac{\sigma+1}{2}}v_2|^{2}+\frac{C\|\varphi\|^2_{L^{\infty}(\partial\Omega)}}{\nu}\|v\|^2_V\label{g2-3}\end{eqnarray}
and
\begin{eqnarray}
|(B(\psi,v),A^{\sigma}v_2)|&\leq&\int_{\Omega}|\psi||\nabla v||A^{\sigma}v_2|dx\nonumber\\
&\leq& C_4\|\varphi\|_{L^{\infty}(\partial\Omega)}\int_{\Omega}|\nabla v||A^{\sigma}v_2|dx\nonumber\\
&\leq& \frac{\nu}{12}|A^{\frac{\sigma+1}{2}}v_2|^{2}+\frac{C\|\varphi\|^2_{L^{\infty}(\partial\Omega)} }{\nu}\|v\|^2_V.\label{g2-4}\end{eqnarray}

(c) Using the method in Lemma \ref{le6.4}, we obtain
\begin{eqnarray}
|(B(\psi,\psi),A^{\sigma}v_2)|&\leq&\int_{\Omega}|\psi||\nabla \psi||A^{\sigma}v_2|dx\leq \frac{\nu}{12}|A^{\frac{\sigma+1}{2}}v_2|^{2}+\frac{1}{\nu\lambda_1}\|\varphi\|^4_{L^{\infty}(\partial\Omega)}|\partial\Omega|.\label{g2-5}\end{eqnarray}

(d) Using the techniques in \eqref{gg17} and \eqref{gg18}, we can derive
\begin{eqnarray}
|<f,A^{\sigma}v_2>|&\leq& \frac{C}{\nu}|f|^2+\frac{\nu}{12}|A^{\frac{\sigma+1}{2}}v_2|^2\label{g2-6}\end{eqnarray}
and
\begin{eqnarray}
\nu|<F,A^{\sigma}v_2>|
&\leq& \frac{C\nu|\partial\Omega|}{\varepsilon}\|\varphi\|^2_{L^{\infty}(\partial\Omega)}+
\frac{\nu}{12}|A^{\frac{\sigma+1}{2}}v_2|^2.\label{g6-1}
\end{eqnarray}

Combining \eqref{g2-1}--\eqref{g6-1}, we conclude
\begin{eqnarray}\frac{d}{dt}|A^{\frac{\sigma}{2}}v_2|^2+\nu|A^{\frac{\sigma+1}{2}}v_2|^2
&\leq&\Big(\frac{C}{\nu\lambda_1^{\frac{2}{3}}}+\frac{4C\|\varphi\|^2_{L^{\infty}(\partial\Omega)} }{\nu}\Big)\|v\|^2_V+\frac{2}{\nu\lambda_1}\|\varphi\|^4_{L^{\infty}(\partial\Omega)}|\partial\Omega| \nonumber\\
&&
+\frac{2C\nu\|\varphi\|^2_{L^{\infty}(\partial\Omega)}|\partial\Omega|}{\varepsilon}+\frac{2C}{\nu}|f|^2.\label{6.41}\end{eqnarray}

By Poincar\'{e}'s inequality, using the definition of fractal power of operator $A$, since $\lambda_j$ is increasing with respect to $j=m+1,\cdots$, for $v_2=\displaystyle{\sum_{j=m+1}}a_j\omega'_j$ and $A^sv_2=\displaystyle{\sum_{j=m+1}}\lambda^sa_j\omega'_j$, we have
\begin{eqnarray}
|A^{\frac{\sigma+1}{2}}v_2|^2=\|v_2\|^2_{V^{\sigma}}
&=&\|\sum_{j=m+1}\lambda^{\frac{\sigma+1}{2}}_ja_j\omega'_j\|^2_H\geq \|\sum_{j=m+1}\lambda^{\frac{\sigma}{2}}_j\lambda^{\frac{1}{2}}_{m+1}a_j\omega'_j\|^2_H\nonumber\\
&\geq& \lambda_{m+1}\|\sum_{j=m+1}\lambda^{\frac{\sigma}{2}}_ja_j\omega'_j\|^2_H
=\lambda_{m+1}|A^{\frac{\sigma}{2}}v_2|^2.\label{6.42}\end{eqnarray}
Thus, we derive from \eqref{6.41}-\eqref{6.42}
\begin{eqnarray}\frac{d}{dt}\|A^{\frac{\sigma}{2}}v_2\|_H^2+\nu\lambda_{m+1}\|A^{\frac{\sigma}{2}}v_2\|^2_H
&\leq&\Big(\frac{C}{\nu\lambda_1^{\frac{2}{3}}}+\frac{4C\|\varphi\|^2_{L^{\infty}(\partial\Omega)} }{\nu}\Big)\|v\|^2_V+\frac{2}{\nu\lambda_1}\|\varphi\|^4_{L^{\infty}(\partial\Omega)}|\partial\Omega| \nonumber\\
&&
+\frac{2C\nu\|\varphi\|^2_{L^{\infty}(\partial\Omega)}|\partial\Omega|}{\varepsilon}+\frac{2C}{\nu}|f|^2\nonumber\\
&=&K^2_4\|v\|^2_V+K^2_5+\frac{2C}{\nu}|f|^2,\label{gg10-1}\end{eqnarray}
where $K_4^2=\Big(\frac{C}{\nu\lambda_1^{\frac{2}{3}}}+\frac{4C\|\varphi\|^2_{L^{\infty}(\partial\Omega)} }{\nu}\Big)$,
$K^2_5=\frac{2}{\nu\lambda_1}\|\varphi\|^4_{L^{\infty}(\partial\Omega)}|\partial\Omega|
+\frac{2C\nu\|\varphi\|^2_{L^{\infty}(\partial\Omega)}|\partial\Omega|}{\varepsilon}$.

Applying the Gronwall inequality in $[\tau,t]$ to \eqref{gg10-1}, we conclude
\begin{eqnarray}
|A^{\frac{\sigma}{2}}v_2(t)|^2&\leq& |A^{\frac{\sigma}{2}}v_2(\tau)|^2e^{-\nu\lambda_{m+1}(t-\tau)}+\int^t_{\tau}e^{-\nu\lambda_{m+1}(t-s)}|f(s)|^2ds\nonumber\\
&&+\int^t_{\tau}\Big(K^2_4\|v(s)\|^2+K^2_5\Big)e^{-\nu\lambda_{m+1}(t-s)}ds\nonumber\\
&\leq& e^{-\nu\lambda_{m+1}t}(\rho'_{\hat{D}}(t))^2e^{\nu\lambda_{m+1}\tau}+e^{-\nu\lambda_{m+1}t}e^{\nu(\lambda_{m+1}-1)\tau}\int^t_{-\infty}e^{\nu s}|f(s)|^2ds\nonumber\\
&&+\Big(K^2_4e^{-\nu\lambda_{m+1}t}e^{\nu\lambda_{m+1}\tau}\int^t_{\tau}\|v\|^2_Vds\Big)+(\frac{K^2_5}{\nu\lambda_{m+1}}-\frac{K^2_5e^{-\nu\lambda_{m+1}t}e^{\nu\lambda_{m+1}\tau}}{\nu\lambda_{m+1}})\nonumber\\
&=&I'_1+I'_2+I'_3+I'_4.\label{g7-1}
\end{eqnarray}

By the definition of the universe $\mathcal{D}'_{\tilde{\mu}}$ in \eqref{gg7-1-1}, using the technique in \eqref{gg7-1}, there exists a pullback time $\hat{\tau}'_1<<\tau_1$, then for any $\tau<\hat{\tau}'_1\leq t$ and $m$ large enough, it follows
\begin{eqnarray}
I'_1&=&\|A^{\frac{\sigma}{2}}v_2(\tau)\|_H^2e^{-\nu\lambda_{m+1}(t-\tau)}\leq e^{-\nu\lambda_{m+1}t}(\rho'_{\hat{D}}(t))^2e^{\nu\lambda_{m+1}\tau}\leq  \frac{\varepsilon}{4}\end{eqnarray}
as $\tau\rightarrow -\infty$ and $m\rightarrow +\infty$.

Since $f\in L^2_{loc}(\mathbb{R};H)$ and $K^2_5$ is bounded, using $\displaystyle{\lim_{m\rightarrow\infty}}\lambda_{m+1}=+\infty$ if we can choose $m$ large enough and letting $\tau\rightarrow -\infty$, i.e., there exists a pullback time $\hat{\tau}_2\leq t$ and positive integer $m_0>0$, for any $\tau<\hat{\tau}_2$ and $m>m_0$, we obtain
\begin{eqnarray}
I'_2=e^{-\nu\lambda_{m+1}t}e^{\nu\lambda_{m+1}\tau}\int^t_{\tau}|f(s)|^2d s\leq\frac{\varepsilon}{4}\end{eqnarray}
and
\begin{eqnarray}
I'_4=\frac{K^2_5}{\nu\lambda_{m+1}}-\frac{K^2_5e^{-\nu\lambda_{m+1}t}e^{\nu\lambda_{m+1}\tau}}{\nu\lambda_{m+1}}\leq \frac{\varepsilon}{4}.
\end{eqnarray}

Moreover, since $K_4^2$ and $\int^t_{\tau}\|v\|^2_Vds$ are bounded from Theorem \ref{th4.1}, using $\displaystyle{\lim_{m\rightarrow\infty}}\lambda_{m+1}=+\infty$ if we can choose $m$ large enough, letting $\tau\rightarrow -\infty$, i.e., there exists a pullback time $\hat{\tau}_3$, for any $\tau\leq\hat{\tau}_3$, we can achieve
\begin{eqnarray}
I'_3=K^2_4e^{-\nu\lambda_{m+1}t}e^{\nu\lambda_{m+1}\tau}\int^t_{\tau}\|v\|^2_Vds< \frac{\varepsilon}{4}.\label{g8-1}\end{eqnarray}

Combining  \eqref{g7-1}--\eqref{g8-1} and choosing a pullback time $$\hat{\tau}'\leq min\{\hat{\tau}_1,\hat{\tau}_2,\hat{\tau}_3,\hat{\tau}'_1,\hat{\tau}(t,\varepsilon_1)\},$$ such that for any $\tau<\hat{\tau}\leq t$, we conclude
\begin{eqnarray}
&&\|(I-P)U(t,\tau)v_{\tau}\|^2_{V^{\sigma}}=\|A^{\frac{\sigma}{2}}v_2\|^2_H\leq \varepsilon\end{eqnarray}
which implies that the norm-to-weak process satisfies pullback $\mathcal{D}'_{\tilde{\mu}}$-condition-(MWZ), i.e., $U(t,\tau)$ is pullback $\mathcal{D}'_{\tilde{\mu}}$-asymptotically compact. Hence, we complete the proof.$\hfill$$\Box$\\

\noindent $\bullet$ {\bf Proof of Theorem \ref{th7.11}: Minimal family of pullback  $\mathcal{D}'_{\tilde{\mu}}$-attractors in $D(A^{\frac{\sigma}{2}})$} Similarly to the proof in Theorem \ref{th5.6}, it is easily to prove the theorem here. $\hfill$$\Box$\\

\noindent $\bullet$ {\bf Proof of Theorem \ref{th7.12}: Unique family of pullback  $\mathcal{D}'_{\tilde{\mu}}$-attractors in $D(A^{\frac{\sigma}{2}})$} Using the same techniques as in the proof of Theorem \ref{th5.7}, we can derive the inclusion closed of the universe $\mathcal{D}'_{\tilde{\mu}}$ in $D(A^{\frac{\sigma}{2}})$ easily guarantees that $\hat{B}'_0\in\mathcal{D}'_{\tilde{\mu}}$, hence
using Theorem \ref{th8.35}, since the process is norm-to-weak continuous and closed, it is easy to achieve the uniqueness of family of pullback $\mathcal{D}'_{\tilde{\mu}}$-absorbing attractors $\mathcal{A}'_{\mathcal{D}'_{\tilde{\mu}}}$ in $D(A^{\frac{\sigma}{2}})$. $\hfill$$\Box$

\section{Appendix}
In Chepyzhov and Vishik \cite{cv3}, they used translation compact external forces which is translation bounded to construct uniform attractors, and Lu \cite{lu} extended this to normal function and translation bounded. Cui, Langa and Li \cite{cll} considered the backward property for the external forces for pullback attractors. We can consider pullback property of external forces in more general functional spaces.
 Let $X$ be a Banach spaces or Hilbert spaces, $L^p(\Omega)$ be the general $p$-power integral space.

 \begin{Definition}
		We call the external force $g(t,x)$ is pullback translation bounded if it satisfies
		\begin{eqnarray}\label{13}
		\sup_{s\leq \tau,\ h\leq 0}\int^s_{s+h}\|g(r,x)\|^p_{X}dr<+\infty
		\ \ \mbox{or}\ \
		\sup_{s\leq \tau}\int^{s}_{s-1}\|g(r,x)\|^p_{X}dr<+\infty
		\end{eqnarray}
		for  all $t \in \mathbb{R}$ and $h\leq 0$. Denote all pullback translation bounded in $L^p_{loc}(\mathbb{R};X)$ as $L^p_{pb}(\mathbb{R};X)$.
	\end{Definition}

 \begin{Definition}
 A function $g\in \mathcal{L}^p_b(\mathbb{R};X)$ is called pullback normal if
 \begin{eqnarray}\label{13}
		\sup_{s\leq \tau,\ h\leq 0}\int^s_{s+h}\|g(r,x)\|^p_{X}dr<\varepsilon
		\end{eqnarray}
holds for arbitrary $\varepsilon>0$ and $1<p<+\infty$. Denote all pullback normal functions in $L^p_{loc}(\mathbb{R};X)$ as $L^p_{pn}(\mathbb{R};X)$.
 \end{Definition}

\begin{Definition}
			The function $g(t,x)$ is called as pullback tempered in $L^p_{loc}(\mathbb{R};X)$, if
			\begin{eqnarray}
			\sup_{s\leq \tau}\int^s_{-\infty}e^{\alpha(r-s)}\|g(r,x)\|^p_{X}dr<+\infty\ \ \mbox{or}\ \ \sup_{s\leq \tau}\int^s_{-\infty}e^{\alpha(r-s)}\|g(r,x)\|^p_{X}dr\rightarrow 0\ \mbox{as}\ \tau\rightarrow-\infty.
				\end{eqnarray}	
Denote all pullback tempered functions in $L^p_{loc}(\mathbb{R};X)$ as $L^p_{pt}(\mathbb{R};X)$.				
\end{Definition}

\begin{Definition}
The function $g(t,x)$ is called as uniformly pullback tempered in $L^p_{loc}(\mathbb{R};X)$, if
			\begin{eqnarray}
			\int^{\tau}_{-\infty}e^{\alpha(r-\tau)}\|g(r,x)\|^p_{X}dr<+\infty\ \ \mbox{or}\ \ \int^t_{-\infty}e^{\alpha(r-t)}\|g(r,x)\|^p_{X}dr<+\infty.
				\end{eqnarray}
Denote all uniformly pullback tempered functions in $L^p_{loc}(\mathbb{R};X)$ as $L^p_{upt}(\mathbb{R};X)$.	
\end{Definition}

\begin{Theorem}
We have the following property:\\
(I) The pullback translation bounded is equivalent with pullback tempered;\\
(II) The pullback normal is stronger than pullback bounded;\\
(III) The uniformly pullback tempered is stronger than pullback translation bounded and pullback tempered.
\end{Theorem}
{\bf Proof.} (I) Assume that $g(t,x)$ is pullback translation bounded, then for $h<0$, we have
\begin{eqnarray}
\int^s_{-\infty}e^{\alpha(r-s)}\|g(r,x)\|^p_{X}dr
&=&\sum^{\infty}_{m=0}\int^{s+n h}_{s+(n+1)h}e^{\alpha(r-s)}\|g(r,x)\|^p_{X}dr\nonumber\\
&\leq& \sum^{\infty}_{n=0}\int^{s+n h}_{s+(n+1)h}e^{-\alpha n h}\|g(r,x)\|^p_{X}dr\nonumber\\
&\leq& \sum^{\infty}_{n=0}e^{-\alpha n h}\Big[\int^{s}_{s+h}\|g(r,x)\|^p_{X}dr\Big]<+\infty
\end{eqnarray}
which means pullback tempered.

Let $g(t,x)$ is pullback tempered, i.e.,
\begin{equation}
\sup_{s\leq \tau}\int^s_{-\infty}e^{\alpha(r-s)}\|g(r,x)\|^p_{X}dr\rightarrow 0\ \mbox{as}\ \tau\rightarrow-\infty,
\end{equation}
then it follows
\begin{eqnarray}
\int^s_{s+h}\|g(r,x)\|^p_Xdr &\leq& e^{-\alpha h}\int^s_{s+h}e^{\alpha(r-s)}\|g(r,x)\|^p_Xdr\nonumber\\
&\leq& e^{-\alpha h}\int^s_{-\infty}e^{\alpha(r-s)}\|g(r,x)\|^p_Xdr
\end{eqnarray}
which implies pullback translation bounded.

(II) and (III): Comparing the definition, it is easy to check the results. $\hfill$$\Box$\\

{\bf Acknowledgements:} This work was initiated when Xin-Guang Yang
was a long term visitor at ICMC-USP, Brazil, from May 2015 to June 2016, supported by FAPESP
(Grant 2014/17080-0). He was also partially supported by NSFC of China (Grant No.11726626), Yongjin Lu was partially supported by NSF (Grant No:
1601127), they were also supported by the Key Project
of Science and Technology of Henan Province (Grant No. 182102410069).
 Yuming Qin was in part supported by the NSFC of China (Grant No. 11671075). T. F. Ma was partially supported by CNPq (Grant No. 310041/2015).

\end{document}